\documentclass[12pt,a4paper]{article}
\usepackage[numbers,comma,sort&compress]{natbib}
\usepackage{amssymb}
\usepackage{lineno}
\usepackage{amsmath}
\usepackage{amssymb}
\usepackage{etoolbox}
\usepackage{graphicx}
\usepackage{dsfont}
\usepackage{txfonts}
\usepackage{mathrsfs}
\usepackage{caption}
\usepackage{subcaption}
\usepackage{enumitem}

\usepackage{amsthm}
\usepackage{cases}
\usepackage{siunitx}
\usepackage{upgreek}
\theoremstyle{remark}
\newtheorem*{remark}{Remark}
\newcommand{\toVect}[1]{{\boldsymbol{#1}}}
\newcommand{\toMat}[1]{{\mathbf{#1}}}
\newcommand{\dd}{\text{\,}\mathrm{d}}
\newcommand{\ds}{\mathrm{s}}
\newcommand{\x}{\toVect{x}}
\newcommand{\Displacement}{{\toVect{u}}}
\newcommand{\displacement}{{\toVect{v}}}
\newcommand{\Traction}{\toVect{t}}
\newcommand{\HighOrderTraction}{\toVect{r}}
\newcommand{\HighOrderJump}{\toVect{j}}
\newcommand{\displacementKnown}{\bar{u}}
\newcommand{\DisplacementKnown}{\toVect{\bar{u}}}
\newcommand{\partialndisplacementKnown}{\overline{v}}
\newcommand{\partialnDisplacementKnown}{\overline{\toVect{v}}}

\newcommand{\TractionKnown}{\toVect{\tractionKnown}}
\newcommand{\tractionKnown}{\bar{t}}
\newcommand{\HighOrderTractionKnown}{\toVect{\highOrderTractionKnown}}
\newcommand{\highOrderTractionKnown}{\bar{r}}
\newcommand{\HighOrderJumpKnown}{\toVect{\highOrderJumpKnown}}
\newcommand{\highOrderJumpKnown}{\bar{j}}
\newcommand{\phiKnown}{\bar{\phi}}
\newcommand{\wKnown}{\bar{w}}
\newcommand{\E}{\toVect{E}}
\newcommand{\n}{\toVect{n}}
\newcommand{\m}{\toVect{m}}
\newcommand{\tv}{\toVect{s}}

\newcommand{\Strain}{\toVect{\strain}}
\newcommand{\strain}{\varepsilon}
\newcommand{\CauchyStress}{\toVect{\cauchyStress}}
\newcommand{\cauchyStress}{\hat{\sigma}}
\newcommand{\HyperStress}{\toVect{\hyperStress}}
\newcommand{\hyperStress}{\tilde{\sigma}}
\newcommand{\ElectricDisp}{\toVect{\electricDisp}}
\newcommand{\electricDisp}{\hat{D}}
\newcommand{\ShapeOperator}{\toMat{\shapeOperator}}
\newcommand{\shapeOperator}{S}
\newcommand{\Projector}{\toMat{\projector}}
\newcommand{\projector}{P}
\newcommand{\curvatureProjector}{\tilde{N}}
\newcommand{\CurvatureProjector}{\toMat{\curvatureProjector}}
\newcommand{\Piezo}{\toVect{\piezo}}
\newcommand{\piezo}{e}
\newcommand{\Flexo}{\toVect{\flexo}}
\newcommand{\flexo}{\mu}
\newcommand{\Dielec}{\toVect{\dielec}}
\newcommand{\dielec}{\kappa}
\newcommand{\Elast}{\varmathbb{C}}
\newcommand{\elast}{\varmathbb{C}}
\newcommand{\StrGr}{\toMat{\strGr}}
\newcommand{\strGr}{h}
\newcommand{\divergence}{\nabla\text{\!}\cdot\text{\!}}
\newcommand{\gradient}{\nabla}
\newcommand{\symmetricGradient}{\nabla^\text{sym}}
\newcommand{\surfaceGradient}{\nabla^{S}}

\newcommand{\surfaceDivergence}{\nabla^S\text{\!}\text{\!}\cdot\text{\!}}

\newcommand{\jump}[1]{\left\llbracket#1\right\rrbracket}
\newcommand{\trace}[1]{Tr(\,#1\,)}

\newcommand{\var}[2][]{\ifstrempty{#1}
                    {{\delta #2}}
                    {{\delta\hspace{-0.1em}_{#1} #2}}}
\newcommand{\vvar}[2][]{\ifstrempty{#1}
                    {{\delta^2 #2}}
                    {{\delta^2\hspace{-0.1em}_{#1} #2}}}
\newcommand{\normLL}[2][]{\ifstrempty{#1}
                    {{\lVert #2 \rVert}_{L_2(\Omega)}}
                    {{\lVert #2 \rVert}_{L_2(#1)}}}

\newcommand{\eq}{Eq.~}
\newcommand{\fig}{Fig.~}
\newcommand{\etal}{et.~al.~}
\newcommand{\ie}{i.e.~}
\newcommand{\eg}{e.g.~}
\newcommand{\um}{~\si{\um}}
\newenvironment{rcases}{\left.\begin{aligned}}{\qquad\end{aligned}\right\rbrace}

\title{An immersed boundary hierarchical B-spline method for flexoelectricity}
\author
{D. Codony$^1$, O. Marco$^1$, S. Fern\'{a}ndez-M\'{e}ndez$^1$, I. Arias$^{1,\ast}$\\
	\\
	\small{
	$^1$ Laboratori de C\`{a}lcul Num\`{e}ric (LaC\`{a}N), Universitat Polit\`{e}cnica de Catalunya (UPC),\vspace{-.4em}}
	\\
	\small{Campus Nord UPC-C2, E-08034 Barcelona, Spain}
	\\
	\small{$^\ast$ Corresponding author; E-mail:  irene.arias@upc.edu.}
}
\date{}

\begin{document}
\maketitle
\newcommand{\sep}{,~}
\begin{abstract}
This paper develops a computational framework with unfitted meshes to solve linear piezoelectricity and flexoelectricity electromechanical boundary value problems including strain gradient elasticity at infinitesimal strains.
The high-order nature of the coupled PDE system is addressed by a sufficiently smooth hierarchical B-spline approximation on a background Cartesian mesh.
The domain of interest is embedded into the background mesh and discretized in an unfitted fashion. 
The immersed boundary approach allows us to use B-splines on arbitrary domain shapes, regardless of their geometrical complexity, and could be directly extended, for instance, to shape and topology optimization.
The domain boundary is represented by NURBS, and exactly integrated by means of the NEFEM mapping.
Local adaptivity is achieved by hierarchical refinement of B-spline basis, which are efficiently evaluated and integrated thanks to their piecewise polynomial definition.
Nitsche's formulation is derived to weakly enforce essential boundary conditions, accounting also for the non-local conditions on the non-smooth portions of the domain boundary (\ie edges in 3D or corners in 2D) arising from Mindlin's strain gradient elasticity theory.
Boundary conditions modeling sensing electrodes are formulated and enforced following the same approach.
Optimal error convergence rates are reported using high-order B-spline approximations.
The method is verified against available analytical solutions and well-known benchmarks from the literature.
\\


\end{abstract}

\section{Introduction}\label{sec_01}
Electroactive materials are able to transform mechanical energy into electrical energy (and viceversa), which can be used for sensing, actuating or energy harvesting applications. A wide range of modern technologies are based on the electromechanical properties of these materials, such as cameras, printers or motors.

Different electromechanical couplings can be found depending on the material.
The most common coupling is \emph{piezoelectricity}, by which the strain $\Strain$ and polarization $\toVect{p}$ are linearly coupled:
\begin{equation}
p_l=d_{lij}\strain_{ij},
\end{equation}
where $\mathbf{d}$ is the third-rank tensor of piezoelectricity.
This is the case of piezoelectric ceramics, which are polarized by deformation, and conversely deform when an electrical field is applied. Some piezoelectrics exhibit further electromechanical couplings, such as \emph{pyroelectricity} (temperature-dependent polarization) or \emph{ferroelectricity} (reversible spontaneous polarization). Soft materials such as piezoelectric polymers or dielectric elastomers exhibit also \emph{electrostriction}, a nonlinear electromechanical coupling between the strain state and the square of the polarization field.

This variety of electromechanical couplings has been largely studied, is quite well understood and is suitable to model electromechanical couplings in materials at a macroscale. However, micro- and nanoscale electromechanics cannot be described by just considering traditional models,
because
additional effects become relevant 
at small scales, prominently \emph{flexoelectricity}.

Flexoelectricity is a two-way linear coupling between electric polarization and \emph{strain gradient}.
The (direct) flexoelectric effect is understood as the material polarization due to inhomogeneous deformation (\eg bending) and is mathematically expressed as
\begin{equation}
p_l=f_{lijk}\frac{\partial\strain_{ij}}{\partial x_k},
\end{equation}
where $\mathbf{f}$ is the fourth-rank tensor of flexoelectricity. There also exists a thermodynamically conjugate converse flexoelectric effect that consists on the generation of stress $\toVect{\sigma}$ due to the application of an inhomogeneous electric field $\toVect{E}$, \ie
\begin{equation}
\sigma_{ij}=f_{lijk}\frac{\partial E_{l}}{\partial x_k}.
\end{equation}

\begin{figure}[b!]\centering
	\begin{subfigure}[c]{0.5\textwidth}\centering
		\captionof{figure}{}\vspace{-0.5em}
		\includegraphics[width=\textwidth]{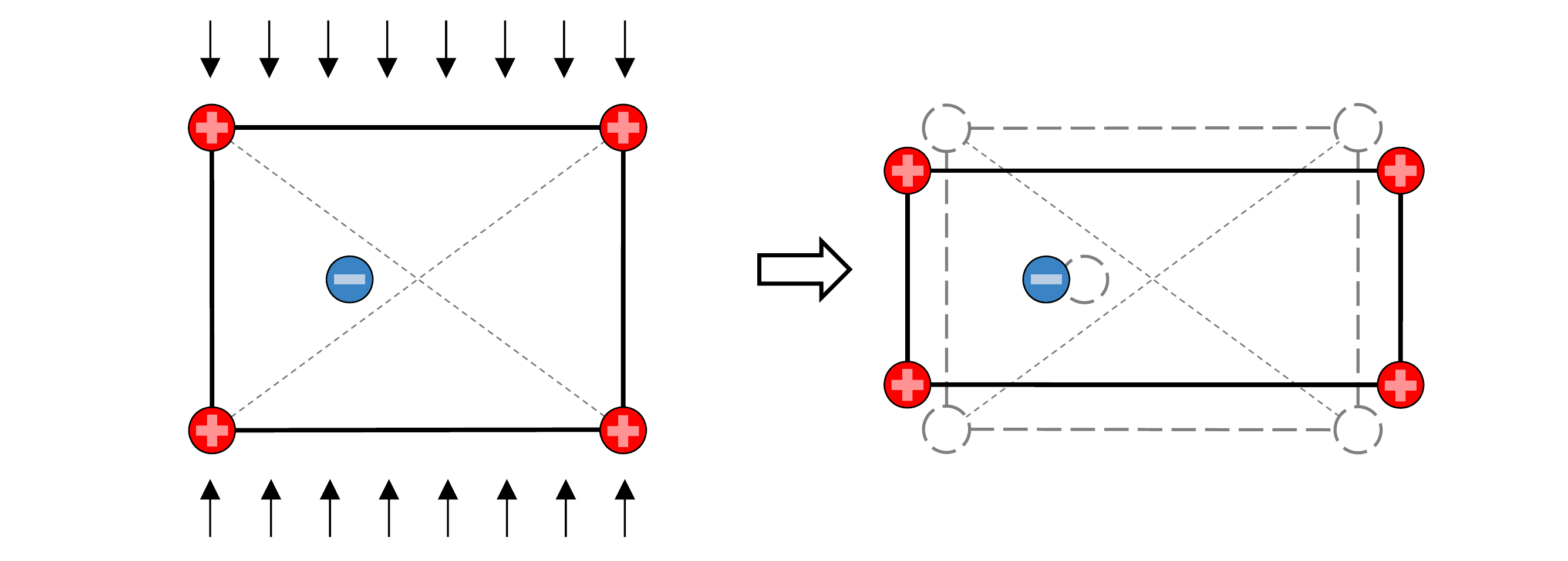}
		\label{fig_sketch1}\vspace{-1.5em}
		\newline
		\captionof{figure}{}\vspace{-0.5em}
		\includegraphics[width=\textwidth]{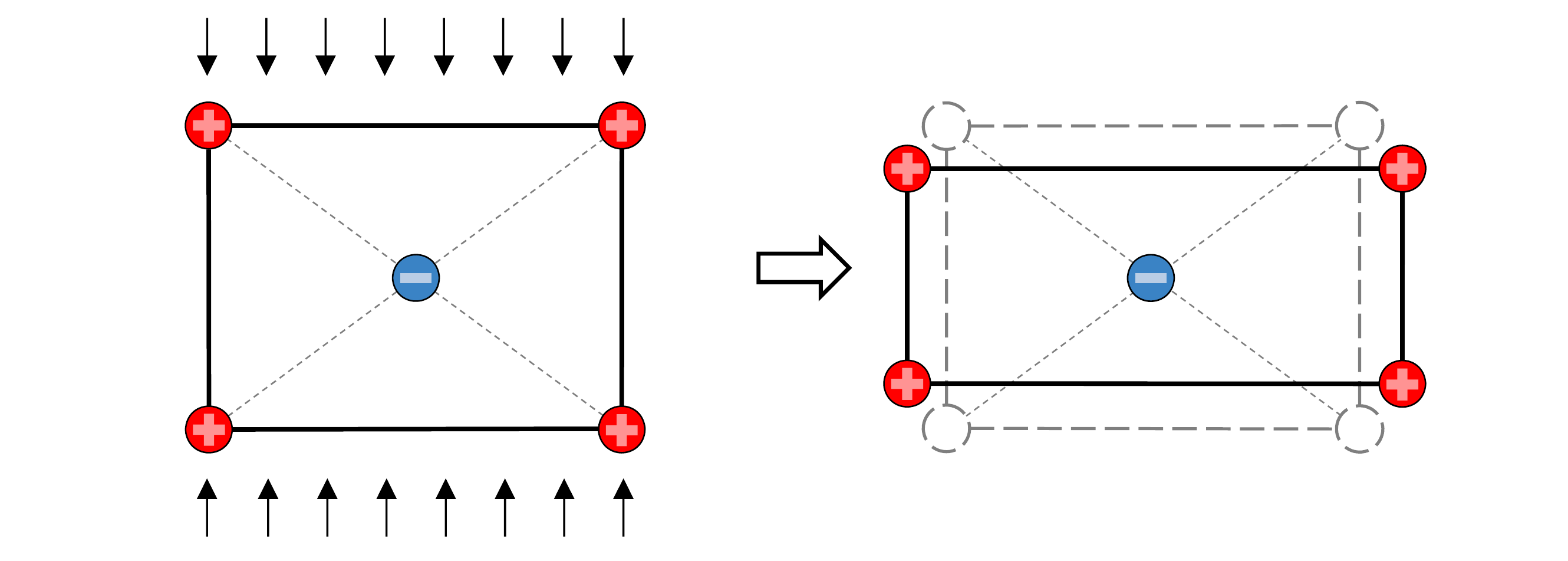}
		\label{fig_sketch2}\vspace{-1.5em}
	\end{subfigure}%
	\begin{subfigure}[c]{0.5\textwidth}\centering
		\captionof{figure}{}\vspace{-0.5em}
		\includegraphics[width=\textwidth]{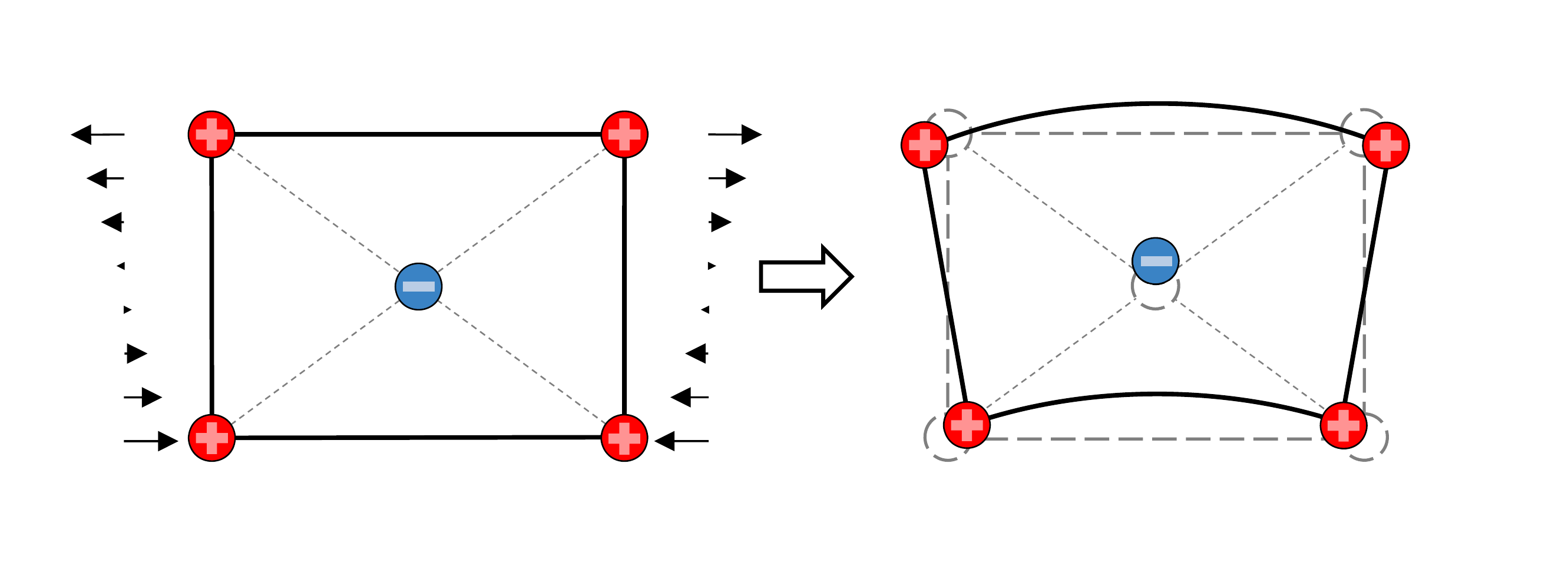}
		\label{fig_sketch3}\vspace{-1.5em}
	\end{subfigure}
	\caption{2D sketches of a crystalline material unit cell. \emph{a)} Compression-induced polarization in a non-centrosymmetric crystal. \emph{b)} Centrosymmetric crystal does not polarize under uniform compression. \emph{c)} A non-uniform deformation (\ie bending) breaks the inversion symmetry of the material and induces a net polarization, regardless of its crystalline structure.}
	\label{fig_sketch}
\end{figure}

Compared to piezoelectricity, flexoelectricity has two distinctive features.
On the one hand, it is universal, meaning it is present in any dielectric material. For a crystalline material to be piezoelectric, its crystalline structure is required to be non-centrosymmetric in order to allow for a net polarization as a result of a uniform deformation (\fig\ref{fig_sketch1}). Otherwise, the relative position of positive and negative ions remains unchanged after deformation and no net polarization is expected (\fig\ref{fig_sketch2}).
However, flexoelectricity generically breaks the inversion symmetry of the material, regardless of the internal crystalline structure, and a net polarization is observed after a non-uniform mechanical stimulus such as bending (\fig\ref{fig_sketch3}).
On the other hand, the flexoelectric material constants are typically small, and therefore sufficiently large strain gradients are required in order to trigger a sizable flexoelectric effect. Since strain-gradients scale inversely to spatial dimension, they are considerably large in the micro- and nanoscale. Therefore, flexoelectricity is by nature a size dependent effect.

Flexoelectricity in crystalline dielectrics was first studied by Mashkevich \cite{Mashkevich1957}, Tolpygo \cite{Tolpygo1963} and Kogan \cite{Kogan1964}, who proposed the first phenomenological model. In 1968, Bursian \etal\cite{bursian1968} performed the first experiment showing evidence of flexoelectricity in ferroelectric films and, in fact, it is not until 1981 that the phenomenon is named flexoelectricity \cite{Indenbom1981}. 
The first comprehensive theoretical works by Tagantsev \cite{Tagantsev1986,Tagantsev1991} clarified the distinction between piezoelectricity and flexoelectricity. However, since its effect is negligible at the macro-scale, it received little attention. In recent years, there has been a renewed interest in the scientific community, motivated by the need to downscale electromechanical transduction, and enabled by recent developments in nanotechnology \cite{Nguyen2013}.
Maranganti proposed the first mathematical framework for the flexoelectric governing equations \cite{Maranganti2006}. Following this work, flexoelectricity has been studied analytically for simple reduced models under restrictive assumptions, such as cantilever beams \cite{Majdoub2008} and thin films \cite{Sharma2010}, to name a few. In this last work by Sharma, atomistic calculations are also performed in order to verify the analytical results.
Cross developed a formulation to measure experimentally the longitudinal flexoelectric effect \cite{Cross2006}.
Other authors consider further physics, such as the flexoelectric effect in ferroelectrics \cite{Catalan2004,Eliseev2009}, the coupling with magnetic fields \cite{Liu2014} and the contributions of surface effects \cite{Shen2010}.
The general variational principles for flexoelectric materials can be found in \cite{Hu2010,Shen2010,Liu2014}.
The reader is referred to \cite{Yudin2013,Nguyen2013,Zubko2013,krichen2016} for recent reviews of flexoelectricity in solids.

Within the continuum flexoelectric theory, the symmetry of the flexoelectric tensor is well understood \cite{Shu2011,LeQuang2011}, although its full characterization is still lacking for most materials \cite{Zubko2013}. The equations are a coupled system of 4th-order partial differential equations, which renders analytical solutions difficult to obtain and precludes the use of conventional $C^0$ finite elements.
Several numerical alternatives have been proposed in the literature, based on smooth approximations with at least $C^1$ continuity \cite{Abdollahi2014,Abdollahi2015a,Abdollahi2015b,Ghasemi2017,Nanthakumar2017,Ghasemi2018,Yvonnet2017} or on mixed formulations \cite{Mao2016,Deng2017}.
The first self-consistent numerical solution of the linear flexoelectric problem was provided by Abdollahi \etal \cite{Abdollahi2014,Abdollahi2015a,Abdollahi2015b} using a mesh-free approach in 3D.
The degrees of freedom correspond only to displacements and electric potential, discretized with a $C^\infty$-continuous approximation to address the high-order nature of the equations.
This method was successfully applied to study the effect of flexoelectricity on the fracture of piezoelectric materials \cite{Abdollahi2015b}, and on the design of bimorph microsensors and microactuators \cite{Abdollahi2015c}.
Later, an alternative 2D continuum approach was proposed by the group of Aravas \etal \cite{Mao2016}, extending the mixed FEM formulation originally developed in \cite{Aravas2011} for strain-gradient elasticity to flexoelectricity. 
Displacement and displacement gradient fields are treated as separate degrees of freedom in order to circumvent the $C^1$-continuity requirement. This approach was also used by Deng \cite{Deng2017}.
Another alternative is the isogeometric approach, which has been used to perform topology optimization on 2D flexoelectric cantilever beams \cite{Ghasemi2017,Nanthakumar2017,Ghasemi2018}.
More recently, the $C^1$ triangular Argyris element was used by Yvonnet \etal in \cite{Yvonnet2017} to model flexoelectricity in soft dielectrics at finite strains.

In this paper, we propose an immersed boundary hierarchical B-spline approach to numerically solve the governing equations of flexoelectricity in 2D and 3D. This method enables simulations on arbitrary geometries within a reasonable computational cost, unlike previous works in the literature. For the sake of simplicity we restrict ourselves to infinitesimal strains, although the same idea applies also to finite strains \cite{CodonyNonlinearFlexoUnpublished}.

In this approach, the domain boundary is immersed into a fixed Cartesian mesh, and a hierarchical B-spline basis is built on top of it to discretize the primal unknowns (\ie displacements and electric potential), fulfilling the smoothness requirement of the equations.
%Due to the embedded boundary, the discretization step does not depend on the physical domain shape, which reports two main advantages.
The computational mesh does not fit to the embedded boundary, overcoming the rigidity of IGA approaches \cite{Ghasemi2017} that require Cartesian-like body-fitted meshes, difficult to generate for non-trivial geometries. In our case, mesh generation is straightforward regardless of the complexity of the domain shape. Moreover, a fixed mesh facilitates shape and topology optimization, avoiding re-meshing and the projection of the solution at each iteration.
In this work, the domain boundary is represented explicitly by NURBS surfaces in 3D and NURBS curves in 2D, which can be exactly integrated by means of the NEFEM mapping \cite{Sevilla2008}, but any other geometrical description, such as \eg level sets \cite{Sukumar2001,Belytschko2003} or subdivision surfaces \cite{Bandara2016}, could also be considered.

Local mesh refinement to resolve local features can be implemented in a B-spline context with several approaches, such as T-Splines \cite{Scott2011} and hierarchical B-splines (HB-splines) \cite{Forsey1988,Kraft1995,Kraft1997}. In this work we consider the latter, mainly due to its straightforward generalization to arbitrary dimensions and its relatively simple implementation.

For the first time to our knowledge, the complete set of boundary conditions is explicitly considered in a numerical solution of the flexoelectric boundary value problem.
In the seminal Mindlin's theory of strain gradient elasticity \cite{Mindlin1964,Mindlin1968a,Mindlin1968b}, which is the basis for deriving a stable flexoelectric theory \cite{Mao2014,Abdollahi2014, Abdollahi2015a, Abdollahi2015b}, additional non-local boundary conditions are required along non-smooth regions of the domain boundary (\ie corners in 2D and edges in 3D).
This is also the case for the flexoelectric theory.
However, in practice, the numerical methods mentioned above neglect these non-local conditions or consider smooth enough domains so that they do not appear \cite{Mao2016,Abdollahi2014,Abdollahi2015a,Abdollahi2015b,Ghasemi2017,Deng2017,Yvonnet2017}.
In this work, we show that non-local boundary conditions are mathematically required and we consider them in the formulation and implementation.
We demonstrate that neglecting them can deteriorate the solution.
In addition, we formulate the boundary conditions corresponding to sensing electrodes, which are common in electromechanical setups.

Within the unfitted framework, a Nitsche's formulation is derived for the flexoelectric equations to weakly enforce essential boundary conditions, accounting also for the non-local condition from Mindlin's theory. We show that not only the normal to the boundary, but also the curvatures play a role in the correct enforcement of boundary conditions. A Nitsche's formulation is also proposed to enforce electrode boundary conditions.

The method converges optimally for high-order approximations of degree $p$ in the $L_2$ norm and $H_s$ semi-norms, for $s=1,\dots,p$. Namely, it achieves the optimal convergence rates $p+1-s$.

The paper is organized as follows.
The variational formulation for flexoelectricity and the associated boundary value problem are presented in Section \ref{sec_02}.
The numerical approximation based on B-spline approximation and the immersed boundary method are presented in Section \ref{sec_04}.
Some illustrative numerical examples are given in Section \ref{sec_05}. The importance of considering non-local boundary conditions is illustrated in the first example. In the second one, we perform a sensitivity analysis with respect to the Nitsche penalty parameters. In the third one, optimal convergence is tested with a synthetic problem. The remaining examples show 2D and 3D simulations, and compare with available analytical solutions and well-known benchmarks from the literature.
%The paper is concluded and summarized in section \ref{sec_06}.
\section{Variational formulation and associated boundary value problem}\label{sec_02}
%%%%%
\subsection{Notation and preliminary definitions}
\begin{figure}[b]\centering
	\begin{subfigure}[t]{.34\textwidth}\centering
		\captionof{figure}{}
		\includegraphics[scale=0.17]{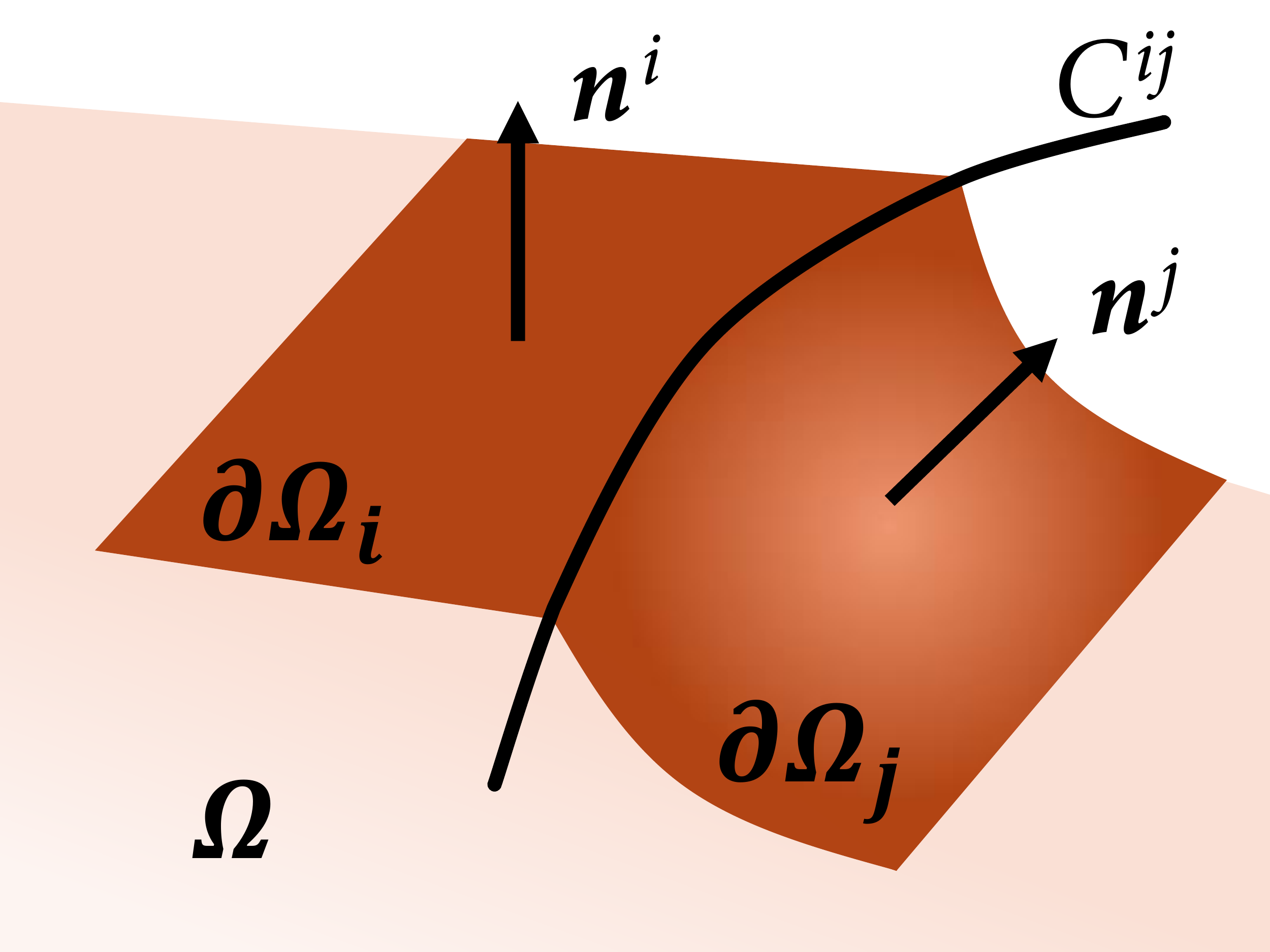}
		\vfill\label{fig_2_1_2}
	\end{subfigure}%
	\begin{subfigure}[t]{.25\textwidth}\centering
		\captionof{figure}{}
		\includegraphics[scale=0.17]{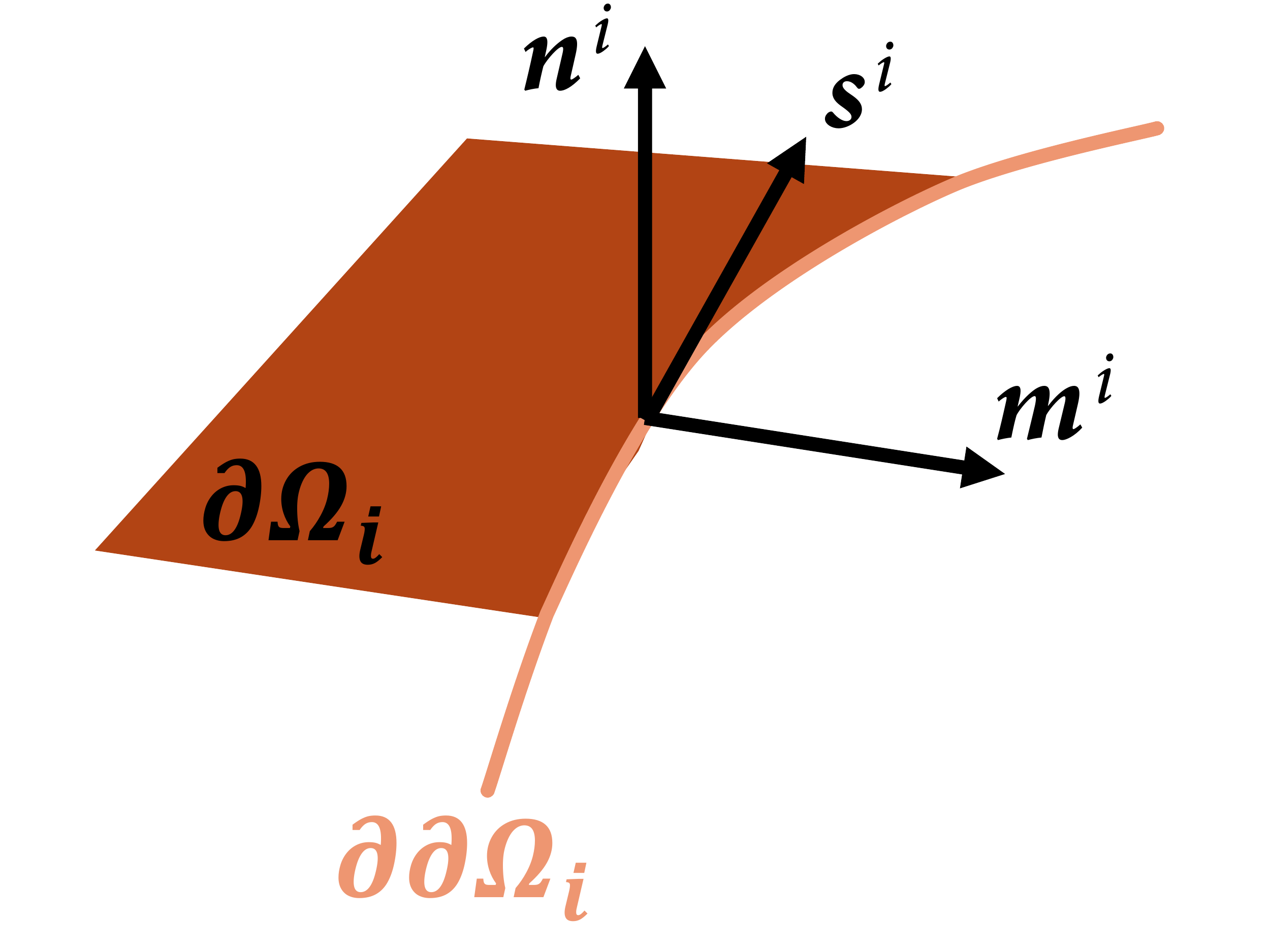}
		\vfill\label{fig_2_1_3}
	\end{subfigure}%
	\begin{subfigure}[t]{.24\textwidth}\centering
		\captionof{figure}{}
		\includegraphics[scale=0.17]{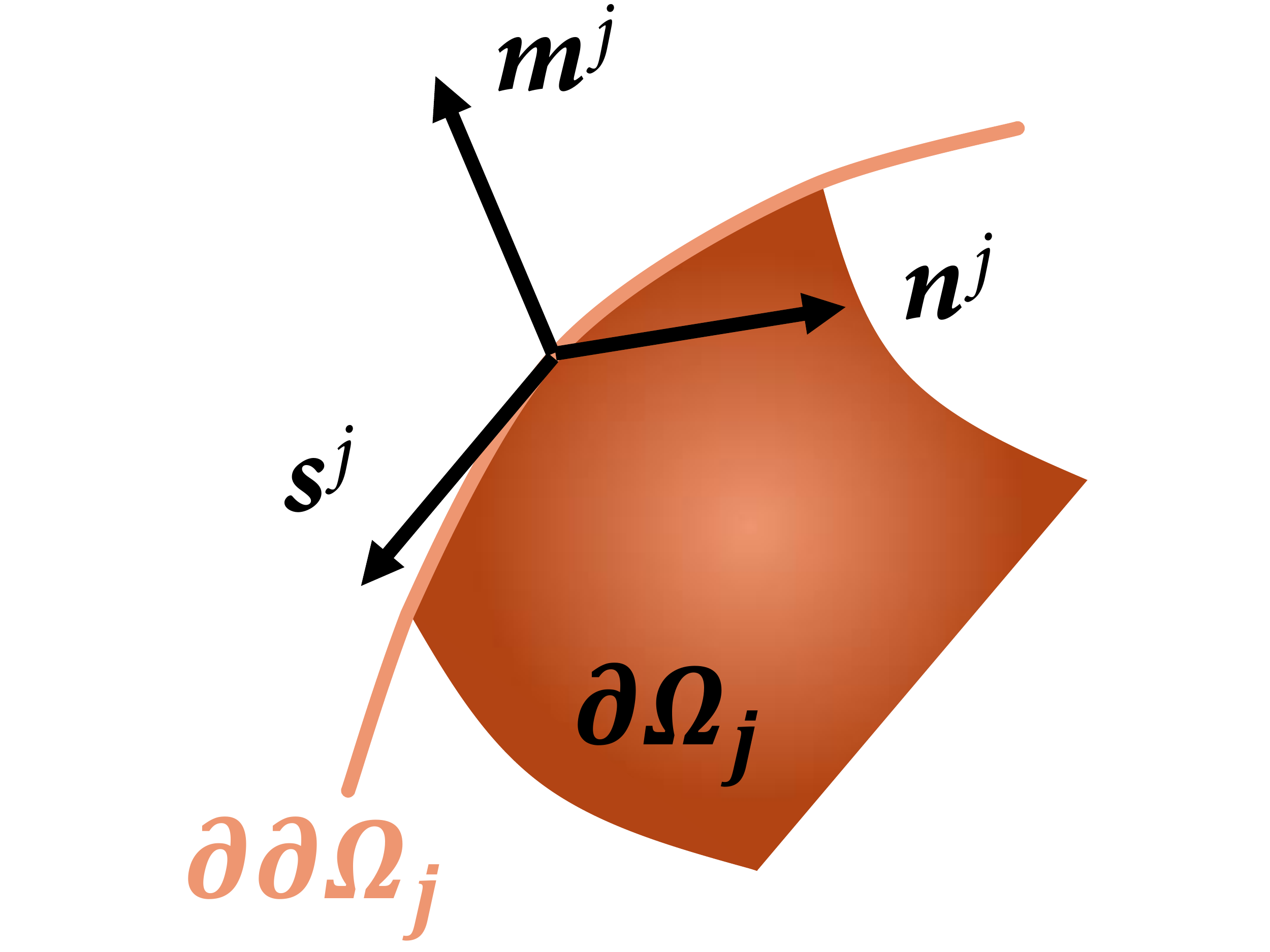}
		\vfill\label{fig_2_1_4}
	\end{subfigure}
	\caption{Sketch of the geometry of $\Omega\in\mathbb{R}^3$.
		%a) Physical domain $\Omega\in\mathbb{R}^3$, b) 
		\emph{a)} Detail of $\partial\Omega$ subdivided in smooth portions $\partial\Omega_i$ and $\partial\Omega_j$, with their corresponding normal vectors $\n^\toVect{i}$ and $\n^\toVect{j}$, \emph{b)} detail of $\partial\Omega_i$, with the triplet $\{\m^\toVect{i},\tv^\toVect{i},\n^\toVect{i}\}$ defined on $\partial\partial\Omega_i$, and \emph{c)} detail of $\partial\Omega_j$, with the triplet $\{\m^\toVect{j},\tv^\toVect{j},\n^\toVect{j}\}$ defined on $\partial\partial\Omega_j$.}
	\label{fig_2_1}
\end{figure}

%%%%%
Let $\Omega$ be a physical domain in $\mathbb{R}^3$
. The domain boundary, $\partial\Omega$, can be conformed by several smooth portions as $\partial\Omega=\bigcup_f\partial\Omega_f$ (\fig\ref{fig_2_1_2}).
At each point $\x\in\partial\Omega_f$ we define $\n^f$ as the outward unit \emph{normal vector}.
The boundary of the $f$-th portion of $\partial\Omega$ is denoted as $\partial\partial\Omega_f$, which is a closed curve. At each point $\x\in \partial\partial\Omega_f$ we define $\m^f$ as the unit \emph{co-normal vector} pointing outwards of $\partial\Omega_f$, which is orthogonal to the normal vector $\n^f$ and to the tangent vector of the curve $\partial\partial\Omega_f$, $\tv^f$ (see \fig\ref{fig_2_1_3} and \ref{fig_2_1_4}). The orientation of $\tv^f$ is arbitrary and not relevant in the derivations next.

The operators that appear throughout this Section are defined next.
The \emph{spatial derivative} of a function $f$ with respect to the coordinate $x_i$ is denoted by $\partial_i$ or the subindex $i$ after a comma, that is, $\frac{\partial}{\partial x_i}(f)=\partial_if=f_{,i}$.
The \emph{gradient} operator is denoted as $\gradient(~~)$, and the \emph{divergence} operator as $\divergence(~~)$. For instance, for a second order tensor $\toMat{A}$, they are defined as $[\gradient(\toMat{A})]_{ijk}\coloneqq\partial_kA_{ij}=A_{ij,k}$, and $[\divergence(\toMat{A})]_i\coloneqq\partial_jA_{ij}=A_{ij,j}$, respectively.\\
The \emph{symmetrized gradient} $\symmetricGradient{(~~)}$ of a vector field $\toVect{a}$ is defined as:
\begin{equation}
\left[\symmetricGradient{(\toVect{a})}\right]_{ij}\coloneqq\frac{1}{2}\left[\gradient{(\toVect{a})}+\gradient{(\toVect{a})}^T\right]_{ij}=\frac{1}{2}(a_{i,j}+a_{j,i}).
\end{equation}
On the domain boundary $\partial\Omega$ the derivative in the normal direction, namely the \emph{normal derivate}, is denoted by $\partial^n$. For instance, the normal derivative of a vector field $\toVect{a}$ is defined as $\partial^n(\toVect{a})\coloneqq\gradient(\toVect{a})\cdot\n$.
The gradient and divergence operators can be decomposed on $\partial\Omega$ into their normal and tangential components as $\gradient(~~)=\partial^n(~~)\n+\surfaceGradient(~~)$ and $\divergence(~~)=\partial^n(~~)\cdot\n+\surfaceDivergence(~~)$, respectively,
where $\surfaceGradient(~~)$ and $\surfaceDivergence(~~)$ denote the \emph{surface gradient}  and \emph{surface divergence} operators, namely the projection of the gradient and divergence operators onto the tangent space of $\partial\Omega$.
For a second order tensor $\toMat{A}$ they are expressed as
\begin{subequations}\begin{align}
\left[\surfaceGradient(\toMat{A})\right]_{ijk}=\surfaceGradient_kA_{ij}\coloneqq&
%\left[\gradient{(\toMat{A})}\cdot\Projector\right]_{ijk}=
A_{ij,l}P_{lk},\\
\left[\surfaceDivergence(\toMat{A})\right]_i=\surfaceGradient_jA_{ij}\coloneqq&
A_{ij,k}P_{jk},
\end{align}\end{subequations}
respectively, where $\Projector$ is the \emph{projection operator} defined on $\partial\Omega$ as 
$\left[\Projector\right]_{ij}\coloneqq%\toMat{I}-\n\otimes\n=
\delta_{ij}-n_in_j$, being $\delta_{ij}$ the Kronecker delta.

The formulation involves second-order measures of the geometry, namely curvatures of $\partial\Omega$. The tensor that contains this information is known as the \emph{shape operator} $\ShapeOperator$ (also known as curvature tensor) \cite{Block1997}, defined on the surface $\partial\Omega$ as
\begin{equation}
[\ShapeOperator]_{ij}\coloneqq-\left[\surfaceGradient(\n)\right]_{ij}=-n_{i,l}P_{lj}.
\end{equation}
The \emph{mean curvature} $H$ of a surface is an invariant, expressed in terms of $\ShapeOperator$ as
\begin{equation}
H\coloneqq\frac{1}{2}\trace{\ShapeOperator}=-\frac{1}{2}\surfaceDivergence\n=-\frac{1}{2}n_{i,j}P_{ij}.
\end{equation}
With $\ShapeOperator$ and $H$ we define a tensor which arises in the formulation 
(see \ref{sec_app01}), that we name \emph{second-order geometry tensor} $\CurvatureProjector$ and is defined as:
\begin{equation}
\left[\CurvatureProjector\right]_{ij}\coloneqq S_{ij}-2Hn_in_j.
\end{equation}
On the curve $C^{ij}\coloneqq\partial\partial\Omega_i\cap\partial\partial\Omega_j$, $(i\neq j)$ we define the \emph{jump} operator $\jump{~~}$ acting on a given quantity $\toVect{a}$ as the sum of that quantity evaluated at both sides of the curve, namely
$\jump{\toVect{a}}\coloneqq\toVect{a}^i+\toVect{a}^j$, where $\toVect{a}^k$ is the value of $\toVect{a}$ from $\partial\Omega_k$ (see \fig\ref{fig_2_1_2}).
Note that a \emph{sum} is considered in the definition of the jump, so that $\jump{\m\times\n}$ vanishes in the case $C^{ij}$ is defined along a smooth region of $\partial\Omega$,
where $\m^i=-\m^j$ and $\n^i=\n^j$.

Finally, we denote the \emph{first variation} of a certain functional $F[f_1,\dots,f_n]$ with respect to the function $f_i$ as the functional
\begin{equation}
\var[f_i]{F}[f_1,\dots,f_n;\psi_i]
\coloneqq
\frac{d}{d\epsilon}F[f_1,\dots,f_{i-1},f_i+\epsilon\psi_i,f_{i+1},\dots,f_n]\Big{|}_{\epsilon=0},
\end{equation}
where $\epsilon\in\mathbb{R}$, and the function $\psi_i$ is the variation of $f_i$ (hence denoted also as ${\delta\!f}_i$), defined on the same functional space as $f_i$. The variation of $F[f_1,\dots,f_n]$ with respect to all the functions $f_1,\dots,f_n$ is denoted as
\begin{equation}
\var{F}[f_1,\dots,f_n;\psi_1,\dots,\psi_n]
\coloneqq\sum_{i=1}^n \var[f_i]{F}[f_1,\dots,f_n;\psi_i].
\end{equation}
The \emph{second variation} of $F[f_1,\dots,f_n]$ with respect to the function $f_i$ in the direction $\psi_i$ is denoted by $\vvar[f_i]{F}[f_1,\dots,f_n;\psi_i]$, and is defined as the first variation of the functional $\var[f_i]{F}[f_1,\dots,f_n;\psi_i]$ with respect to $f_i$, \ie
\begin{equation}
\vvar[f_i]{F}[f_1,\dots,f_n;\psi_i]
\coloneqq
\var[f_i]{\Big(\var[f_i]{F}[f_1,\dots,f_n;\psi_i]\Big)}[f_1,\dots,f_n;\psi_i],
\end{equation}
where $\psi_i$ is the variation of the function $f_i$ for both variations of $F$.
%%%
\subsection{Standard variational formulation}
%%%
A continuum model for flexoelectric materials can be obtained by coupling a strain-gradient elasticity model with classical electrostatics, through the piezoelectric and flexoelectric effects.
The state variables are the displacement field $\Displacement$ and the electric potential $\phi$.
The strain tensor and the electric field are given by 
\begin{equation}
[\Strain(\Displacement)]_{ij}=[\Strain(\Displacement)]_{ji}\coloneqq\left[\symmetricGradient(\Displacement)\right]_{ij}=\frac{1}{2}(u_{i,j}+u_{j,i}),\label{eq_str}
\end{equation}
\begin{equation}
[\E(\phi)]_l\coloneqq-[\gradient\phi]_l=-\phi_{,l}.
\end{equation}
The bulk energy density $\mathcal{H}^\Omega$ in a flexoelectric material can be stated as \cite{Mao2014,Abdollahi2014, Abdollahi2015a, Abdollahi2015b}
\begin{multline}\label{eq_FlexRawEnergy}
\mathcal{H}^\Omega[\Displacement,\phi]=\mathcal{H}^\Omega[\Strain,\gradient\Strain,\E]\coloneqq\frac{1}{2}\strain_{ij}\elast_{ijkl}\strain_{kl}+\frac{1}{2}\strain_{ij,k}\strGr_{ijklmn}\strain_{lm,n}\\-\frac{1}{2}E_{l}\dielec_{lm}E_{m}-E_{l}\piezo_{lij}\strain_{ij}-E_{l}\flexo_{lijk}\strain_{ij,k}.
\end{multline}
The first two terms correspond to the mechanical energy density of a strain-gradient elastic material, in the Form II of the original paper of Mindlin \cite{Mindlin1964} about strain gradient elasticity. The tensor $\elast_{ijkl}=\elast_{klij}=\elast_{jikl}=\elast_{ijlk}$ is the fourth-order \emph{elasticity tensor} and $\strGr_{ijklmn}=\strGr_{lmnijk}=\strGr_{jiklmn}=\strGr_{ijkmln}$ is the sixth-order \emph{strain-gradient elasticity tensor}.
The third term is the electrostatic energy density, where $\dielec_{lm}=\dielec_{ml}$ is the second-order \emph{dielectricity tensor}. The last two terms correspond to the piezoelectric and flexoelectric effects, where $\piezo_{lij}=\piezo_{lji}$ is the third-order \emph{piezoelectric tensor} and $\flexo_{lijk}=\flexo_{ljik}$ the fourth-order \emph{flexoelectric tensor}.

Alternative descriptions could also be considered, writing the bulk energy density in terms of other quantities instead of the electric field $\E$, such as the polarization $\toVect{p}$ or the electric displacement $\ElectricDisp=\varepsilon_0\E+\toVect{p}$, where $\varepsilon_0$ is the vacuum permittivity constant. As argued in \cite{Liu2014}, all of them are valid, but the choice of $\E$ facilitates the derivation of equilibrium equations, which are simpler than those in terms of $\toVect{p}$, allowing for simpler numerical methods for solving the associated boundary value problems.
Some authors \cite{Shen2010,Sharma2010,Majdoub2008,Majdoub2009} describe the energy density as $\mathcal{H}^\Omega[\Strain,\gradient\Strain,\toVect{p},\gradient\toVect{p}]$ accounting for the \emph{polarization gradient} theory instead of strain gradient elasticity. We consider the simplified model in \eq\eqref{eq_FlexRawEnergy} since \emph{i)} the quadratic term to the electric field gradient is assumed to be negligible for the problems we are interested in, and \emph{ii)} the last term in \eq\eqref{eq_FlexRawEnergy} describes both the \emph{direct} and \emph{converse} flexoelectric effects as a Lifshitz invariant (see \cite{Landau2013,Sharma2010,Abdollahi2014} for details), leading to the same equilibrium equations as if both terms were treated separately.
We also note that either non-local mechanical or electrical effects (or both) must be introduced in the flexoelectric problem in order to get a meaningful energy density from a physical point of view, but also to get numerically stable formulations. It's worth noting that incorporating the polarization gradient theory to the present formulation is straightforward by following the same derivations highlighted in the present paper.

The contribution from external loads is presented next. Being $\toVect{b}$ the \emph{body force} per unit volume, and $q$ the \emph{free charge} per unit volume, their work per unit volume is
\begin{equation}\label{eq_ExternalBulk}
\mathcal{W}^{\Omega}[\Displacement,\phi]\coloneqq -b_iu_i+q\phi.
\end{equation}

Additional external loads are present on the domain boundary. In a strain gradient elasticity formulation \cite{Mindlin1964}, those are the \emph{traction} $\Traction$ and the \emph{double traction} $\HighOrderTraction$, which are the conjugates of the displacement $\Displacement$ and the normal derivative of the displacement $\partial^n(\Displacement)$ on $\partial\Omega$, respectively. The electrical boundary load is the \emph{surface charge density} $w$, which is conjugate of the electric potential $\phi$ on $\partial\Omega$. The work of the external loads per unit area is
\begin{equation}\label{eq_ExternalBoundary}
\mathcal{W}_0^{\partial\Omega}[\Displacement,\phi]\coloneqq -t_iu_i-r_i\partial^nu_i+w\phi.
\end{equation}

Moreover, as dictated by strain gradient elasticity theory \cite{Mindlin1964}, an additional \emph{force per unit length} $\HighOrderJump$ arises at the edges $C$ of the boundary, \ie at the union of the edges formed by the intersection of the portions of the boundary, in case $\partial\Omega$ is not smooth (see \fig\ref{fig_2_1_2}). That is, at $C=\bigcup_f\partial\partial\Omega_f$. The force $\HighOrderJump$ is the conjugate of the displacement $\Displacement$ on $C$. Hence, its work per unit length is
\begin{equation}\label{eq_ExternalLine}
\mathcal{W}_0^C[\Displacement,\phi]\coloneqq -j_iu_i.
\end{equation}

The total energy $\Pi_0[\Displacement,\phi]$ of a flexoelectric material is found by collecting all the internal and external energy densities as follows:

\begin{equation}\label{eq_Functional}
\Pi_0[\Displacement,\phi]=
\int_{\Omega}\Big(\mathcal{H}^\Omega[\Displacement,\phi]+\mathcal{W}^{\Omega}[\Displacement,\phi]\Big)\dd\Omega
+
\int_{\partial\Omega}\mathcal{W}^{\partial\Omega}_0[\Displacement,\phi]\dd\Gamma
+
\int_{C}\mathcal{W}^C_0[\Displacement,\phi]\dd\ds.
\end{equation}

The boundary of the domain $\partial\Omega$ can be split into several disjoint regions, corresponding to the different \emph{Dirichlet} and \emph{Neumann} boundaries. The external load is prescribed on the latter, whereas on the former its conjugate is prescribed.
For the mechanical loads, we have 
$\partial\Omega=\partial\Omega_u\cup\partial\Omega_t$
and
$\partial\Omega=\partial\Omega_v\cup\partial\Omega_r$.
The electrical part is also split into 
$\partial\Omega=\partial\Omega_\phi\cup\partial\Omega_w$.
The corresponding boundary conditions are
\begin{align}\label{strgr_classicalBC}
\Displacement - \DisplacementKnown &= \mathbf{0} \quad\text{on }{\partial\Omega_u},
&&&
\Traction(\Displacement,\phi) - \TractionKnown &= \mathbf{0} \quad\text{on }{\partial\Omega_t};
\\\label{strgr_nonlocalBC}
\partial^n(\Displacement) - \partialnDisplacementKnown &= \mathbf{0} \quad\text{on }{\partial\Omega_v},
&&&
\HighOrderTraction(\Displacement,\phi) - \HighOrderTractionKnown &= \mathbf{0} \quad\text{on }{\partial\Omega_r};
\\\label{elecBC}
\phi-\phiKnown &=0\quad\text{on }{\partial\Omega_\phi},
&&&
w(\Displacement,\phi)-\wKnown &=0\quad\text{on }{\partial\Omega_w};
\end{align}
where $\DisplacementKnown$, $\partialnDisplacementKnown$ and $\phiKnown$
are the prescribed
\emph{displacement}, \emph{normal derivative of the displacement} and \emph{electric potential} at the Dirichlet boundaries, and 
$\TractionKnown$, $\HighOrderTractionKnown$ and $\wKnown$
the prescribed \emph{traction}, \emph{double traction} and \emph{surface charge density} at the Neumann boundaries. The expressions $\Traction(\Displacement,\phi)$, $\HighOrderTraction(\Displacement,\phi)$ and $w(\Displacement,\phi)$ will be derived later as a result of the variational principle in \eq\eqref{eq_variationalPrinciple}.

The edges $C$ of $\partial\Omega$ are also split into $C=C_u\cup C_j$ corresponding to the Dirichlet and Neumann edge partitions, respectively. Here, $C_u$ is assumed to correspond to the curves within the classical Dirichlet boundary, namely $C_u=C\cap\overline{\partial\Omega_u}$, and $C_v=C\setminus C_u$. Edge boundary conditions are:
\begin{align}\label{strgr_edgeBC}
\Displacement - \DisplacementKnown &= \mathbf{0}\quad\text{on }{C_u},
&&&
\HighOrderJump(\Displacement,\phi) - \HighOrderJumpKnown &= \mathbf{0}\quad\text{on }{C_j},
\end{align}
where $\HighOrderJumpKnown$ is the prescribed \emph{force per unit length} at the Neumann edges, and the expression $\HighOrderJump(\Displacement,\phi)$ will be derived later from \eq\eqref{eq_variationalPrinciple}.
Many authors in the literature neglect the edge conditions in \eq\eqref{strgr_edgeBC} \cite{Aravas2011,Mao2014,Abdollahi2014,Abdollahi2015a}. It is important to note that dismissing them is equivalent to considering homogeneous Neumann edge conditions, which may not be true on $C_u$ (Dirichlet edges).
In this work, the edge conditions are kept in the formulation to ensure self-consistency and a well-defined boundary value problem.

The energy functional in \eq\eqref{eq_Functional} can be rewritten as follows, according to \eq(\ref{strgr_classicalBC}b)-(\ref{strgr_edgeBC}b):
\begin{equation}\label{eq_Functional2}
\Pi_0[\Displacement,\phi]=
\Pi^\Omega[\Displacement,\phi]
+
\Pi_0^\text{Dirichlet}[\Displacement,\phi]
+
\Pi^\text{Neumann}[\Displacement,\phi],
\end{equation}
where
\begin{equation}
\Pi^\Omega[\Displacement,\phi]=\int_{\Omega}\Big(\mathcal{H}^\Omega[\Displacement,\phi]+\mathcal{W}^{\Omega}[\Displacement,\phi]\Big)\dd\Omega,
\end{equation}
\begin{multline}\label{Dirichlet_0}
\Pi_0^\text{Dirichlet}[\Displacement,\phi]=
  \int_{\partial\Omega_u}-u_it_i(\Displacement,\phi)\dd\Gamma
  +
  \int_{\partial\Omega_v}-\partial^nu_ir_i(\Displacement,\phi)\dd\Gamma
  \\+
  \int_{\partial\Omega_\phi}\phi w(\Displacement,\phi)\dd\Gamma
  +
  \int_{C_u}-u_ij_i(\Displacement,\phi)\dd\ds
\end{multline}
and
\begin{equation}
\Pi^\text{Neumann}[\Displacement,\phi]=
\int_{\partial\Omega_t}-u_i\tractionKnown_i\dd\Gamma
+
\int_{\partial\Omega_r}-\partial^nu_i\highOrderTractionKnown_i\dd\Gamma
+
\int_{\partial\Omega_w}\phi\wKnown\dd\Gamma
+
\int_{C_j}-u_i\highOrderJumpKnown_i\dd\ds.
\end{equation}
The state variables 
$(\Displacement,\phi)\in\mathcal{U}_0\otimes\mathcal{P}_0$,
where
\begin{align}
\mathcal{U}_0 &\coloneqq\{\Displacement\in[\mathcal{H}^2(\Omega)]^3\text{ $|$ }
\Displacement-\DisplacementKnown=\mathbf{0}\text{ on }\partial\Omega_u,
\Displacement-\DisplacementKnown=\mathbf{0}\text{ on }C_u
\text{ and }
\partial^n\Displacement-\partialnDisplacementKnown=\mathbf{0}\text{ on }\partial\Omega_v\},
\\
\mathcal{P}_0 &\coloneqq\{\phi\in\mathcal{H}^1(\Omega)\text{ $|$ }
\phi-\phiKnown=0\text{ on }\partial\Omega_\phi\};
\end{align}
fulfilling Dirichlet boundary conditions in \eq (\ref{strgr_classicalBC}a)-(\ref{strgr_edgeBC}a).

The equilibrium states $(\Displacement^*,\phi^*)$ of the body correspond to the following variational principle:
\begin{equation}\label{eq_variationalPrinciple}
(\Displacement^\text*,\phi^\text*)=\arg\min_{\Displacement\in\mathcal{U}_0}\max_{\phi\in\mathcal{P}_0}\Pi_0[\Displacement,\phi].
\end{equation}

The Euler-Lagrange equations associated with the variational principle in \eq\eqref{eq_variationalPrinciple} and the expressions $\Traction(\Displacement,\phi)$, $\HighOrderTraction(\Displacement,\phi)$, $w(\Displacement,\phi)$ and $\HighOrderJump(\Displacement,\phi)$ from the Neumann boundary conditions in \eq (\ref{strgr_classicalBC}b)-(\ref{strgr_edgeBC}b) are found by enforcing
\begin{subequations}\begin{align}
\var{\Pi_0}=\var[\Displacement]{\Pi_0}+\var[\phi]{\Pi_0}=0
;\label{eq_1var}
\\
\vvar[\Displacement]{\Pi_0}>0,\quad
\vvar[\phi]{\Pi_0}<0
,\label{eq_2var}
\end{align}\end{subequations}
for all admissible variations
$\delta\Displacement\in\mathcal{U}_0$ and $\delta\phi\in\mathcal{P}_0$. 
The full derivation can be found in \cite{Liu2014}, and the resulting equations are given next \cite{Liu2014,Mao2014,Abdollahi2014, Abdollahi2015a, Abdollahi2015b}:
\begin{equation}\label{eq_EulerLagrange}
\begin{cases}
\left(\cauchyStress_{ij}(\Displacement,\phi)-\hyperStress_{ijk,k}(\Displacement,\phi)\right)_{,j}+b_i=0_i &\text{in }\Omega,\\
\hfill\electricDisp_{l,l}(\Displacement,\phi)-q=0\hphantom{{}_i} & \text{in }\Omega;
\end{cases}
\end{equation}
and
\begin{subnumcases}{\label{eq_FlexoForces}}
t_i(\Displacement,\phi) \coloneqq \left(\cauchyStress_{ij}(\Displacement,\phi)-\hyperStress_{ijk,k}(\Displacement,\phi)-\surfaceGradient_{k}\hyperStress_{ikj}(\Displacement,\phi)\right)n_j+\hyperStress_{ijk}(\Displacement,\phi)\curvatureProjector_{jk}
& on $\partial\Omega$ \label{eq_Traction},\\
r_i(\Displacement,\phi) \coloneqq \hyperStress_{ijk}(\Displacement,\phi)n_jn_k
& on $\partial\Omega$, \\
w(\Displacement,\phi) \coloneqq -\electricDisp_l(\Displacement,\phi)n_l
& on $\partial\Omega$, \\
j_i(\Displacement,\phi) \coloneqq \jump{\hyperStress_{ijk}(\Displacement,\phi)m_jn_k}
& on $C$.
\end{subnumcases}

Note that the traction $\Traction(\Displacement,\phi)$ in \eq\eqref{eq_Traction} is an alternative expression to the one in \cite{Mindlin1964}, whose derivation can be found in \ref{sec_app01}.

In \eq\eqref{eq_EulerLagrange} and \eqref{eq_FlexoForces}, the \emph{stress} $\CauchyStress(\Displacement,\phi)$, the \emph{double stress} $\HyperStress(\Displacement,\phi)$ and the \emph{electric displacement} $\ElectricDisp(\Displacement,\phi)$ are defined as the conjugates to the strain $\Strain(\Displacement)$, the strain gradient $\gradient\Strain(\Displacement)$ and the electric field $\E(\phi)$, respectively, as follows:
\begin{subequations}\label{elec_tensors}\begin{align}
\cauchyStress_{ij}(\Displacement,\phi)=\cauchyStress_{ji}(\Displacement,\phi)
&\coloneqq
\left.\frac{\partial \mathcal{H}^\Omega[\Strain,\gradient\Strain,\E]}{\partial\strain_{ij}}\right|_{\substack{\Strain=\Strain(\Displacement)\\\E=\E(\phi)}}
=\elast_{ijkl}\strain_{kl}(\Displacement)-\piezo_{lij}E_{l}(\phi);
\\
\hyperStress_{ijk}(\Displacement,\phi)=\hyperStress_{jik}(\Displacement,\phi)
&\coloneqq
\left.\frac{\partial\mathcal{H}^\Omega[\Strain,\gradient\Strain,\E]}{\partial\strain_{ij,k}}\right|_{\substack{\gradient\Strain=\gradient\Strain(\Displacement)\\\E=\E(\phi)}}
=\strGr_{ijklmn}\strain_{lm,n}(\Displacement)-\flexo_{lijk}E_{l}(\phi);
\\
\electricDisp_l(\Displacement,\phi)
&\coloneqq
\left.-\frac{\partial\mathcal{H}^\Omega[\Strain,\gradient\Strain,\E]}{\partial E_{l}}\right|_{\substack{\Strain=\Strain(\Displacement)\\\E=\E(\phi)\\\gradient\Strain=\gradient\Strain(\Displacement)}}
=\dielec_{lm}E_{m}(\phi)+\piezo_{lij}\strain_{ij}(\Displacement)+\flexo_{lijk}\strain_{ij,k}(\Displacement).
\end{align}\end{subequations}

The positivity and negativity conditions on the second variations in \eq\eqref{eq_2var} lead to the following restrictions on the material tensors:
\begin{align}\label{restrictions}
\dielec_{ii}>0, \qquad
\elast_{ijij}>0, \qquad
\strGr_{ijkijk}\geq0
; &&
&& \quad i,j,k
%,l
&=1,\dots,n_d,
\end{align}
being $n_d$ the number of spatial dimensions. 

\begin{figure}[t]\centering
	\begin{subfigure}[t]{.3\textwidth}\centering
		\captionof{figure}{}\vspace{-0.5em}
		\includegraphics[scale=0.17]{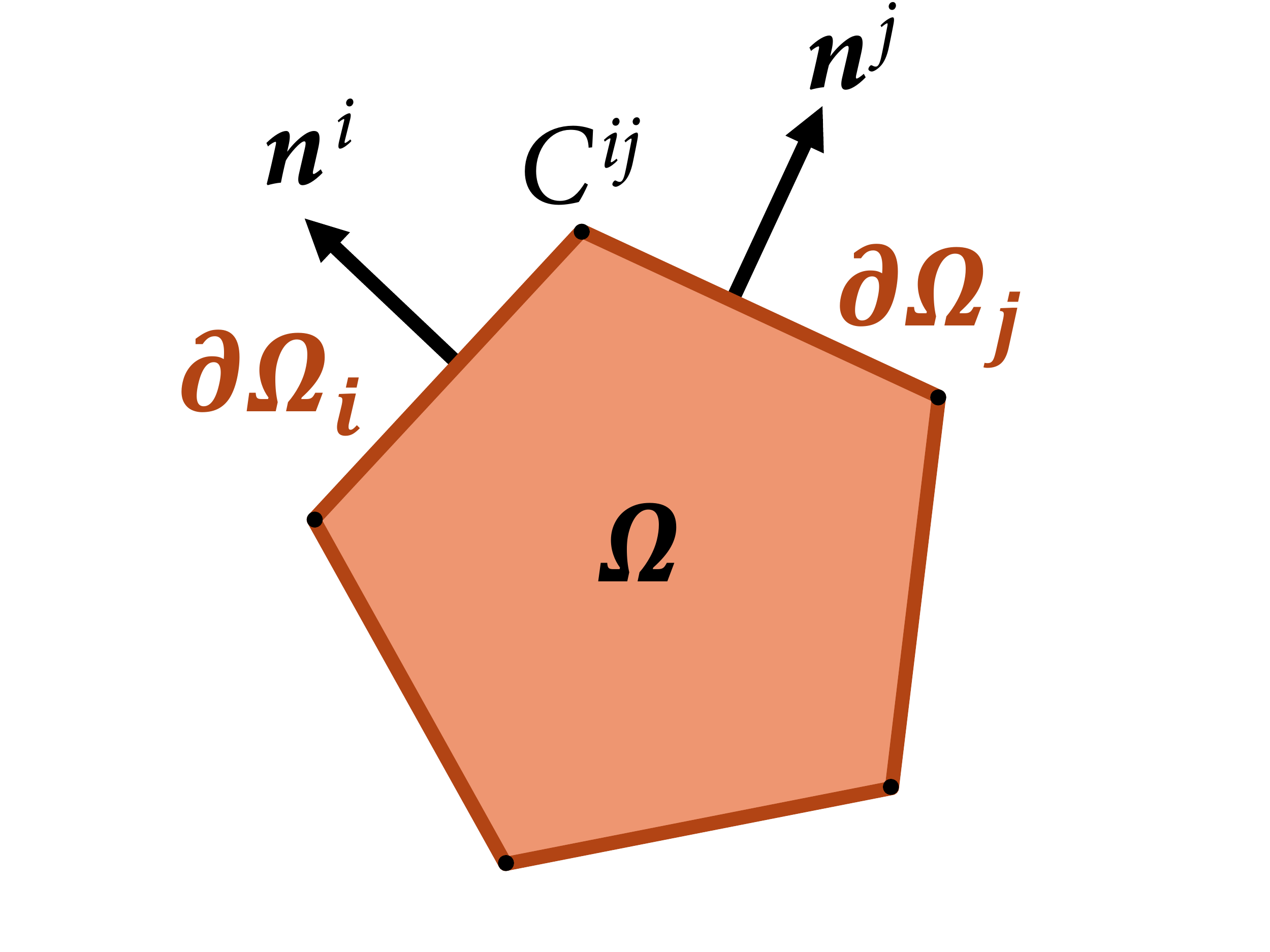}
		\vfill\label{fig_3_1_2}
	\end{subfigure}%
	\begin{subfigure}[t]{.3\textwidth}\centering
		\captionof{figure}{}\vspace{-0.5em}
		\includegraphics[scale=0.17]{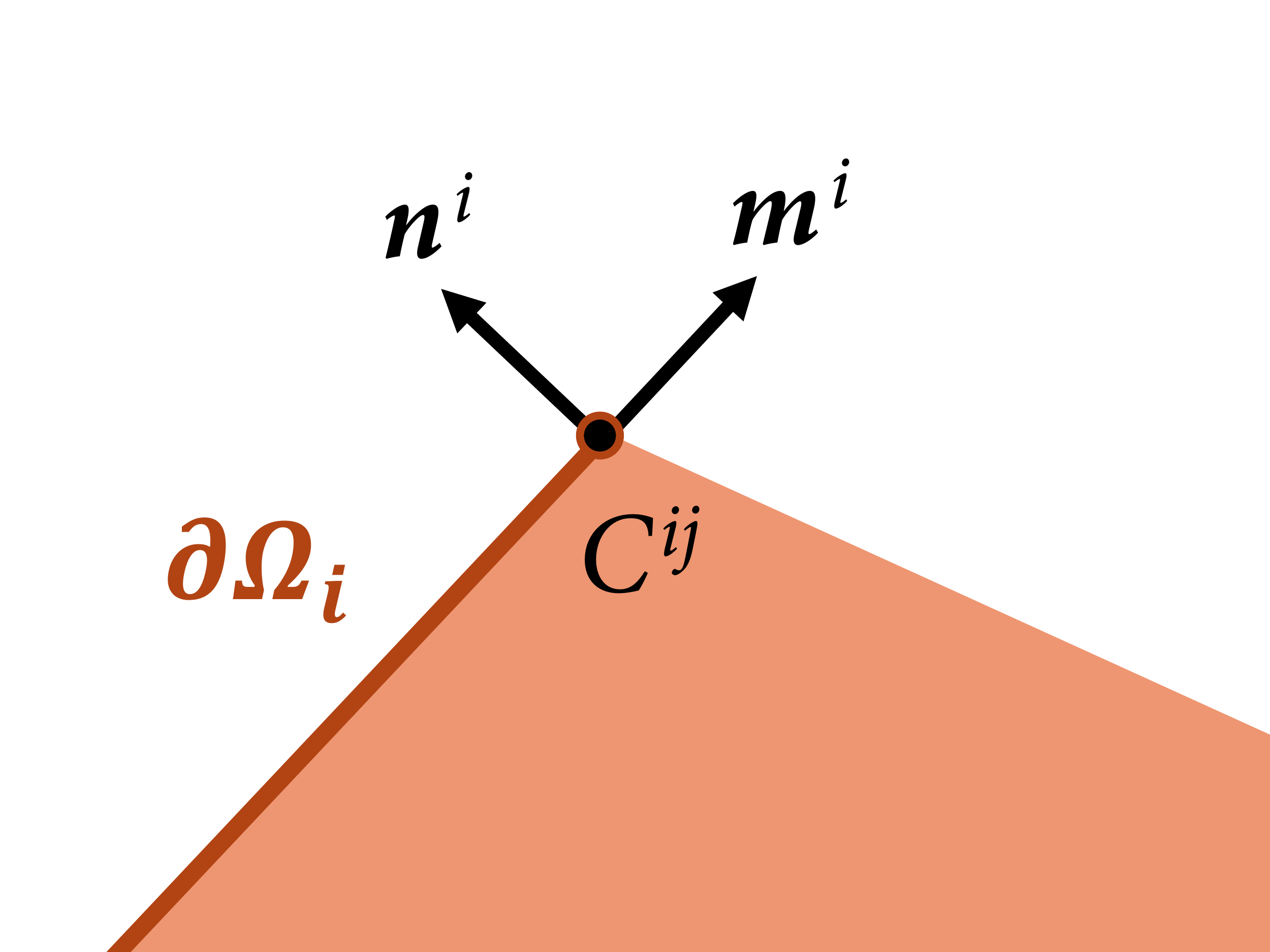}
		\vfill\label{fig_3_1_3}
	\end{subfigure}%
	\begin{subfigure}[t]{.3\textwidth}\centering
		\captionof{figure}{}\vspace{-0.5em}
		\includegraphics[scale=0.17]{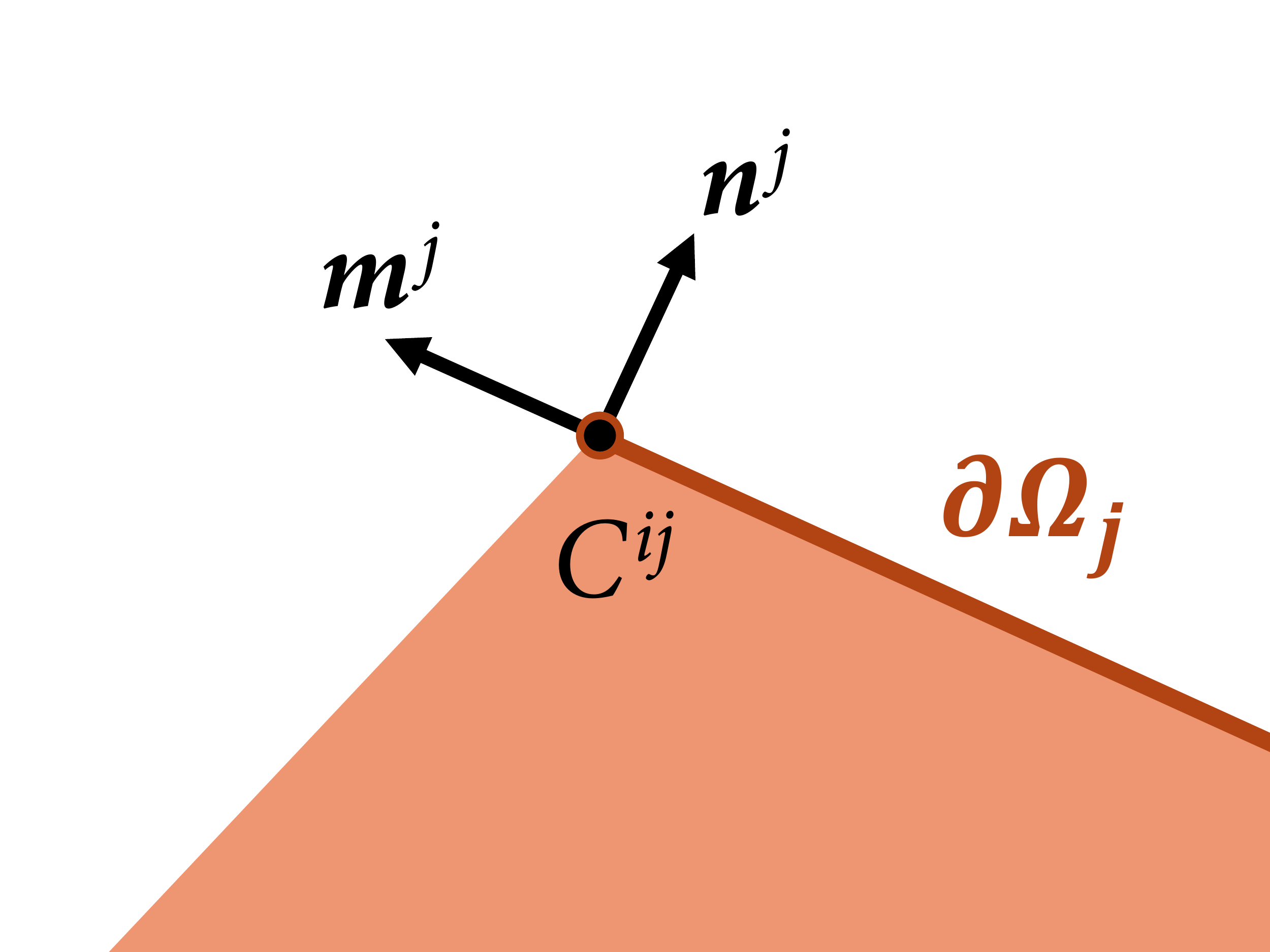}
		\vfill\label{fig_3_1_4}
	\end{subfigure}
	\caption{Sketch of the geometry of $\Omega\in\mathbb{R}^2$.
		\emph{a)} $\Omega$ and the boundary $\partial\Omega$ subdivided in several smooth portions $\partial\Omega_f$, \eg $\partial\Omega_i$ and $\partial\Omega_j$ with their corresponding normal vectors $\n^\toVect{i}$ and $\n^\toVect{j}$. The intersection between $\partial\Omega_i$ and $\partial\Omega_j$ is the corner $C^{ij}$. \emph{b)} Detail of $\partial\Omega_i$, with the pair $\{\m^\toVect{i},\n^\toVect{i}$ defined on $C^{ij}$, and \emph{c)} detail of $\partial\Omega_j$, with the pair $\{\m^\toVect{j},\n^\toVect{j}\}$ defined on $C^{ij}$.}
	\label{fig_3_1}
\end{figure}

\begin{remark}
	This formulation corresponds to the 3D case. It also holds for 2D with the following considerations:
	\vspace{-0.5em}\begin{itemize}\setlength\itemsep{-0.3em}
		\item The vectors $\n$ and $\m$ defined on the \emph{curve} $\partial\Omega_f$ refer to the outward unit normal and tangent vectors, respectively (see \fig\ref{fig_3_1}).
		\item The edges $C$ in \eq\eqref{eq_Functional} correspond to \emph{corners} of $\partial\Omega$, and therefore have dimension 0. As a consequence, their contribution to the external work $\int_{C}\HighOrderJump\delta\Displacement\dd\ds$ is $\sum_C\HighOrderJump\delta\Displacement$ in \eq\eqref{eq_Functional} and subsequent terms.
	\end{itemize}
\end{remark}
%%%%%
\subsection{Variational formulation within an unfitted framework: The Nitsche's method}
%%%%%
The admissible space $\mathcal{U}_0\otimes\mathcal{P}_0$ of the state variables $(\Displacement,\phi)$ is constrained on the Dirichlet boundaries, and therefore is not suitable for an unfitted formulation. In order to overcome this requirement, an alternative energy functional $\Pi[\Displacement,\phi]$ is proposed following Nitsche's approach \cite{Nitsche1971}:
\begin{equation}\label{eq_Functionalimm}
\Pi[\Displacement,\phi]=
\Pi^\Omega[\Displacement,\phi]
+
\Pi^\text{Dirichlet}[\Displacement,\phi]
+
\Pi^\text{Neumann}[\Displacement,\phi],
\end{equation} 
where $\Pi^\text{Dirichlet}[\Displacement,\phi]$ acts on the Dirichlet boundaries instead of $\Pi_0^\text{Dirichlet}[\Displacement,\phi]$, and incorporates Dirichlet boundary conditions in \eq(\ref{strgr_classicalBC}a)-(\ref{strgr_edgeBC}a) weakly as follows:
\begin{multline}\label{Dirichlet}
\Pi^\text{Dirichlet}[\Displacement,\phi]=
\int_{\partial\Omega_u}\left(\frac{1}{2}\beta_u\Big(u_i-\displacementKnown_i\Big)^2-\Big(u_i-\displacementKnown_i\Big)t_i(\Displacement,\phi)\right)\dd\Gamma
+{}\\{}+
\int_{\partial\Omega_v}\left(\frac{1}{2}\beta_v\Big(\partial^nu_i-\partialndisplacementKnown_i\Big)^2-\Big(\partial^nu_i-\partialndisplacementKnown_i\Big)r_i(\Displacement,\phi)\right)\dd\Gamma
+{}\\{}+
\int_{\partial\Omega_\phi}\left(-\frac{1}{2}\beta_\phi\Big(\phi-\phiKnown\Big)^2+\Big(\phi-\phiKnown\Big) w(\Displacement,\phi)\right)\dd\Gamma
+{}\\{}+
\int_{C_u}\left(\frac{1}{2}\beta_{C_u}\Big(u_i-\displacementKnown_i\Big)^2-\Big(u_i-\displacementKnown_i\Big)j_i(\Displacement,\phi)\right)\dd\ds,
\end{multline}
with the numerical parameters $\beta_u,\beta_v,\beta_{C_u},\beta_\phi\in\mathds{R}^+$.

Comparing \eq\eqref{Dirichlet} against \eq\eqref{Dirichlet_0}, one can readily see that the expressions $\Traction(\Displacement,\phi)$, $\HighOrderTraction(\Displacement,\phi)$, $w(\Displacement,\phi)$ and $\HighOrderJump(\Displacement,\phi)$ are now conjugate to the Dirichlet boundary conditions.
The penalty terms inserted in each boundary integral are quadratic in the Dirichlet boundary conditions, and its purpose is to weakly enforce Dirichlet boundary conditions and to ensure equilibrium states $(\Displacement^*,\phi^*)$ being, respectively, actual minima and maxima of the energy functional with respect to $\Displacement$ and $\phi$.

The variational principle associated to $\Pi[\Displacement,\phi]$ for the equilibrium states $(\Displacement^*,\phi^*)$ of the body is stated next:
\begin{equation}\label{variationalPrinciple}
(\Displacement^\text*,\phi^\text*)=\arg\min_{\Displacement\in\mathcal{U}}\max_{\phi\in\mathcal{P}}\Pi[\Displacement,\phi],
\end{equation}
where $\mathcal{P}\coloneqq\mathcal{H}^1(\Omega)$, and $\mathcal{U}$ is the space of functions belonging to $\big[\mathcal{H}^2(\Omega)\big]^3$ with $\mathcal{L}^2$-integrable third derivatives on the boundary $\partial\Omega_u$, to account for the integrals involving $\Traction(\Displacement,\phi)$ in \eq\eqref{eq_Traction}.
The variational principle in \eq\eqref{variationalPrinciple} leads to the same Euler-Lagrange equations in \eq\eqref{eq_EulerLagrange} and definitions of $\Traction(\Displacement,\phi)$, $\HighOrderTraction(\Displacement,\phi)$, $w(\Displacement,\phi)$ and $\HighOrderJump(\Displacement,\phi)$ in \eq\eqref{eq_FlexoForces} as the constrained variational principle in \eq\eqref{eq_variationalPrinciple}.

The penalty parameters $\beta_u,\beta_v,\beta_{C_u},\beta_\phi$ in \eq\eqref{Dirichlet} have to be chosen large enough, but too large values would lead to ill-conditioning problems when finding the equilibrium states $(\Displacement^*,\phi^*)$ numerically.
The derivation of stability lower bounds of the penalty parameters can be found in \ref{sec_app02}. However, moderate values of the penalty parameters are enough to ensure convergence and enforce boundary conditions properly \cite{Fernandez2004,Schillinger2016,badia2018,Ruberg2016}.

%%%%%
\subsection{Weak form of the boundary value problem}
%%%%%
The weak form of the unfitted variational formulation is presented next. 
Vanishing of the first variation of \eq\eqref{eq_Functionalimm} yields
\begin{multline}\label{firstvariation}
0=\var{\Pi}[\Displacement,\phi;\delta\Displacement,\delta\phi]=
\var{\Pi^\Omega}[\Displacement,\phi;\delta\Displacement,\delta\phi]
+
\var{\Pi^\text{Dirichlet}}[\Displacement,\phi;\delta\Displacement,\delta\phi]
\\+
\var{\Pi^\text{Neumann}}[\delta\Displacement,\delta\phi],
\quad
\forall(\delta\Displacement,\delta\phi)\in\mathcal{U}\otimes\mathcal{P};
\end{multline}
where
\begin{subequations}\begin{align}
\var{\Pi^\Omega}[\Displacement,\phi;\delta\Displacement,\delta\phi]&=
\int_{\Omega}\Big(\cauchyStress_{ij}(\Displacement,\phi)\strain_{ij}(\delta\Displacement)+\hyperStress_{ijk}(\Displacement,\phi)\strain_{ij,k}(\delta\Displacement)\nonumber\\&\qquad-\electricDisp_l(\Displacement,\phi)E_l(\delta\phi)-b_i\delta u_{i}+q\delta \phi\Big)\dd\Omega,
\\
\var{\Pi^\text{Dirichlet}}[\Displacement,\phi;\delta\Displacement,\delta\phi]&=
\int_{\partial\Omega_u}\Big(\beta_u\Big(u_i-\displacementKnown_i\Big)\delta u_i
-t_i(\Displacement,\phi)\delta u_i
-\Big(u_i-\displacementKnown_i\Big)t_i(\delta\Displacement,\delta\phi)\Big)\dd\Gamma
+{}\nonumber\\&{}+
\int_{\partial\Omega_v}\Big(\beta_v\Big(\partial^nu_i-\partialndisplacementKnown_i\Big)\partial^n\delta u_i
-r_i(\Displacement,\phi)\partial^n\delta u_i
-\Big(\partial^nu_i-\partialndisplacementKnown_i\Big)r_i(\delta\Displacement,\delta\phi)\Big)\dd\Gamma
+{}\nonumber\\&{}+
\int_{\partial\Omega_\phi}\Big(-\beta_\phi\Big(\phi-\phiKnown\Big)\delta\phi
+w(\Displacement,\phi)\delta\phi 
+\Big(\phi-\phiKnown\Big) w(\delta\Displacement,\delta\phi)\Big)\dd\Gamma
+{}\nonumber\\&{}+
\int_{C_u}\Big(\beta_{C_u}\Big(u_i-\displacementKnown_i\Big)\delta u_i
-j_i(\Displacement,\phi)\delta u_i 
-\Big(u_i-\displacementKnown_i\Big)j_i(\delta\Displacement,\delta\phi)\Big)\dd\ds,
\\
\var{\Pi^\text{Neumann}}[\delta\Displacement,\delta\phi]&=
\int_{\partial\Omega_t}-\tractionKnown_i\delta u_i\dd\Gamma
+
\int_{\partial\Omega_r}-\highOrderTractionKnown_i\partial^n\delta u_i\dd\Gamma
\nonumber\\&+
\int_{\partial\Omega_w}\wKnown\delta\phi\dd\Gamma
+
\int_{C_j}-\highOrderJumpKnown\delta u_i\dd\ds;
\end{align}\end{subequations}
being $\delta\Displacement\in\mathcal{U}$ and $\delta\phi\in\mathcal{P}$ admissible variations of $\Displacement$ and $\phi$, respectively.

The functionals $\var{\Pi^\Omega}$ and $\var{\Pi^\text{Dirichlet}}$ are conveniently rearranged as $\var{\Pi^\Omega}[\Displacement,\phi;\delta\Displacement,\delta\phi]=\var{\Pi^\Omega_\mathscr{B}}[\Displacement,\phi;\delta\Displacement,\delta\phi]-\var{\Pi^\Omega_\mathscr{L}}[\delta\Displacement,\delta\phi]$ and $\var{\Pi^\text{Dirichlet}}[\Displacement,\phi;\delta\Displacement,\delta\phi]=\var{\Pi^\text{Dirichlet}_\mathscr{B}}[\Displacement,\phi;\delta\Displacement,\delta\phi]-\var{\Pi^\text{Dirichlet}_\mathscr{L}}[\delta\Displacement,\delta\phi]$, with
\begin{subequations}\begin{align}
\var{\Pi^\Omega_\mathscr{B}}[\Displacement,\phi;\delta\Displacement,\delta\phi]&\coloneqq
\int_{\Omega}\Big(\cauchyStress_{ij}(\Displacement,\phi)\strain_{ij}(\delta\Displacement)+\hyperStress_{ijk}(\Displacement,\phi)\strain_{ij,k}(\delta\Displacement)-\electricDisp_l(\Displacement,\phi)E_l(\delta\phi)\Big)\dd\Omega,
\\
\var{\Pi^\Omega_\mathscr{L}}[\delta\Displacement,\delta\phi]&\coloneqq
\int_{\Omega}\Big(b_i\delta u_{i}-q\delta \phi\Big)\dd\Omega,
\\
\var{\Pi^\text{Dirichlet}_\mathscr{B}}[\Displacement,\phi;\delta\Displacement,\delta\phi]&\coloneqq
\int_{\partial\Omega_u}\Big(\Big(\beta_uu_i-t_i(\Displacement,\phi)\Big)\delta u_i
-u_it_i(\delta\Displacement,\delta\phi)\Big)\dd\Gamma
+{}\nonumber\\&{}+
\int_{\partial\Omega_v}\Big(\Big(\beta_v\partial^nu_i-r_i(\Displacement,\phi)\Big)\partial^n\delta u_i
-\partial^nu_ir_i(\delta\Displacement,\delta\phi)\Big)\dd\Gamma
+{}\nonumber\\&{}+
\int_{\partial\Omega_\phi}\Big(\Big(-\beta_\phi\phi+w(\Displacement,\phi)\Big)\delta\phi
+\phi w(\delta\Displacement,\delta\phi)\Big)\dd\Gamma
+{}\nonumber\\&{}+
\int_{C_u}\Big(\Big(\beta_{C_u}u_i-j_i(\Displacement,\phi)\Big)\delta u_i
-u_ij_i(\delta\Displacement,\delta\phi)\Big)\dd\ds,
\\
\var{\Pi^\text{Dirichlet}_\mathscr{L}}[\delta\Displacement,\delta\phi]&\coloneqq
\int_{\partial\Omega_u}\displacementKnown_i\Big(\beta_u\delta u_i-t_i(\delta\Displacement,\delta\phi)\Big)\dd\Gamma
+
\int_{\partial\Omega_v}\partialndisplacementKnown_i\Big(\beta_v\partial^n\delta u_i-r_i(\delta\Displacement,\delta\phi)\Big)\dd\Gamma
+{}\nonumber\\&{}+
\int_{\partial\Omega_\phi}\phiKnown\Big(-\beta_\phi\delta\phi+w(\delta\Displacement,\delta\phi)\Big)\dd\Gamma
+
\int_{C_u}\displacementKnown_i\Big(\beta_{C_u}\delta u_i-j_i(\delta\Displacement,\delta\phi)\Big)\dd\ds,
\end{align}\end{subequations}

The weak form of the unfitted formulation for flexoelectricity reads:
\begin{equation}\label{eq_FlexWeakSymmNitsche}
\emph{Find $(\Displacement,\phi)\in\mathcal{U}\otimes\mathcal{P}$ such that $\forall(\delta\Displacement,\delta\phi)\in\mathcal{U}\otimes\mathcal{P}:$}\quad
\mathscr{B}[\{\Displacement,\phi\},\{\delta\Displacement,\delta\phi\}]
=
\mathscr{L}[\{\delta\Displacement,\delta\phi\}]
;
\end{equation}
where
\begin{subequations}\begin{align}
\mathscr{B}[\{\Displacement,\phi\},\{\delta\Displacement,\delta\phi\}]&\coloneqq
\var{\Pi^\Omega_\mathscr{B}}[\Displacement,\phi;\delta\Displacement,\delta\phi]+\var{\Pi^\text{Dirichlet}_\mathscr{B}}[\Displacement,\phi;\delta\Displacement,\delta\phi],
\label{bil}\\
\mathscr{L}[\{\delta\Displacement,\delta\phi\}]&\coloneqq
\var{\Pi^\Omega_\mathscr{L}}[\delta\Displacement,\delta\phi]+\var{\Pi^\text{Dirichlet}_\mathscr{L}}[\delta\Displacement,\delta\phi]-\var{\Pi^\text{Neumann}}[\delta\Displacement,\delta\phi].
\end{align}\end{subequations}

%%%%%%%
\subsection{Formulation including sensing electrodes}\label{sec_03bis}
%%%%%%%
In electromechanics, conducting electrodes are frequently attached to the surface of the devices to enable either \emph{actuation} or \emph{sensing}. Actuators induce a deformation due to a prescribed electric potential, whereas sensors infer the deformation state by the measured change in the electric potential.
In both cases, as the electrodes are made of conducting material, the electric potential in the electrode is uniform.
The electrical Dirichlet boundary condition in \eq\eqref{elecBC} corresponds to actuating electrodes where the uniform electric potential is prescribed. In the case of sensing electrodes, the uniform electric potential is unknown and thus requires a special treatment as described next.

Let us consider a partition of the boundary distinguishing actuating and sensing electrodes, \ie
\begin{equation}
\partial\Omega=\partial\Omega_{\phi}\cup\partial\Omega_{\Phi}\cup\partial\Omega_{w},
\end{equation}
where $\partial\Omega_{\phi}$ and $\partial\Omega_{\Phi}$ correspond, respectively, to actuating and sensing electrodes on the boundary, respectively, and $\partial\Omega_{w}$ to the electrical Neumann boundary. The sensing boundary $\partial\Omega_{\Phi}$ is conformed by $N_\text{sensing}$ electrodes, namely
$
\partial\Omega_{\Phi}=\bigcup_{i=1}^{N_\text{sensing}}{\partial\Omega_{\Phi}}^i.
$

The electric potential $\phi$ on sensing electrodes is constant, but unknown. Thus, at each electrode ${\partial\Omega_{\Phi}}^i$ a new state variable ${\Phi}^i\in\mathbb{R}$ is introduced, which is a scalar denoting the unknown constant value of the electric potential. In other words,
\begin{equation}\label{eq_sensor1}
\phi-\Phi^i=0 \quad\text{on }{{\partial\Omega_{\Phi}}^i},\quad\forall i=1,\dots,N_\text{sensing}.
\end{equation}
Boundary conditions in \eq\eqref{eq_sensor1} are weakly enforced by adding to the energy potential $\Pi[\Displacement,\phi]$ in \eq\eqref{eq_Functionalimm} the work $\Pi^\text{Sensing}_i[\Displacement,\phi,\Phi^i]$ of each sensing electrode:
\begin{equation}\label{Sensor}
\Pi^\text{Sensing}_i[\Displacement,\phi,\Phi^i]=
\int_{\partial\Omega_\Phi^i}\Big(\phi-\Phi^i\Big) w(\Displacement,\phi)\dd\Gamma,
\end{equation}
and the associated variational principle for the equilibrium states $(\Displacement^*,\phi^*,\Phi^{1*},\dots,\Phi^{{N_\text{sensing}}*})$ of the body is
\begin{equation}\label{variationalPrincipleSensor}
(\Displacement^*,\phi^*,\Phi^{1*},\dots,\Phi^{{N_\text{sensing}}*})=\arg\min_{\Displacement\in\mathcal{U}}\max_{\phi\in\mathcal{P}}\min_{\Phi^1\in\mathbb{R}}\dots\min_{\Phi^{N_\text{sensing}}\in\mathbb{R}}\left(\Pi[\Displacement,\phi]+\sum_{i=1}^{N_\text{sensing}}\Pi^\text{Sensing}_i[\Displacement,\phi,\Phi^i]\right).
\end{equation}
Equation \eqref{Sensor} has a similar form to the Nitsche terms $\Pi^\text{Dirichlet}[\Displacement,\phi]$ in \eq\eqref{Dirichlet}, but the penalty term quadratic in \eq\eqref{eq_sensor1} is omitted here, because if it is positive in sign then it could be made arbitrarily large with respect to $\phi$,
and, conversely, if it is negative in sign then it could be made arbitrarily large with respect to $\Phi^i$.
Vanishing of the first variation of the energy functional in \eq\eqref{variationalPrincipleSensor} yields
\begin{equation}
0=\var{\Pi}[\Displacement,\phi;\delta\Displacement,\delta\phi]
+\sum_{i=1}^{N_\text{sensing}}\var{\Pi^\text{Sensing}_i}[\Displacement,\phi,\Phi^i;\delta\Displacement,\delta\phi,\delta\Phi^i],
\end{equation}
where
\begin{equation}
\var{\Pi^\text{Sensing}_i}[\Displacement,\phi,\Phi^i;\delta\Displacement,\delta\phi,\delta\Phi^i]\coloneqq
\int_{\partial\Omega_\Phi^i}\Big(\phi-\Phi^i\Big) w(\delta\Displacement,\delta\phi)\dd\Gamma+
\int_{\partial\Omega_\Phi^i}w(\Displacement,\phi)\Big(\delta\phi-\delta\Phi^i\Big)\dd\Gamma,
\end{equation}
being $\delta\Phi^i\in\mathbb{R}$ admissible variations of each $\Phi^i$.
Finally, the weak form of the unfitted formulation for flexoelectricity accounting for sensing electrodes reads:
\begin{multline}\label{eq_FlexWeakSymmNitscheSensor}
\emph{Find $(\Displacement,\phi,\Phi^1,...,\Phi^{N_\text{sensing}})\in\mathcal{U}\otimes\mathcal{P}\otimes\mathbb{R}^{N_\text{sensing}}$ such that}\\\emph{$\forall(\delta\Displacement,\delta\phi,\delta\Phi^1,...,\delta\Phi^{N_\text{sensing}})\in\mathcal{U}\otimes\mathcal{P}\otimes\mathbb{R}^{N_\text{sensing}}:$}\\
\mathscr{B}[\{\Displacement,\phi\},\{\delta\Displacement,\delta\phi\}]
+
\sum_{i=1}^{N_\text{sensing}}\var{\Pi^\text{Sensing}_i}[\Displacement,\phi,\Phi^i;\delta\Displacement,\delta\phi,\delta\Phi^i]
=
\mathscr{L}[\{\delta\Displacement,\delta\phi\}]
.\end{multline}
\section{Numerical approximation}\label{sec_04}
	% % %
\subsection{B-spline basis}\label{sec_BSplineBasis}
% % %
Fourth-order PDEs demand high-order continuity of the functional space for the numerical solution, \ie the displacement and  electric potential fields in the case of electromechanics.
Usually, $C^1$-continuity of the solution is enough; however, in the unfitted approach presented in this work, many boundary integrals in the Nitsche's weak forms (see Section \ref{sec_02}) involve \emph{third-order} derivatives, since the test functions do not vanish at the Dirichlet boundaries. It is clear that, in this case, $C^1$-continuous solutions are not smooth enough, and therefore we consider approximations of (at least) $C^2$-continuity.

A family of functions that provide high-order continuity is that of B-spline functions \cite{deBoor2001,Rogers2001,Piegl2012}, which are smooth piecewise polynomials. Being $p$ the polynomial degree, they are by construction $C^{p-1}$-continuous throughout the domain. Therefore, cubic ($p=3$) or higher-order B-spline basis are suitable for the numerical approximation of the formulation in Section \ref{sec_02}.

\begin{figure}[p]\centering
	\begin{subfigure}[t]{.33\textwidth}
		\captionof{figure}{Linear ($p=1$)}\vspace{-0.8em}
		\includegraphics[width=\textwidth]{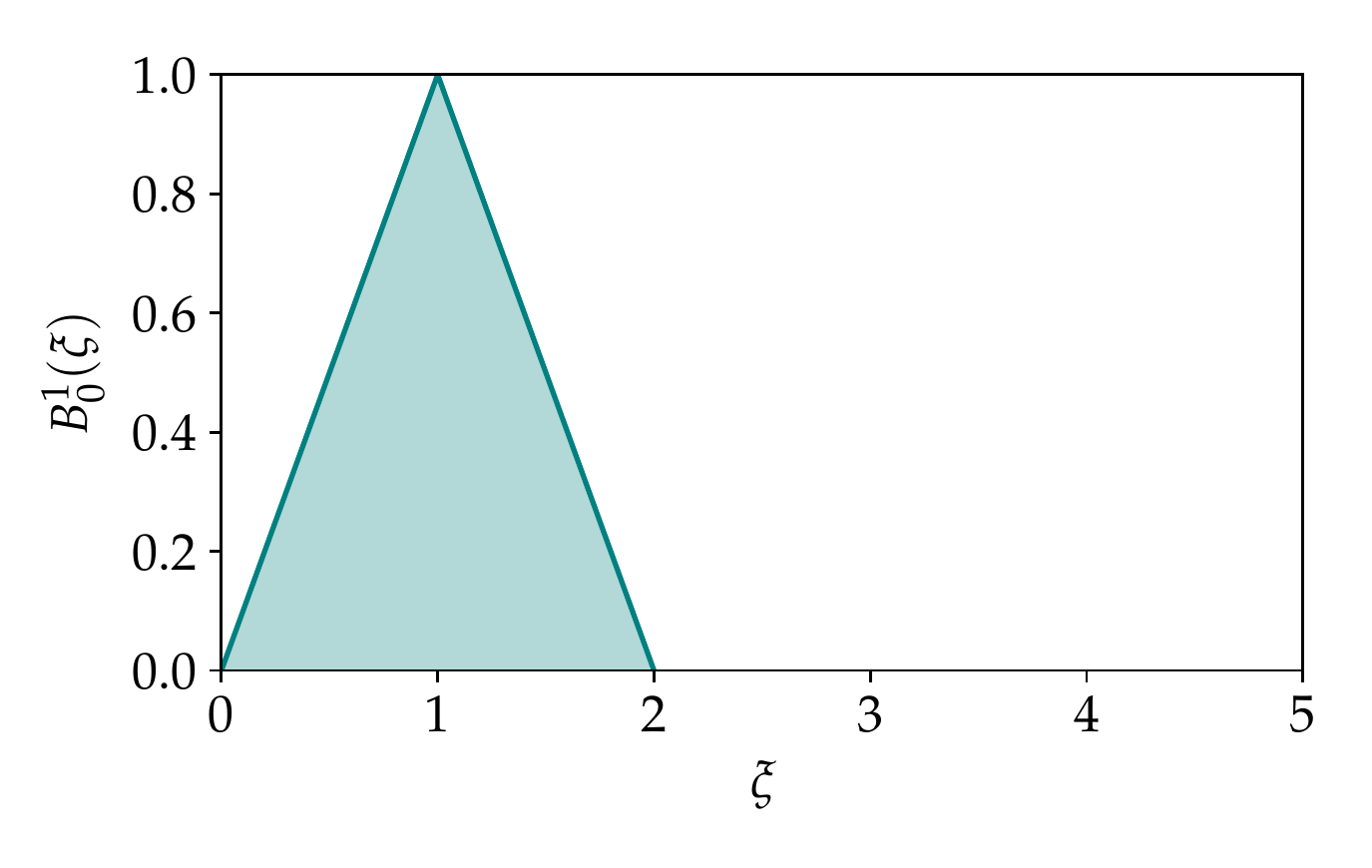}
		\vfill\label{fig_uni1}
	\end{subfigure}%
	\begin{subfigure}[t]{.33\textwidth}
		\captionof{figure}{Quadratic ($p=2$)}\vspace{-0.8em}
		\includegraphics[width=\textwidth]{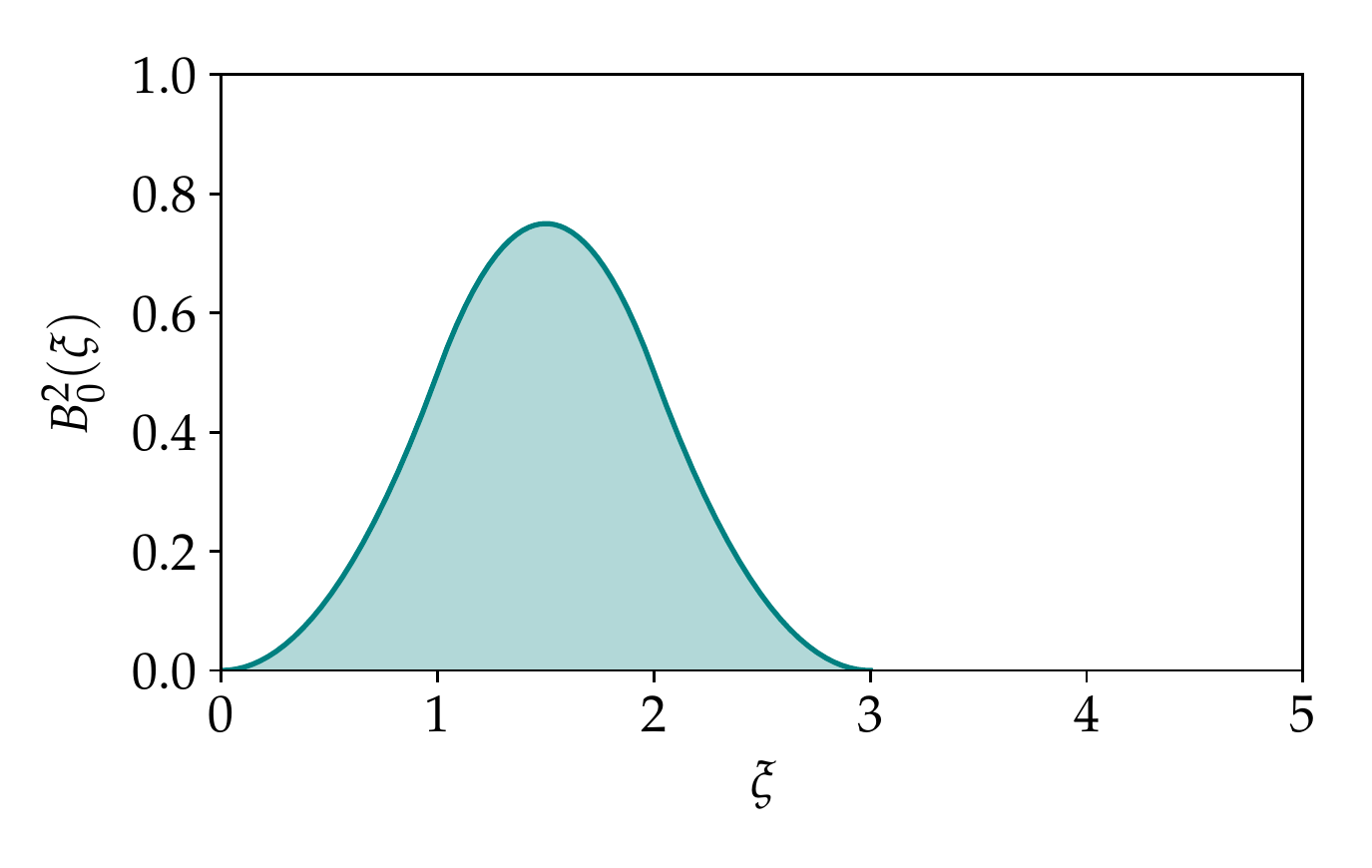}
		\vfill\label{fig_uni2}
	\end{subfigure}\\
	\begin{subfigure}[t]{.33\textwidth}
		\captionof{figure}{Cubic ($p=3$)}\vspace{-0.8em}
		\includegraphics[width=\textwidth]{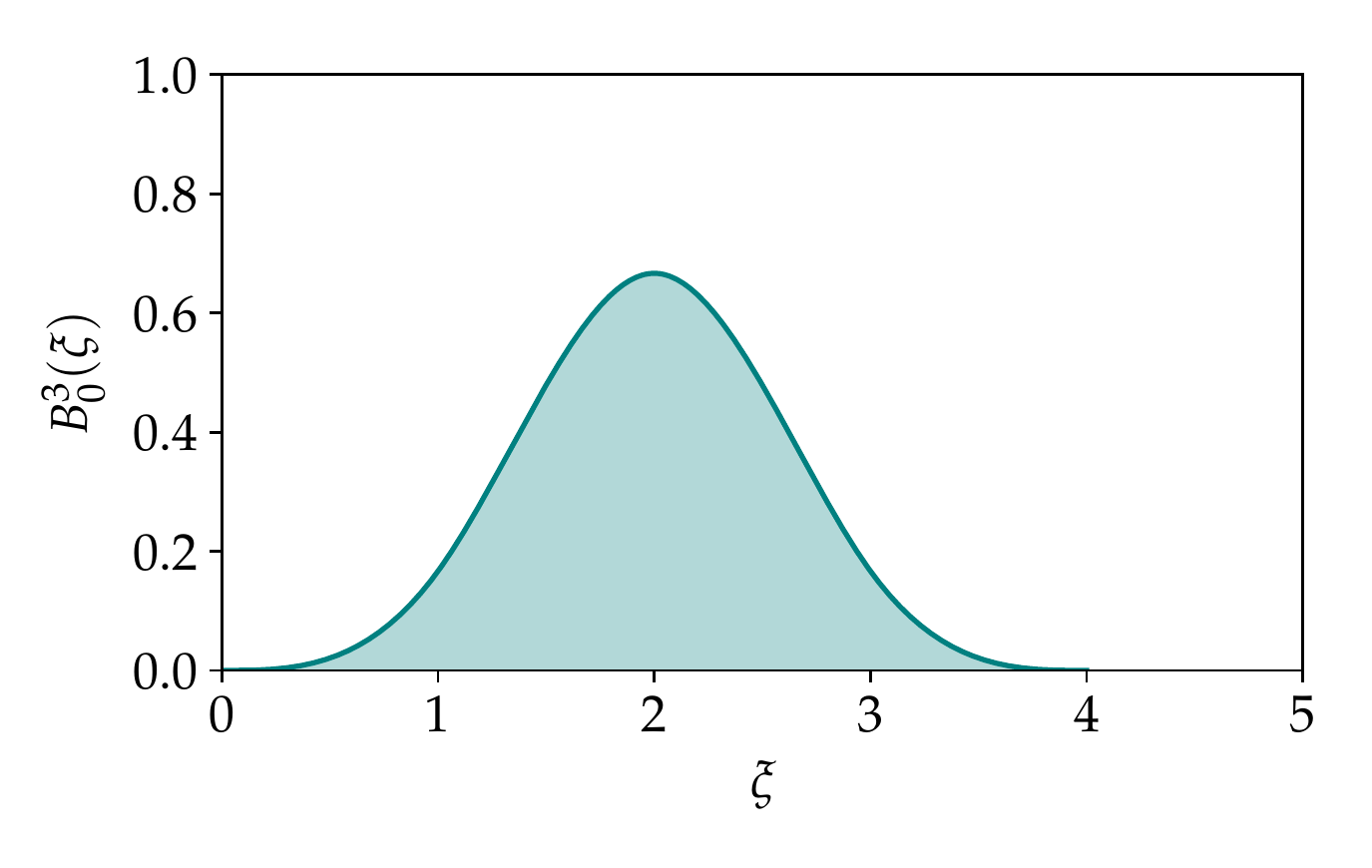}
		\vfill\label{fig_uni3}
	\end{subfigure}%
	\begin{subfigure}[t]{.33\textwidth}
		\captionof{figure}{Quartic ($p=4$)}\vspace{-0.8em}
		\includegraphics[width=\textwidth]{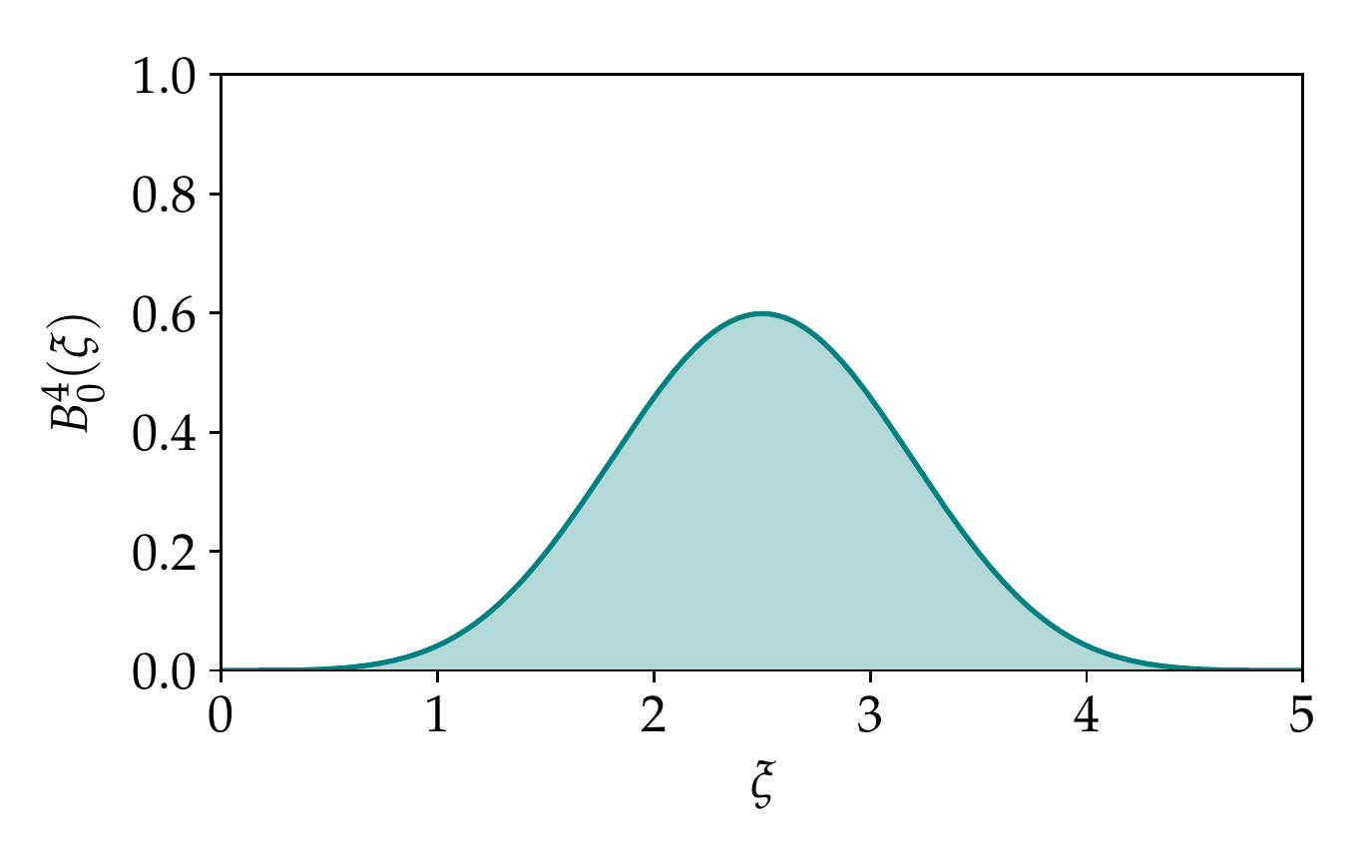}
		\vfill\label{fig_uni4}
	\end{subfigure}
	\vspace{-0.8em}\caption{First univariate B-spline basis function $B_0^p(\xi)$ of degree $p$.}
	\label{fig_BSplines}
\end{figure}

Let us consider a \emph{uniform} B-spline basis; local mesh refinement will be further considered in Section \ref{sec_adapt}. Without going into detail, the univariate uniform B-spline basis of degree $p$ consisting of $n$ basis functions is defined on the unidimensional parametric space $\xi\in[0,n+p]$ in terms of the \emph{uniform knot vector} $\{\xi_i\}=[0,1,2,3,\dots,n+p]$. The $i$-th function of this basis is defined recursively as \cite{deBoor2001}:
\begin{align}\label{eq_BSpline}
		B_i^0(\xi)&=
		\begin{cases}
			1 & \xi_i\leq\xi<\xi_{i+1}\\
			0 & \text{otherwise}
		\end{cases};
		\nonumber\\
		B_i^p(\xi)&=\frac{\xi-\xi_i}{\xi_{i+p}-\xi_i}B_i^{p-1}(\xi)+\frac{\xi_{i+p+1}-\xi}{\xi_{i+p+1}-\xi_{i+1}}B_{i+1}^{p-1}(\xi);
		\qquad
		i= 0,\dots,n-1.
	\end{align}
Due to the uniformity of the knot vector, the $i$-th B-spline function can be expressed as a translation of the first ($0$th) one as $B_i^p(\xi)=B_0^p(\xi-i)$. Figure \ref{fig_BSplines} shows the function $B_0^p(\xi)$ of the basis for degrees $p=\{1,\dots,4\}$.
\begin{figure}[p]\centering
	\begin{subfigure}[t]{.44\textwidth}
		\captionof{figure}{}
		\includegraphics[width=\textwidth]{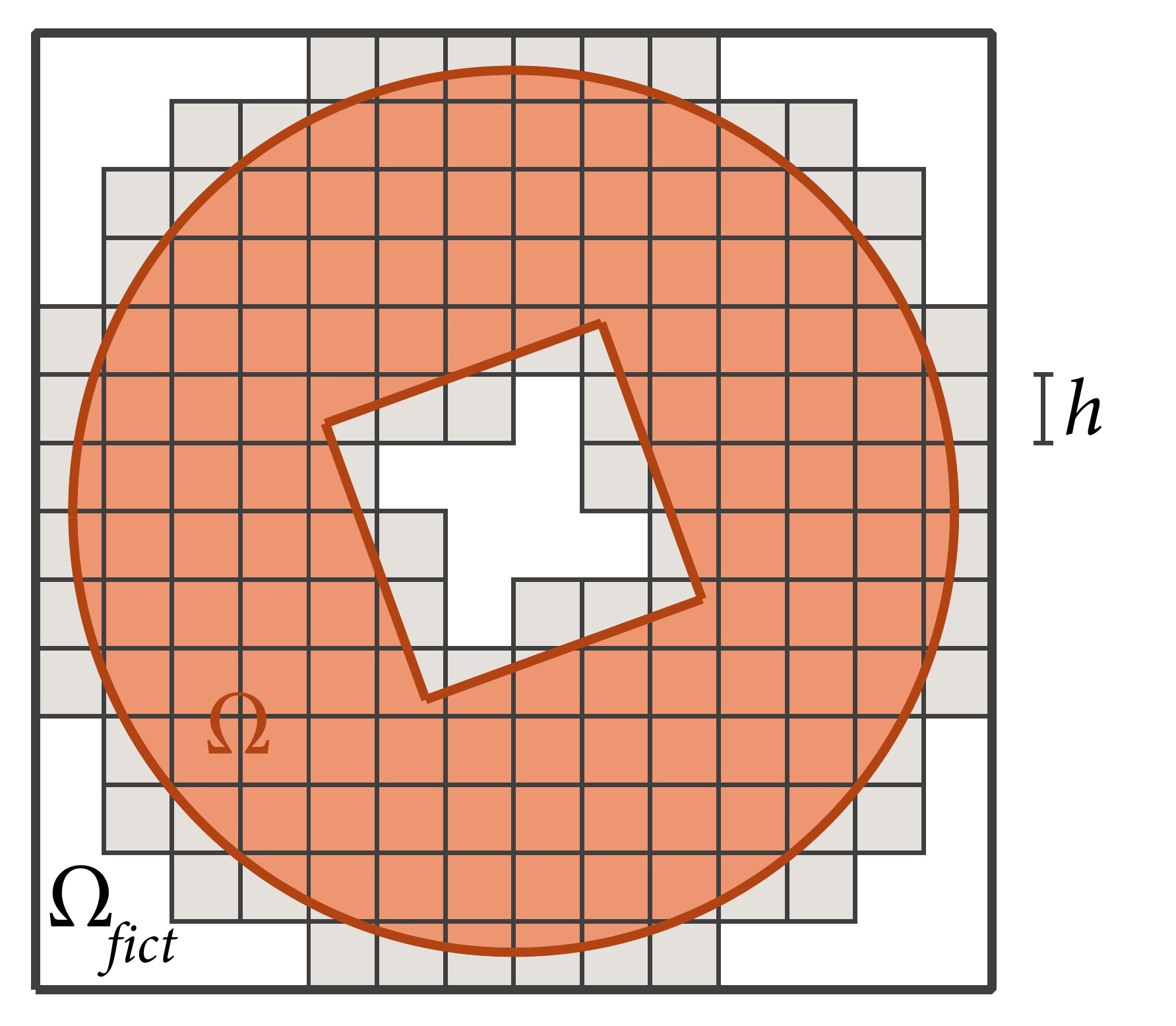}
		\label{fig_imm1}
	\end{subfigure}%
	\begin{subfigure}[t]{.155\textwidth}
		\captionof{figure}{}
		\includegraphics[width=\textwidth]{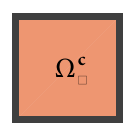}
		\label{fig_imm2}
		\\[-1em]
		\captionof{figure}{}
		\includegraphics[width=\textwidth]{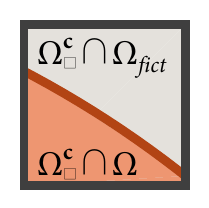}
		\label{fig_imm3}
	\end{subfigure}
	\caption{2D sketch of the immersed boundary method. \emph{a)} Physical domain $\Omega$ (red) immersed in the discretization of the embedding domain $\Omega_\square=\Omega\cup\Omega_\text{fict}$ (grey) with cell size $h$; outer cells are not depicted, \emph{b)} detail of an inner cell $\Omega_\square^{\toVect{c}}\in\mathcal{I}$ and \emph{c)} detail of a cut cell $\Omega_\square^{\toVect{c}}\in\mathcal{C}$. These cells are used to build the B-spline approximation spaces and for integration.}
	\label{fig_ImmersedMesh}
\end{figure}
In the three-dimensional space, the $\boldsymbol{i}$-th B-spline function $B_{\boldsymbol{i}}^p(\boldsymbol{\xi})$ of a trivariate B-spline basis (where $\boldsymbol{i}$ is the trivariate index $[i_{\xi},i_{\eta},i_{\tau}]$) is defined as the tensor product of three univariate B-spline functions as
\begin{multline}\label{eq_trivariate}
B_{\boldsymbol{i}}^p(\boldsymbol{\xi})=B_{[i_{\xi},i_{\eta},i_{\tau}]}^p([\xi,\eta,\tau])
\coloneqq
B_{i_{\xi}}^p(\xi) B_{i_{\eta}}^p(\eta) B_{i_{\tau}}^p(\tau);\\\emph{ with }\quad
i_{\xi}  = 0,\dots,n_{\xi}-1; \quad
i_{\eta} = 0,\dots,n_{\eta}-1;\quad
i_{\tau} = 0,\dots,n_{\tau}-1,
\end{multline}
which is defined on the three-dimensional parametric space $\boldsymbol{\xi}\in[0,n_{\xi}+p]\otimes[0,n_{\eta}+p]\otimes[0,n_{\tau}+p]$.
Therefore, the parametric space is a cuboid which is defined globally on a Cartesian grid,
in contrast with traditional Lagrangian basis present in standard FEM implementations, whose parametric space is defined elementwise.

In the physical space, the problem unknowns $\Displacement$ and $\phi$ are approximated as
\begin{subequations}\label{eq_approx}\begin{align}
	\left[\Displacement(\boldsymbol{x})\right]_d\simeq N_{\toVect{i}}(\boldsymbol{x}){a^u}_{\toVect{i}d}
	=&[N_{\toVect{i}}\circ\boldsymbol{\varphi}](\boldsymbol{\xi}){a^u}_{\toVect{i}d}
	=B_{\toVect{i}}^p(\boldsymbol{\xi}){a^u}_{\toVect{i}d},\qquad d=1,2,3;
	\\
	\phi(\boldsymbol{x})\simeq N_{\toVect{i}}(\boldsymbol{x}){a^\phi}_{\toVect{i}}
	=&[N_{\toVect{i}}\circ\boldsymbol{\varphi}](\boldsymbol{\xi}){a^\phi}_{\toVect{i}}
	=B_{\toVect{i}}^p(\boldsymbol{\xi}){a^\phi}_{\toVect{i}};
\end{align}\end{subequations}
where $\{\boldsymbol{{a^u}},\boldsymbol{{a^\phi}}\}$ are the degrees of freedom of the numerical solution (known as the \emph{control variables} in B-spline nomenclature), $\boldsymbol{N}=[\boldsymbol{B}^p\circ\boldsymbol{\varphi}^{-1}]$ are the basis functions at the physical space, and $\boldsymbol{\varphi}(\boldsymbol{\xi})$ is the \emph{geometric mapping}, which is a bijection that maps a given point $\boldsymbol{\xi}$ in the parametric space to a given point $\boldsymbol{x}$ in the physical space.
Typically, the map $\boldsymbol{\varphi}(\boldsymbol{\xi})$ is expressed as the interpolation of a discretization of the physical space, namely:
\begin{equation}
	\left[\boldsymbol{\varphi}(\boldsymbol{\xi})\right]_d\simeq S_{\toVect{i}}(\boldsymbol{\xi}){\hat{x}}_{{\toVect{i}}d},\qquad d=1,2,3;
\end{equation}
where $\boldsymbol{S}(\boldsymbol{\xi})$ are the basis functions for the interpolation of the geometry, and $\boldsymbol{\hat{x}}$ are points on the physical space defining the map (known as the \emph{control points} in B-spline nomenclature).

Different choices of $\boldsymbol{S}(\boldsymbol{\xi})$ and $\boldsymbol{\hat{x}}$ are possible. However, since we want $\boldsymbol{N}(\boldsymbol{x})$ to be $C^{p-1}$-continuous, $\boldsymbol{S}(\boldsymbol{\xi})$ has to be $C^{p-1}$-continuous too, and the most natural choice is $\boldsymbol{S}(\boldsymbol{\xi})\coloneqq\boldsymbol{B}^p(\boldsymbol{\xi})$.
Therefore, the map $\boldsymbol{x}=\boldsymbol{\varphi}(\boldsymbol{\xi})$ is defined globally. This fact hinders a \emph{conforming} discretization of the physical space, since it requires an underlying rigid, Cartesian-like mesh in order to be mapped to the parametric space (as done in Isogeometric Analysis \cite{Hughes2005} and related works).
In order to circumvent this requirement on the discretization of the physical space, we consider a different approach where the parametric space of the B-spline basis is not mapped to a conforming discretization of the physical space, but rather to a \emph{non-conforming} one, naturally providing high-order continuity of the spanned functional space on arbitrary geometries. This concept is known as the immersed boundary approach and is introduced next.
% % %
\subsection{Immersed boundary method}\label{sec_imm}
% % %

The main idea of the \emph{immersed boundary} method, also known as the embedded domain method \cite{Peskin2002,Mittal2005}, is to extend the physical domain $\Omega$ to a larger \emph{embedding domain} $\Omega_\square=\Omega\cup\Omega_\text{fict}$ which is the one to be discretized, \ie $\Omega_\square=\bigcup_{\toVect{c}}\Omega^{\toVect{c}}_\square$ (see \fig\ref{fig_imm1}).
In this way, the discretization onto the physical space does not depend on the physical domain ($\Omega$) shape.
In order to combine this approach with a B-spline basis, the embedding domain $\Omega_\square$ is defined as a cuboid, and discretized using a Cartesian-like mesh.

The physical boundary $\partial\Omega$ is allowed to intersect the cells $\Omega_\square^{\toVect{c}}$ of the embedding mesh arbitrarily. Cells are classified into three different sets $\mathcal{I}$, $\mathcal{C}$ and $\mathcal{O}$, depending on their intersection with the physical domain $\Omega$:
\begin{enumerate}[label=\emph{\roman*)}]
	\item $\mathcal{I}\coloneqq\{\Omega_\square^{\toVect{c}}:\Omega_\square^{\toVect{c}}\subseteq\Omega\}$, the set of inner cells which remain uncut within the domain (\fig\ref{fig_imm2}),
	\item $\mathcal{C}\coloneqq\{\Omega_\square^{\toVect{c}}:\Omega_\square^{\toVect{c}}\nsubseteq\Omega\text{ and }\Omega_\square^{\toVect{c}}\cap\Omega\neq\emptyset\}$, the set of cells cut by the boundary (\fig\ref{fig_imm3}),
	 \item $\mathcal{O}\coloneqq\{\Omega_\square^{\toVect{c}}:\Omega_\square^{\toVect{c}}\cap\Omega=\emptyset\}$, the set of outer cells, which are neglected.
\end{enumerate}

For the sake of convenience, and without loss of generality, in this work we consider \emph{structured cubic Cartesian} meshes, which present several practical advantages. On the one hand, all cells are cubes, 
and a linear mapping can be considered in each cell, namely $\boldsymbol{\varphi}(\boldsymbol{\xi})=\boldsymbol{\breve{x}}_{\toVect{1}}+h\boldsymbol{\xi}$, being $h$ the physical cell size and $\boldsymbol{\breve{x}}_{\toVect{1}}$ the first corner of the cell.
Thus, the Jacobian of the geometric mapping is constant, namely
$\boldsymbol{J}(\boldsymbol{\xi})=\gradient_\xi\boldsymbol{\varphi}(\boldsymbol{\xi})=h\toMat{I}_3$, where $\toMat{I}_3$ is the identity matrix of rank 3.
On the other hand, in a \emph{linear} problem all the inner cells of the same size lead to the same elemental matrix, which is computed just once.

At this point, the geometry of $\Omega$ is \emph{only} used for cell classification. Without going into details, this is usually accomplished by checking whether all vertices of each cell lie within the domain (inner cell), only part of them (cut cell) or none of them (outer cell).
In the case of implicit boundary representation (\eg level set approaches) it is enough to evaluate the level set function on the vertices of each cell (see \cite{Fries2015,Fries2016,Kudela2016,Legrain2012}). For explicit boundary representation (\eg CAD descriptions), this can be achieved by ray-tracing procedures (see \cite{Marco2015,Marco2017}). In this work we restrict ourselves to explicit boundary representation by means of NURBS surfaces in 3D and NURBS curves in 2D.

% % %
\subsection{Integration on cut cells}
% % %
Bulk integrals are numerically performed in each cell, \ie inner ones $\Omega_\square^{\toVect{c}}\in\mathcal{I}$ and also the physical part of cut ones $\Omega_\square^{\toVect{c}}\cap\Omega$, for $\Omega_\square^{\toVect{c}}\in\mathcal{C}$ (see \fig\ref{fig_imm2} and \ref{fig_imm3}). Standard cubature rules \cite{Vincent2015} apply for the former, but not for the latter which can have arbitrary shape.
To this end, the physical part $\Omega_\square^{\toVect{c}}\cap\Omega$ of every cut cell $\Omega_\square^{\toVect{c}}\in\mathcal{C}$ is divided into several sub-domains (\eg cuboids or tetrahedra) which are easily integrated. To sub-divide cut cells we rely on the \emph{marching cubes algorithm} \cite{Lorensen1987}, which splits each cell into several conforming tetrahedra, although other conforming \cite{Kudela2016,Fries2016} or non-conforming \cite{Duster2008,Schillinger2015} subdivision schemes are also possible.
See \cite{Marco2015} for details of our current implementation.
Surface and line integrals are similarly performed on each corresponding sub-domain boundary.

Note that integration sub-domains in contact with the physical domain boundary $\partial\Omega$ might have curved faces or edges in the case $\partial\Omega$ is not flat. Hence, a linear cellwise approximation of the geometry leads to a geometric error of order $2$ which might spoil the optimal convergence of the method. Therefore, cell-wise polynomial approximations of the geometry of degree $p$ are required in general. Alternatively, we exploit the explicit NURBS representation of the geometry by resorting to the NEFEM approach \cite{Sevilla2008,Sevilla2011Seamless,Sevilla2011Numerical,Sevilla20113D,Legrain2013} which captures the \emph{exact} geometry without the need of any polynomial approximation \cite{Marco2015}.

% % %
\begin{remark}\label{stab}
% % %
The discretization of the weak forms in \eq\eqref{eq_FlexWeakSymmNitsche} and \eqref{eq_FlexWeakSymmNitscheSensor}, with B-spline basis functions and the mentioned numerical integration, leads to a linear system of equations for the coefficients of the approximation of the unknowns $\{\Displacement,\phi\}$, namely $\{\boldsymbol{{a^u}},\boldsymbol{{a^\phi}}\}$.
This linear system typically suffers from ill-conditioning in the presence of cut cells with a small portion in the domain, \ie when $|\Omega_\square^c\cap\Omega|\ll|\Omega_\square^c|$ for a given cell.
Ill-conditioning arises basically due to: \emph{i)} basis functions on the trimmed cell having very small contribution to the integral terms, and \emph{ii)} basis functions being quasi-linearly dependent on the trimmed cell \cite{dePrenter2016}. Moreover, ill-conditioning is more severe for high-order basis \cite{dePrenter2016}.
A detailed investigation on ill-conditioning of immersed boundary methods can be found in \cite{dePrenter2016}.

Several strategies have been proposed to alleviate ill-conditioning of trimmed cells, such as
the ghost penalty method \cite{Burman2010b},
the artificial stiffness approach \cite{Schillinger2015,Duster2008}, 
the extended B-spline method \cite{Hollig2001,Hollig2012,Ruberg2012,Ruberg2016} or
special preconditioning techniques specifically designed for immersed boundary methods \cite{dePrenter2016}, among others.

For uniform meshes, the extended B-spline approach by H\"{o}llig et al. \cite{Hollig2001,Hollig2012,Ruberg2012,Ruberg2016} is considered, due to its simple form and good performance. The main idea is to express the critical basis functions on the boundary as linear combinations of inner ones. The constrained basis has less degrees of freedom, but the conditioning and approximation properties are equivalent to those of body-fitted methods \cite{Hollig2001}.
The extension to hierarchical meshes (see Section \ref{sec_adapt}) follows the same idea but involves a more sophisticated implementation. In the numerical tests, for the sake of simplicity, hierarchical meshes are stabilized by means of a simple diagonal scaling  preconditioning.
\end{remark}

\subsection{Local mesh refinement: Hierarchical B-spline basis}\label{sec_adapt}
% % %
\begin{figure}[p]\centering
	\begin{subfigure}{.5\textwidth}\centering
		\vspace{1em}\captionof{figure}{Linear \emph{($p=1$)}. $\quad\toVect{s}^1=\frac{1}{2}[1,2,1]$}\vspace{-0.5em}
		\includegraphics[height=8em]{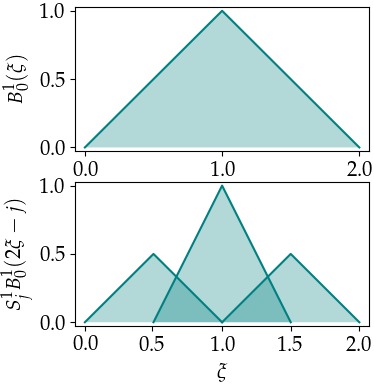}
		\label{fig_ref1}
	\end{subfigure}%
	\begin{subfigure}{.5\textwidth}\centering
		\vspace{1em}\captionof{figure}{Quadratic \emph{($p=2$). $\quad\toVect{s}^2=\frac{1}{4}[1,3,3,1]$}}\vspace{-0.5em}
		\includegraphics[height=8em]{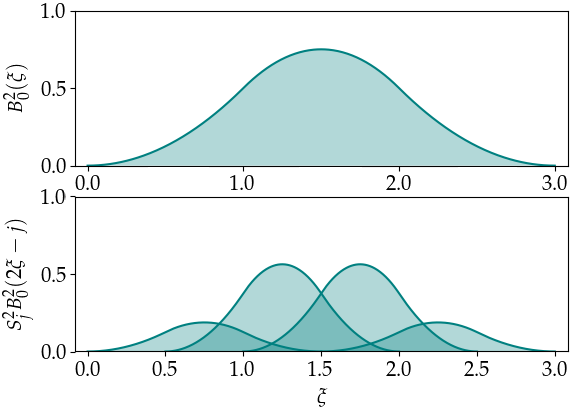}
		\label{fig_ref2}
	\end{subfigure}
	\newline\vfill
	\begin{subfigure}{.5\textwidth}\centering
		\vspace{1em}\captionof{figure}{Cubic \emph{($p=3$). $\quad\toVect{s}^3=\frac{1}{8}[1,4,6,4,1]$}}\vspace{-0.5em}
		\includegraphics[height=8em]{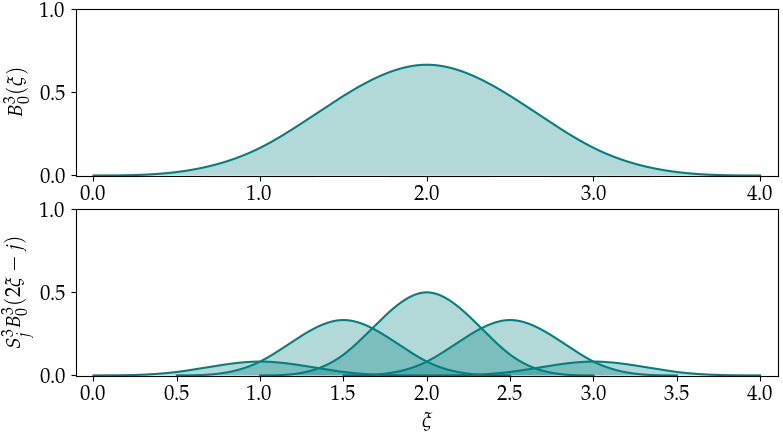}
		\label{fig_ref3}
	\end{subfigure}%
	\begin{subfigure}{.5\textwidth}\centering
		\vspace{1em}\captionof{figure}{Quartic \emph{($p=4$). $\quad\toVect{s}^4=\frac{1}{16}[1,5,10,10,5,1]$}}\vspace{-0.5em}
		\includegraphics[height=8em]{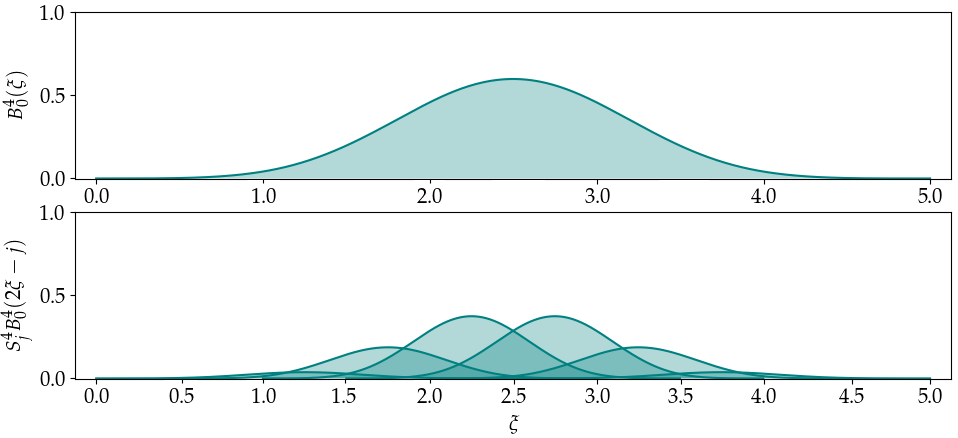}
		\label{fig_ref4}
	\end{subfigure}
	\caption{Hierarchical refinement of univariate B-spline basis function $B_0^p(\xi)$ of degree $p$. \emph{Top:} Original (parent) B-spline. \emph{Bottom:} The $j$-th children B-spline basis, $j=\{0,\dots,p+1\}$, that arise from the two-scale relation.}
	\label{fig_hierarchy}
\end{figure}

\begin{figure}[p]\centering
	\begin{subfigure}{.47\textwidth}\centering
		\vspace{1em}\captionof{figure}{}\vspace{-0.5em}
		\includegraphics[width=\linewidth]{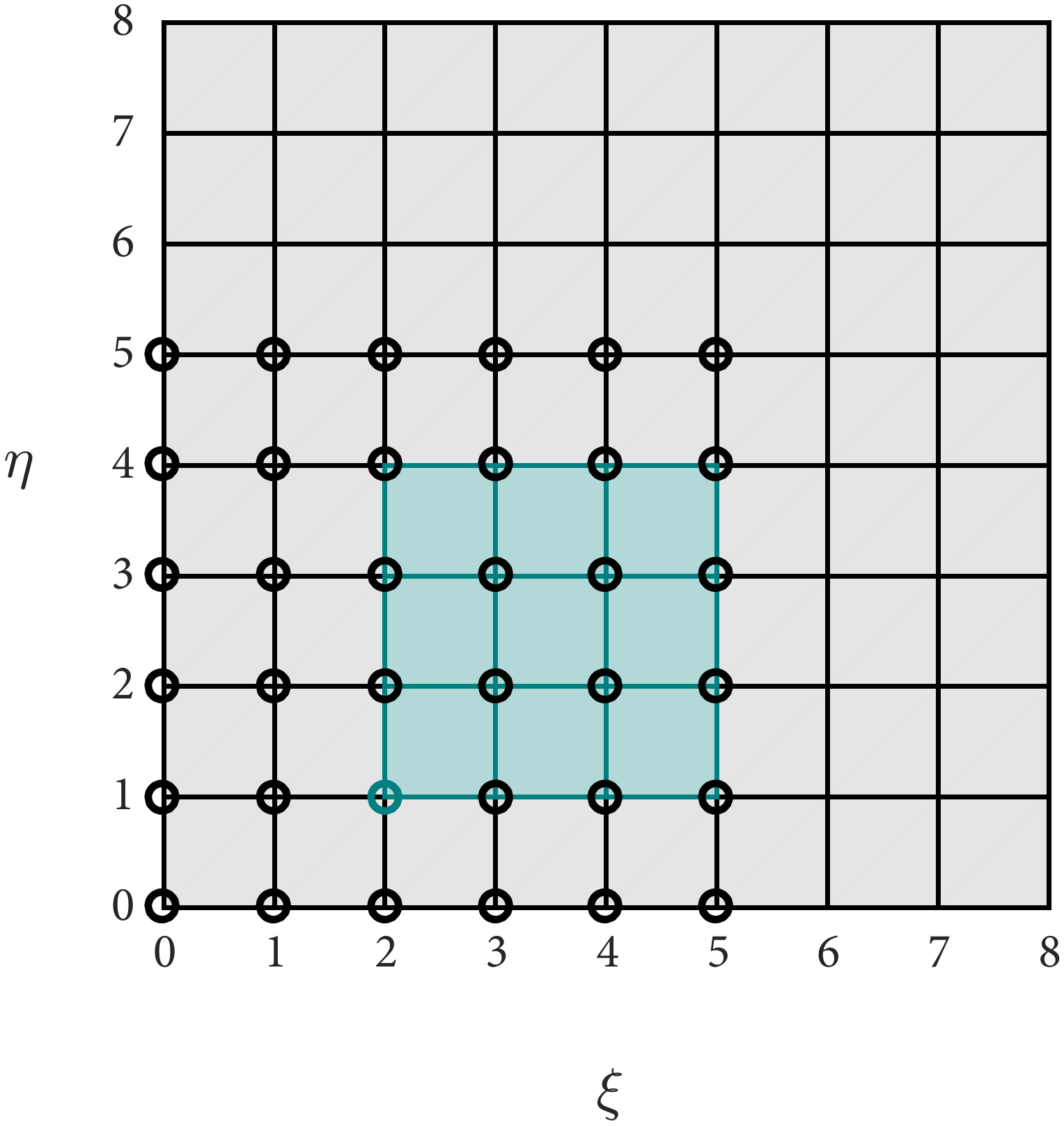}
		\label{fig_href1}
	\end{subfigure}\hfill
	\begin{subfigure}{.47\textwidth}\centering
		\vspace{1em}\captionof{figure}{}\vspace{-0.5em}
		\includegraphics[width=\linewidth]{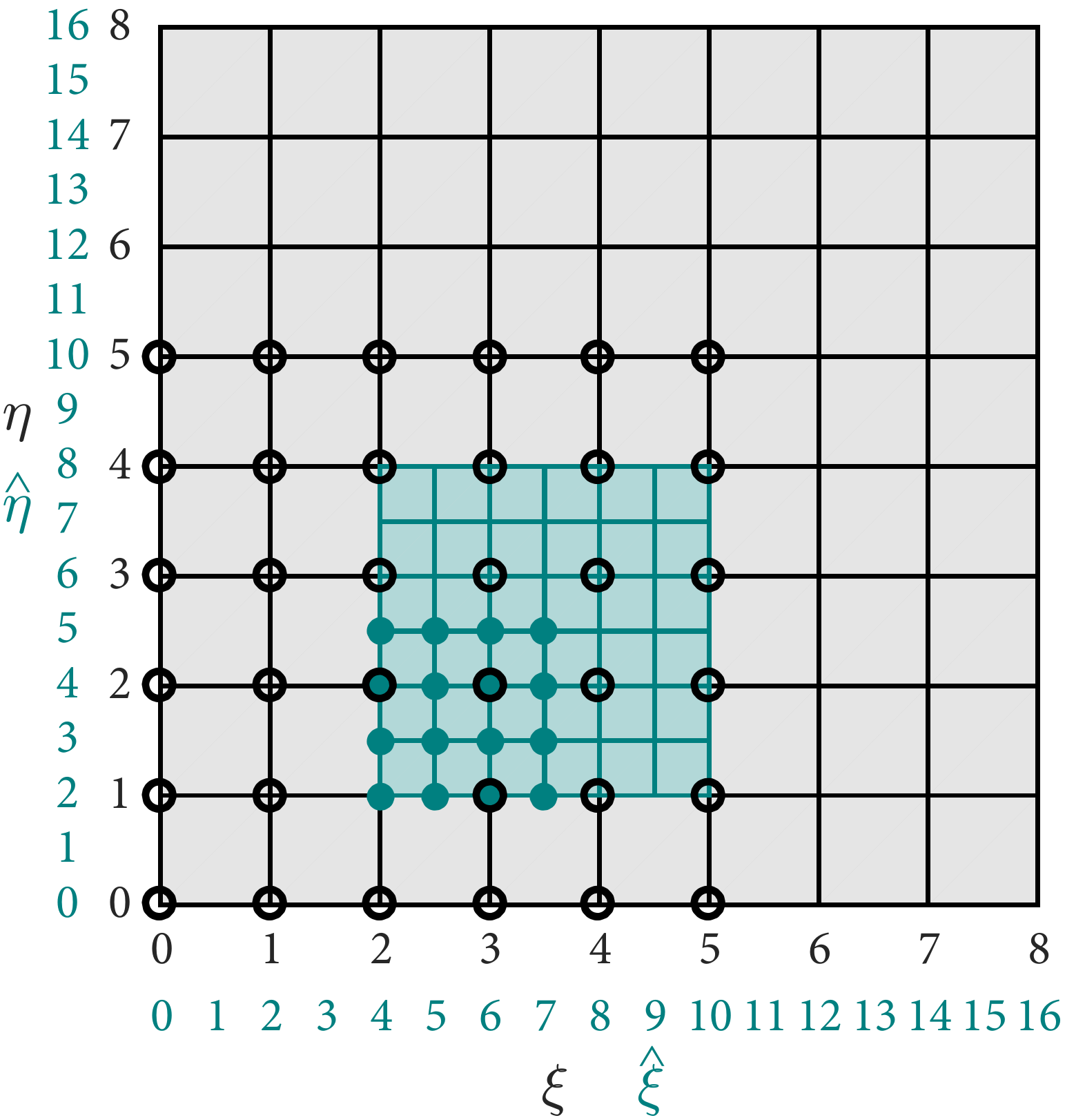}
		\label{fig_href2}
	\end{subfigure}
	\caption{Hierarchical refinement of a quadratic \emph{(p=2)} bivariate B-spline basis. \emph{a)} Uniform mesh; B-spline basis function $B_{2,1}^2(\xi,\eta)$ (blue) is selected for refinement. \emph{b)} Hierarchical mesh; basis function $B_{2,1}^2(\xi,\eta)$ is replaced by their 16 children $\hat{B}_{i,j}^2(\hat\xi,\hat\eta)$ (blue), $\forall \{i,j\}=\{4,5,6,7\}\otimes\{2,3,4,5\}$.}
	\label{fig_hierarchy2}
\end{figure}

Hierarchical B-spline refinement was first introduced by Forsey and Bartels \cite{Forsey1988}. It can be understood as a technique for locally enriching the approximation space by replacing selected coarse B-splines (parents) with finer ones (children). It is based on a remarkable property of uniform B-splines: their natural refinement by subdivision. For a univariate B-spline basis of degree $p$, the subdivision property leads to the following \emph{two-scale relation} \cite{Zorin2000}:
\begin{equation}\label{eq_twoscale}
B_i^p(\xi)= \sum_{j=0}^{p+1}s_j^pB_i^p(2\xi-j)\coloneqq \sum_{j=2i}^{2i+p+1}\hat{B}_j^p(\hat\xi),
\quad\text{with}\quad s_j^p=\frac{1}{2^p}\binom{p+1}{j}=\frac{2^{-p}(p+1)!}{j!(p+1-j)!};
\end{equation}
where $\hat\xi\coloneqq2\xi$.

In other words, a B-spline function $B_i^p(\xi)$ can be expressed as a linear combination of contracted, translated and scaled copies $\hat{B}_j^p(\hat\xi)$ of itself \cite{Schillinger2012}, as illustrated in \fig\ref{fig_hierarchy} for B-splines of different polynomial degree $p$. The extension to higher dimensions is trivial by means of the tensor product of univariate bases.

Without going into details, a hierarchical B-spline basis is defined from a uniform B-spline basis by replacing some basis functions with their corresponding children (see \fig\ref{fig_hierarchy2}). This process can be performed recursively, leading to a parent-children hierarchy spanning several levels of refinement. Since each basis function spans several cells, basis refinement implies refinement of multiple cells. The change of focus from \emph{element} refinement (as in conventional FE) to \emph{basis} refinement is the key point, which allows maintaining the smoothness of the functional space. Further details can be found in \cite{Bornemann2013,Kraft1995,Kraft1997,Schillinger2012,Vuong2011} and references therein.

At the implementation level, the elemental matrices of inner cells can be computed just once \emph{per level of hierarchy} by means of the \emph{subdivision projection technique} developed in \cite{Bornemann2013}. Therefore, hierarchical B-spline bases maintain the computational benefits of uniform meshes, as explained in Section \ref{sec_imm}, while allowing local mesh refinement.
\section{Numerical results}\label{sec_05}
Several numerical simulations are presented next to illustrate the performance of the method.
The first example shows the effect of disregarding edge boundary conditions. A synthetic polynomial solution is considered, which can be exactly captured only if edge boundary conditions are enforced.
In the second example we perform a sensitivity analysis of the solution and the condition number of the system matrix with respect to the Nitsche penalty parameters on the Dirichlet boundaries.
The third example consists on an error convergence analysis considering a non-trivial 2D geometry with curved boundaries and corners, where optimal convergence rates are achieved for different approximation degrees.
The fourth and fifth examples deal with two typical setups for flexoelectric characterization, namely a cantilever beam \cite{Majdoub2008} and a truncated pyramid \cite{Cross2006}. We compare our simulation results to previous solutions obtained in our group with the maximum entropy meshless method \cite{Abdollahi2014} and with approximate analytical solutions from the literature.
In the sixth example we present a 3D simulation of a rod with varying semi-circular cross section under torsion, which could be used to measure the shear flexoelectric coefficient \cite{MocciShearUnpublished}.

The material tensors in this section are defined next.
The mechanics are described by an isotropic elasticity model, with a Young modulus $E$ and a poisson ratio $\nu$, and enriched with an isotropic strain-gradient elasticity model depending on a single length scale parameter $l$. In the 2D case, plane strain conditions are assumed.
Electrostatics are described by an isotropic model with dielectricity constant $\kappa_L$.
Piezoelectricity is described by a tetragonal symmetry model oriented in a certain principal direction $\toVect{d_\text{piezo}}$. It depends on the longitudinal, transversal and shear piezoelectric coefficients $e_L$, $e_T$ and $e_S$, respectively.
Flexoelectricity is described by a cubic symmetry model oriented in the Cartesian axes, with longitudinal, transversal and shear flexoelectric coefficients $\mu_L$, $\mu_T$ and $\mu_S$, respectively.
The complete form of every material tensor can be found in \ref{sec_app03}.

Finally, we briefly comment on the choice of penalty parameters in \eq\eqref{Dirichlet}. \ref{sec_app02} shows the derivation of theoretical stability lower bounds. However, for the sake of simplicity, in the following numerical examples, and analogously to other works in the literature \cite{Fernandez2004,Schillinger2016,badia2018,Ruberg2016}, we consider the penalty parameters $\beta_u,\beta_v,\beta_{C_u},\beta_\phi\in\mathds{R}^+$ in terms of a \emph{dimensionless} parameter $\zeta\in\mathds{R}^+$ as follows:
\begin{equation}\label{parameterzeta}
\beta_u    \coloneqq \frac{E       }{h  }\zeta;\qquad
\beta_v    \coloneqq \frac{l^2E    }{h  }\zeta;\qquad
\beta_{C_u}\coloneqq \frac{l^2E    }{h^2}\zeta;\qquad
\beta_\phi \coloneqq \frac{\kappa_L}{h  }\zeta;
\end{equation}
where $h$ denotes the physical cell size. We choose $\zeta=100$ for all numerical examples in this work, which is a suitable value as confirmed by the sensitivity analysis in Section \ref{SensitivityAnalysis}.
% % % % %
\subsection{Effect of non-local corner conditions}\label{sec_edgebc}
% % % % %
\begin{figure}[b!]
	\begin{subfigure}[t]{0.25\textwidth}\centering
		\vspace{1em}\captionof{figure}{}\hspace{1.05em}
		\includegraphics[width=7em]{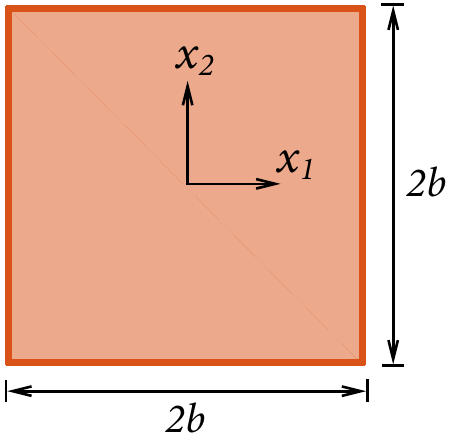}
		\label{fig_ex_corner_setup}
	\end{subfigure}%
	\begin{subfigure}[t]{0.35\textwidth}\centering
		\vspace{1em}\captionof{figure}{}
		\includegraphics[height=15em]{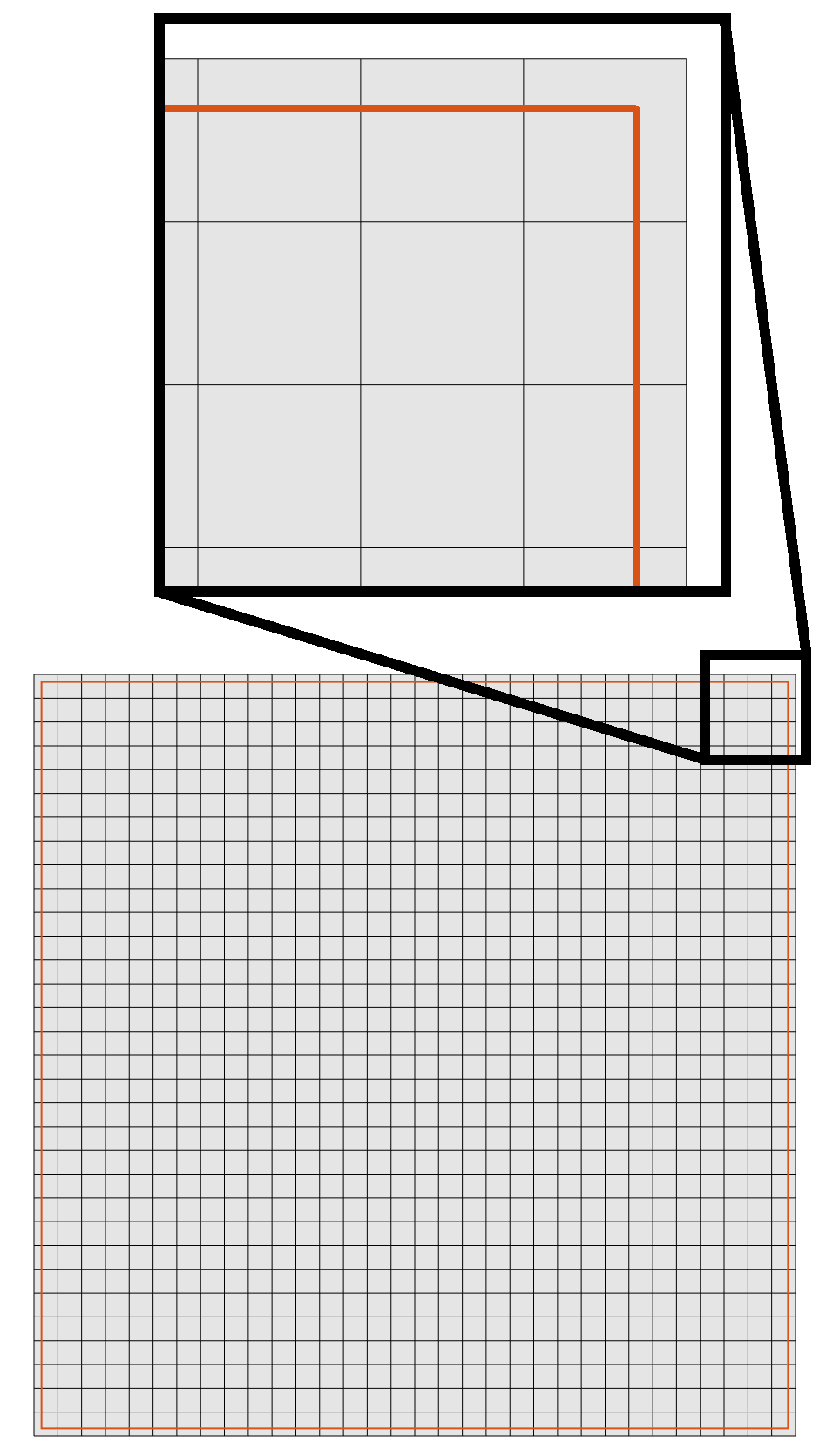}
		\label{fig_ex_corner_mesh}
	\end{subfigure}%
	\begin{subfigure}[t]{0.4\textwidth}\centering
        \vspace{1em}\captionof{figure}{}
		\includegraphics[height=15em]{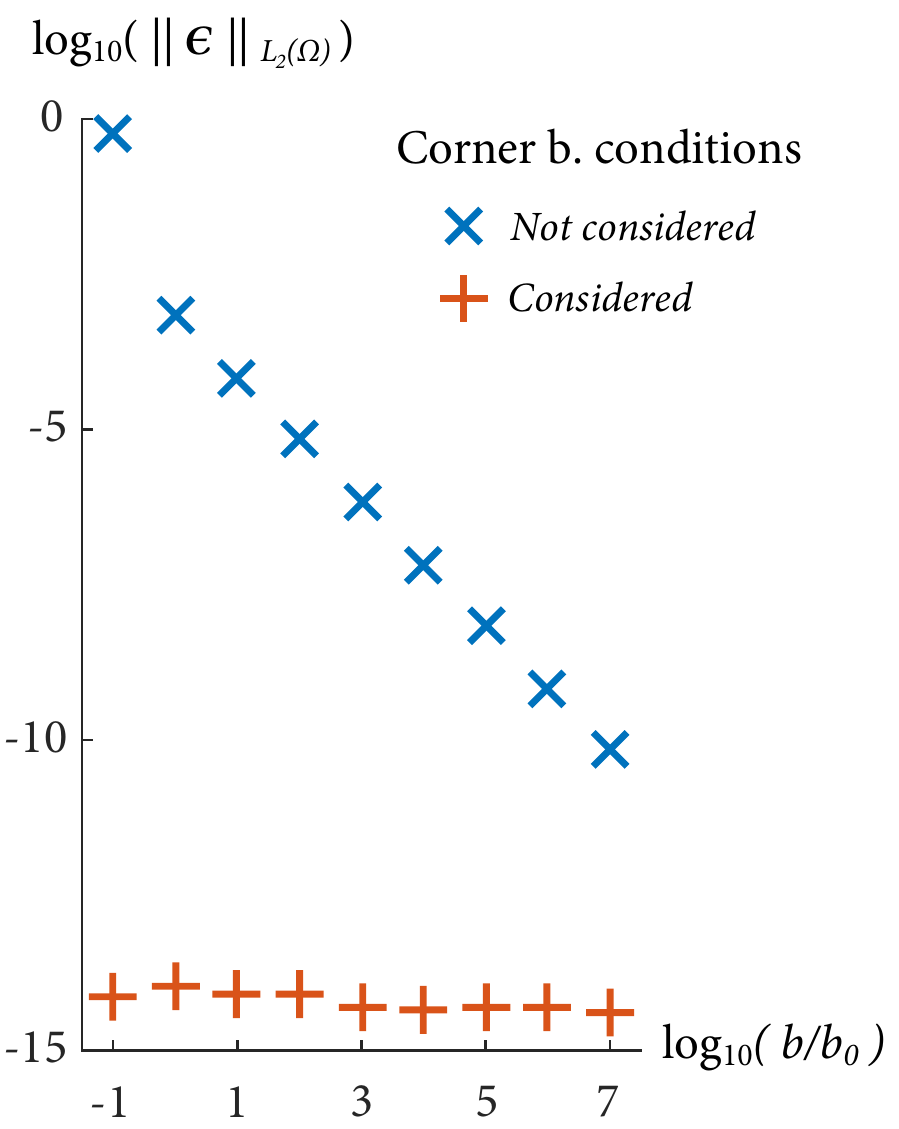}
		\label{fig_ex_corner_result}
	\end{subfigure}%
	\\
	\begin{subfigure}[t]{0.25\textwidth}\centering
		\vspace{-7.5em}\captionof{figure}{}
		\includegraphics[width=6em]{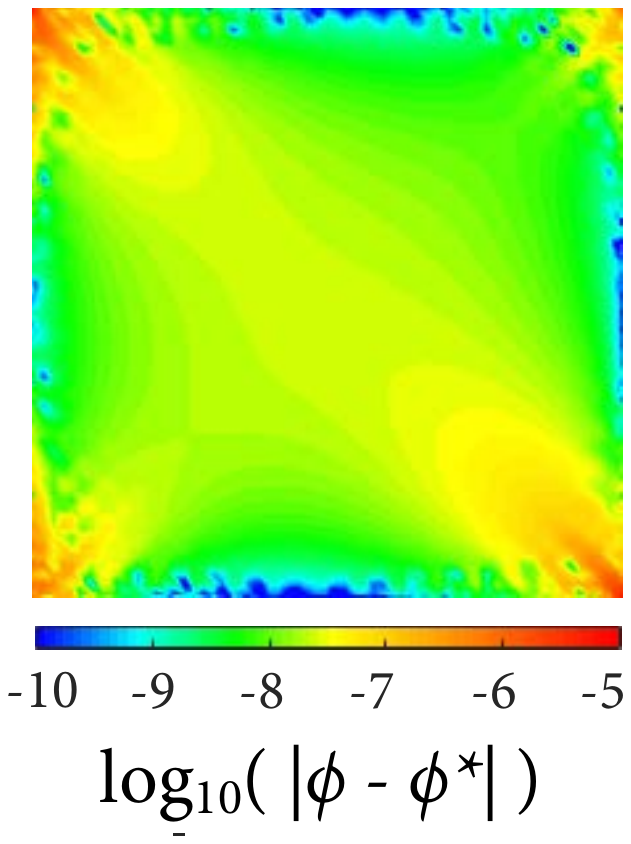}
		\label{fig_ex_corner_plot}
	\end{subfigure}
	\caption{Numerical example \ref{sec_edgebc} on corner boundary conditions. \emph{a)} Geometrical 2D model, \emph{b)} Unfitted uniform B-spline mesh, \emph{c)} $L_2$ norm of the error $\toVect{\upepsilon}$ as a function of the size of the domain for different non-local corner condition setups, and \emph{d)} point-wise error of the electric potential $\phi$ for $b/b_0=100$ without non-local corner conditions (results in $\log_{10}$-scale).}
\end{figure}
The effect of neglecting non-local corner conditions is illustrated in the following example. To this end, we simulate the flexoelectric effect on a $[-b,b]\times[-b,b]$ domain (\fig\ref{fig_ex_corner_setup}).
We consider the following synthetic solution $\{\Displacement^*,\phi^*\}$:
\begin{subequations}\label{bcond1}\begin{align}
{u_1}^*(\x)&\coloneqq \bar x_1+{\bar x_1}^2-2\bar x_1\bar x_2+{\bar x_1}^3-3\bar x_1{\bar x_2}^2+{\bar x_1}^2\bar x_2;\\
{u_2}^*(\x)&\coloneqq -\bar x_2+{\bar x_2}^2-2\bar x_1\bar x_2+{\bar x_2}^3-3{\bar x_1}^2\bar x_2-\bar x_1{\bar x_2}^2;\\
{\phi}^*(\x)&\coloneqq {\bar x_1}^3+{\bar x_2}^2-2{\bar x_1}^2\bar x_2;
\end{align}\end{subequations}
which can be exactly represented by a cubic ($p=3$) B-spline basis, where $\bar x_i\coloneqq x_i/b$. The prescribed terms on the boundary and source terms on the bulk required for \eq\eqref{bcond1} to be the solution of the flexoelectric problem are computed by inserting \eq\eqref{bcond1} into \eq\eqref{eq_EulerLagrange} and \eqref{eq_FlexoForces}.

Dirichlet boundary conditions are considered at the boundary of the square for mechanics in \eq\eqref{strgr_classicalBC},\eqref{strgr_nonlocalBC} and electrostatics in \eq\eqref{elecBC}, \ie
\begin{equation}\begin{rcases}
		\Displacement &=\Displacement^*                     \\
		\partial^n\Displacement &=\partial^n\Displacement^* \\
		\phi &=\phi^*                                     
\end{rcases}
\quad\text{at}\quad \{\bar x_1=-1;\quad \bar x_1=1;\quad \bar x_2=-1;\quad \bar x_2=1\}.
\end{equation}
According to the formulation in Section \ref{sec_02}, additional mechanical Dirichlet conditions in \eq\eqref{strgr_edgeBC} arise at the \emph{corners} of $\partial\Omega$, namely
\begin{equation}\label{eq_ex1_corners}
\Displacement =\Displacement^*                                    
\quad\text{at}\quad (\bar x_1,\bar x_2)=\{(-1,-1)\cup(1,-1)\cup(1,1)\cup(-1,1)\}.
\end{equation}
Two cases are considered, depending whether corner conditions in \eq\eqref{eq_ex1_corners} are enforced or not. In both cases the size $b$ of the domain is varied to assess the effect of corner conditions at different length scales.
For the numerical approximation of the solution we consider a cubic ($p=3$) uniform B-spline basis on an unfitted mesh of $32\times32$ cells (\fig\ref{fig_ex_corner_mesh}).
The material parameters are chosen as follows:
\begin{center}
$\nonumber
\quad
E=\SI{152}{GPa};
\quad
\nu=0.33;
\quad
l=\SI{1}{\nm};
\quad
\dielec_L=\SI{141}{\nano\joule\per\square\volt\per\meter};
\quad
\toVect{d_\text{piezo}}=x_2;
\quad
\piezo_L=\SI{8.8}{\joule\per\volt\per\square\meter};
$

$
\piezo_T=\SI{-4.4}{\joule\per\volt\per\square\meter};
\quad
\piezo_S=\SI{4.4}{\joule\per\volt\per\square\meter};
\quad
\flexo_L=\SI{150}{\micro\joule\per\volt\per\meter};
\quad
\flexo_T=\SI{110}{\micro\joule\per\volt\per\meter};
\quad
\flexo_S=\SI{110}{\micro\joule\per\volt\per\meter}.
$\end{center}

Figure \ref{fig_ex_corner_result} shows the $L_2$ error of the numerical solution at different length scales, computed as the $L_2$ norm of the error vector $\toVect{\upepsilon}\coloneqq(u_1-{u_1}^*,u_2-{u_2}^*,\phi-\phi^*)$ on $\Omega$. Results reveal that the synthetic solution is exactly captured (up to round-off errors) only in the case corner conditions are considered. Otherwise, a significant error is found, which grows as the domain size is reduced, \ie at scales where strain-gradient elasticity and flexoelectricity couplings are relevant. The source of this error is not numerical, since the geometry is exactly represented and integrated, and the approximation space captures the analytical solution exactly. 

The error introduced by neglecting non-local corner conditions is illustrated in \fig\ref{fig_ex_corner_plot} for the scale $b/b_0=100$, where $b_0\coloneqq l=\SI{1}{\nm}$ is a characteristic length of the problem used for normalization. The absolute value of the electric potential component of the error $\toVect{\upepsilon}$, \ie $|\phi-\phi^*|$, is depicted in $\log$-scale within $\Omega$. The error is concentrated around the corners of the domain, where the two boundary value problems defer. Similar behavior is observed for the mechanical components of the error $\toVect{\upepsilon}$.

We conclude that non-local corner conditions are a fundamental part of the mathematical prescription of the physical problem of linear flexoelectricity in the presence of non-smooth domains. Ignoring them leads to solving a different boundary value problem, and the resulting discrepancies can be very significant below length scales where strain-gradient effects play a relevant role. Therefore, in the following numerical examples, non-local corner conditions are always properly considered.

% % % % %
\subsection{Sensitivity analysis on the penalty parameter}\label{SensitivityAnalysis}
% % % % %
\begin{figure}[t!]\centering
	\begin{subfigure}[t]{0.45\textwidth}\centering
		\vspace{1em}\captionof{figure}{Mechanical submatrix $\toMat{K}_{\Displacement\Displacement}$}\hspace{1.05em}
		\includegraphics[height=10em]{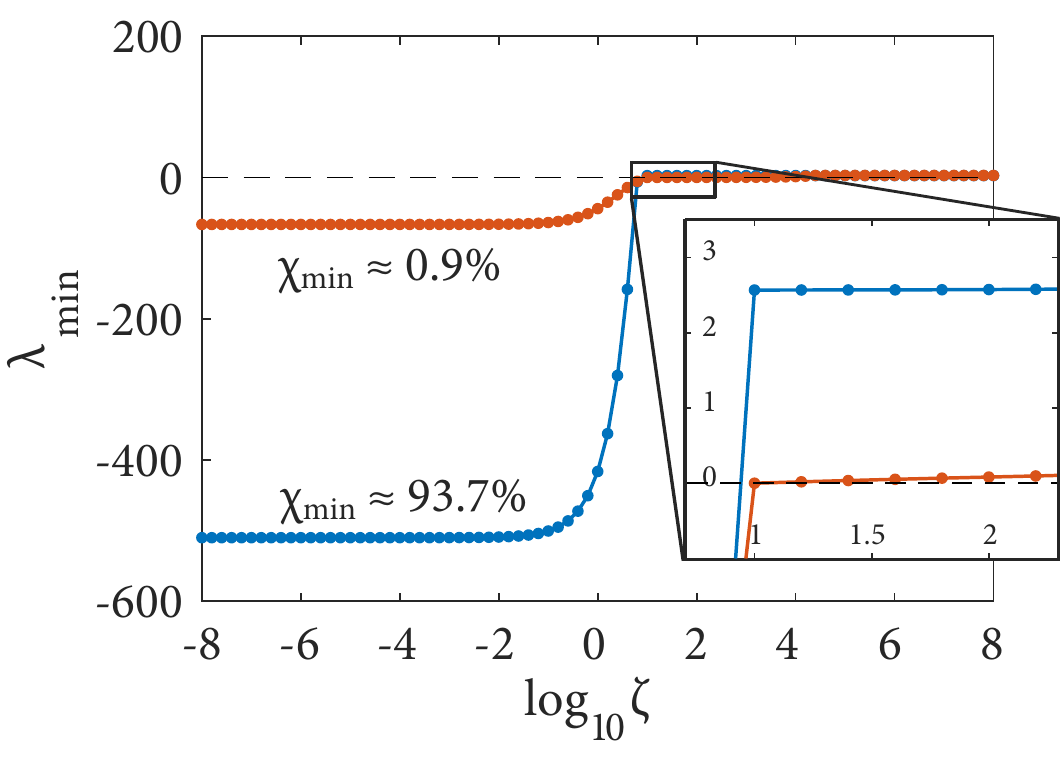}
		\label{EVmech}
	\end{subfigure}\hfill
	\begin{subfigure}[t]{0.45\textwidth}\centering
		\vspace{1em}\captionof{figure}{Electrical submatrix $\toMat{K}_{\phi\phi}$}
		\includegraphics[height=10em]{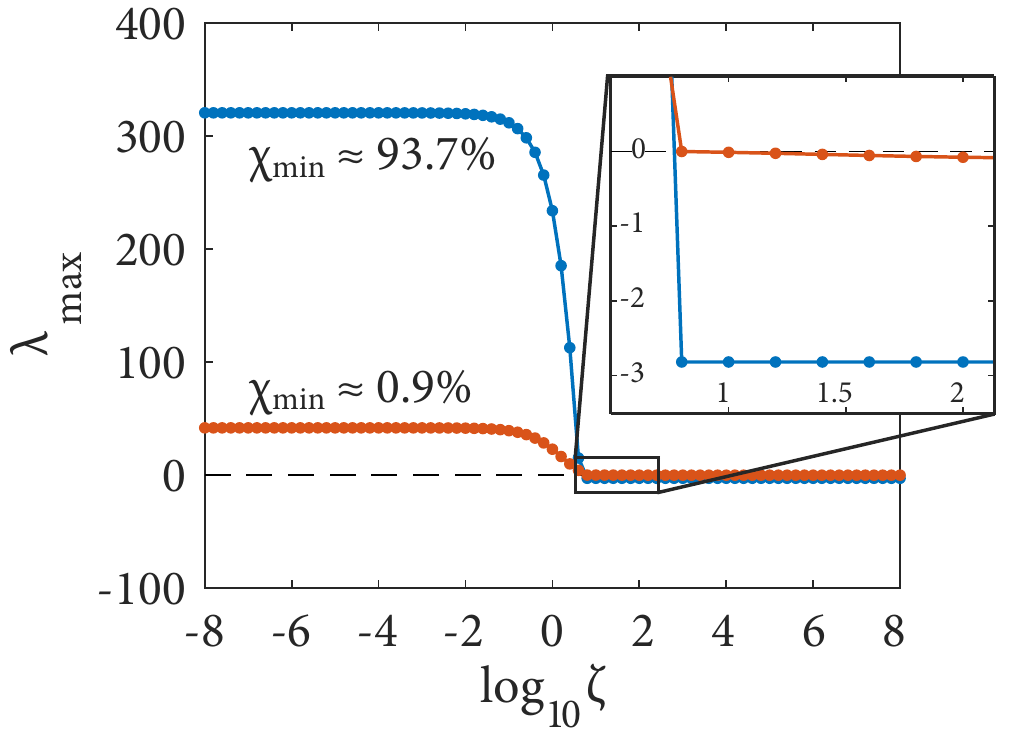}
		\label{EVelec}
	\end{subfigure}
	\caption{Sensitivity of the eigenvalues $\lambda$ with respect to the penalty parameter $\zeta$ for meshes with minimum volume fractions $\chi_\text{min}\approx93.7\%$ (blue) and $\chi_\text{min}\approx0.9\%$ (red). \emph{a)} Minimum eigenvalue of the mechanical submatrix $\toMat{K}_{\Displacement\Displacement}$. \emph{b)} Maximum eigenvalue of the electrical submatrix $\toMat{K}_{\phi\phi}$. The dashed lines indicate $\lambda=0$ for reference. The insets show a zoom around $\log_{10}\zeta=1.5$.}
	\label{EV}	
\end{figure}

\begin{figure}[t!]\centering
	\begin{subfigure}[t]{0.45\textwidth}\centering
		\vspace{1em}\captionof{figure}{Numerical error}
		\includegraphics[height=20em]{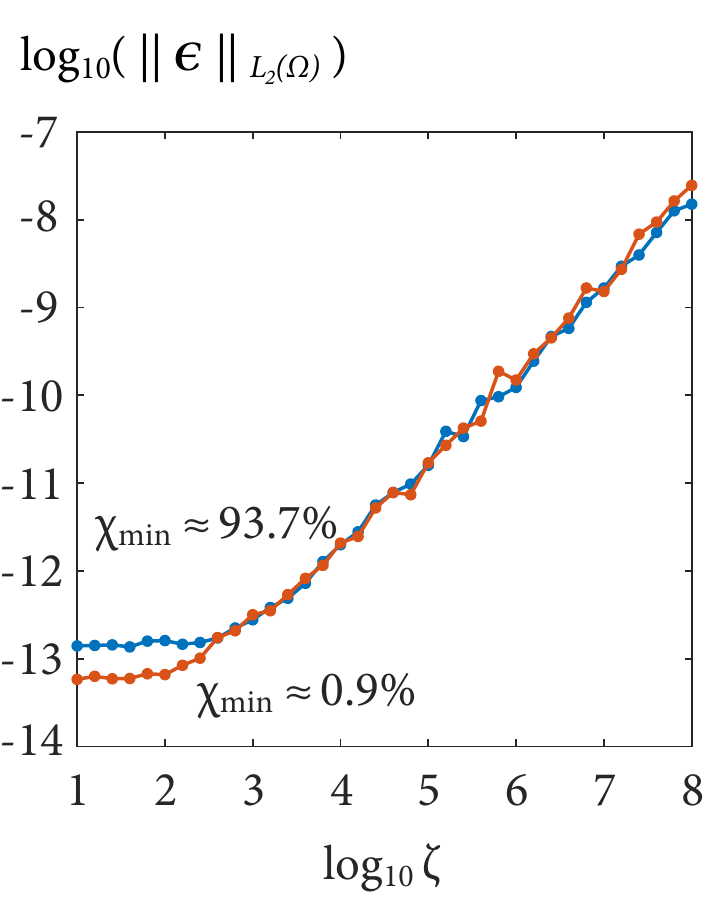}
		\label{ERR}
	\end{subfigure}\qquad
	\begin{subfigure}[t]{0.45\textwidth}\centering
		\vspace{1em}\captionof{figure}{Condition number}
		\includegraphics[height=20em]{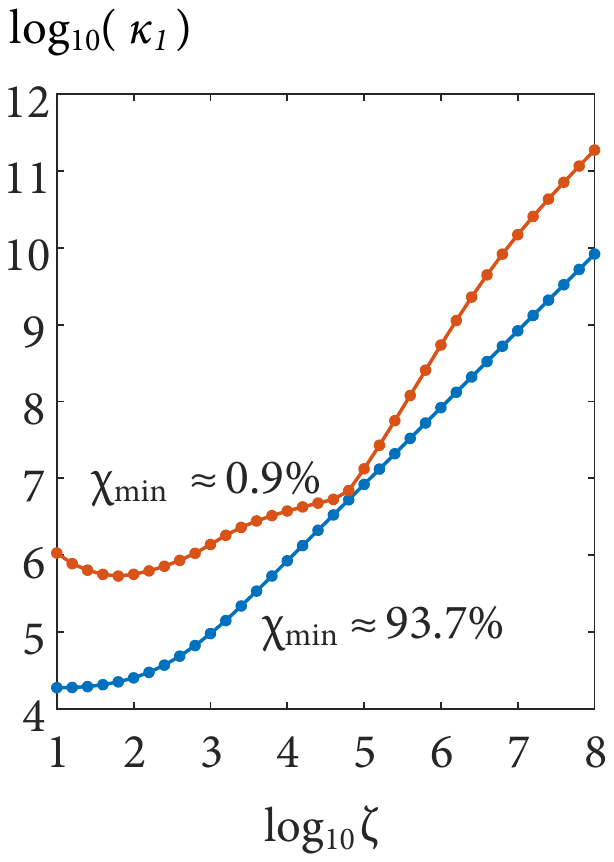}
		\label{CN}
	\end{subfigure}
	\caption{Sensitivity of \emph{a)} $L_2$ norm of the numerical error $\toVect{\upepsilon}$ and \emph{b)} condition number $\kappa_1$ of the system matrix, with respect to the penalty parameter $\zeta$, for meshes with minimum volume fractions $\chi_\text{min}\approx93.7\%$ (blue) and $\chi_\text{min}\approx0.9\%$ (red).}
	\label{ERRCN}	
\end{figure}

We perform a sensitivity analysis on the penalty parameter $\zeta$ in \eq\eqref{parameterzeta}. We consider the setup in the previous example, with two different meshes leading to minimum volume fractions $\chi_\text{min}\approx93.7\%$ and $\chi_\text{min}\approx0.9\%$ respectively, where $\chi_\text{min}\coloneqq\min_c\left(\frac{|\Omega_\square^c\cap\Omega|}{|\Omega_\square^c|}\right)$ denotes the volume fraction of the smallest cut cell in the mesh.

According to \eq\eqref{eq_2var}, the second variations of the energy functional with respect to the mechanical and electrical unknowns are required to be positive and negative, respectively, consistent with the $\min$-$\max$ nature of the problem and to ensure stability of the formulation. Numerically, this is equivalent to checking that the minimum eigenvalue of $\toMat{K}_{\Displacement\Displacement}$ and the maximum eigenvalue of $\toMat{K}_{\phi\phi}$ are positive and negative respectively, where $\toMat{K}_{\Displacement\Displacement}$ and $\toMat{K}_{\phi\phi}$ are the submatrices of the system related to mechanical and electrical equations and unknowns, respectively.

\fig\ref{EV} shows the minimum eigenvalue of $\toMat{K}_{\Displacement\Displacement}$ and the maximum eigenvalue of $\toMat{K}_{\phi\phi}$ as a function of $\zeta$. The $\min$-$\max$ condition stated above is fulfilled for $\zeta>10$. It is important noting that this threshold does not depend on the minimum volume fraction $\chi_\text{min}$ of the mesh. This result is in agreement with \cite{badia2018}, where the threshold of the penalty parameter $\zeta$ is analytically proven to be mesh-independent after using cell-aggregation stabilization schemes on unfitted meshes, as the extended B-splines method we consider in this work (see remark in Section \ref{stab}).

The sensitivity of both the numerical error, $\toVect{\upepsilon}$, and of the condition number of the system matrix, ${\kappa_1}$, on the penalty parameter $\zeta$ is shown in \fig\ref{ERRCN}. For both meshes, machine-precision accuracy is achieved for moderate values of $\zeta$, \ie $10<\zeta<1000$. However, for larger values of $\zeta$ the errors increase due to the increase in the condition number, which grows proportionally to $\zeta$. In view of the reported results, we consider $\zeta=100$ for the following numerical examples. 
% % % % %
\subsection{Accuracy and convergence properties of the method}
% % % % %

\begin{figure}[t!]\centering
	\begin{subfigure}[t]{0.23\textwidth}\centering
		\captionof{figure}{}
		\includegraphics[height=18em]{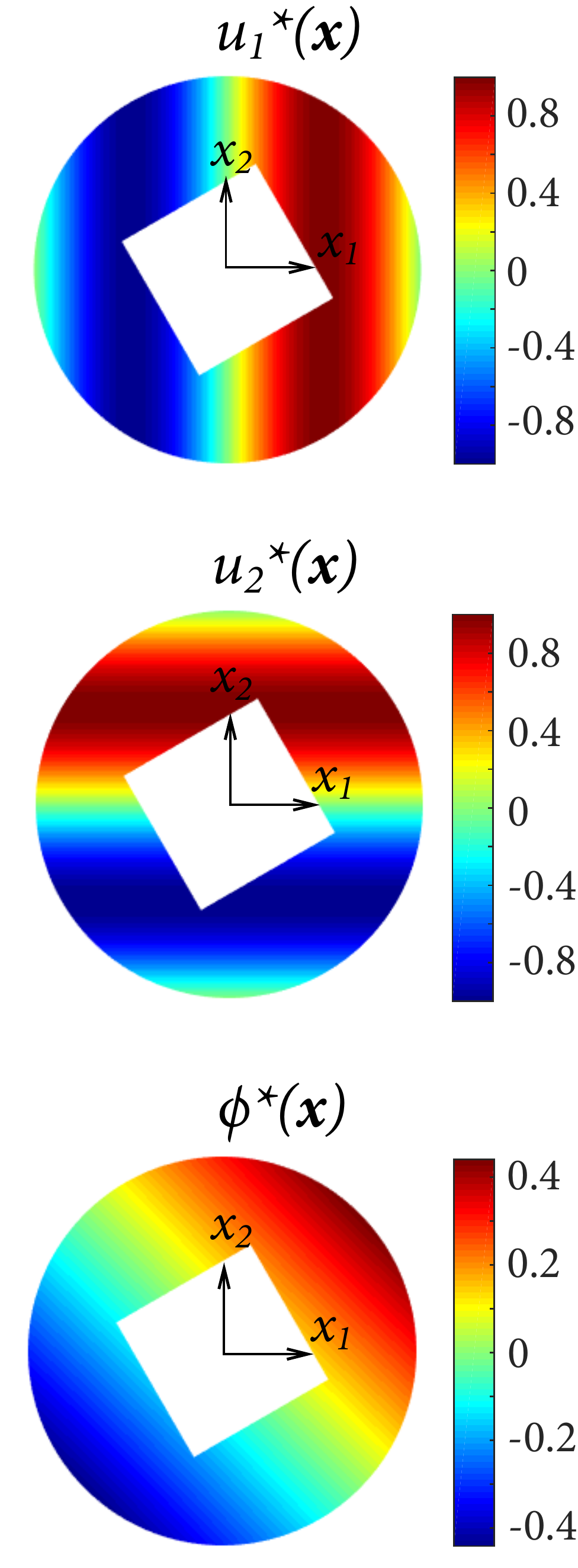}
		\label{fig_ex_conv_analytical}
	\end{subfigure}%
	\begin{subfigure}[t]{0.77\textwidth}\centering
		\captionof{figure}{}
		\includegraphics[height=18em]{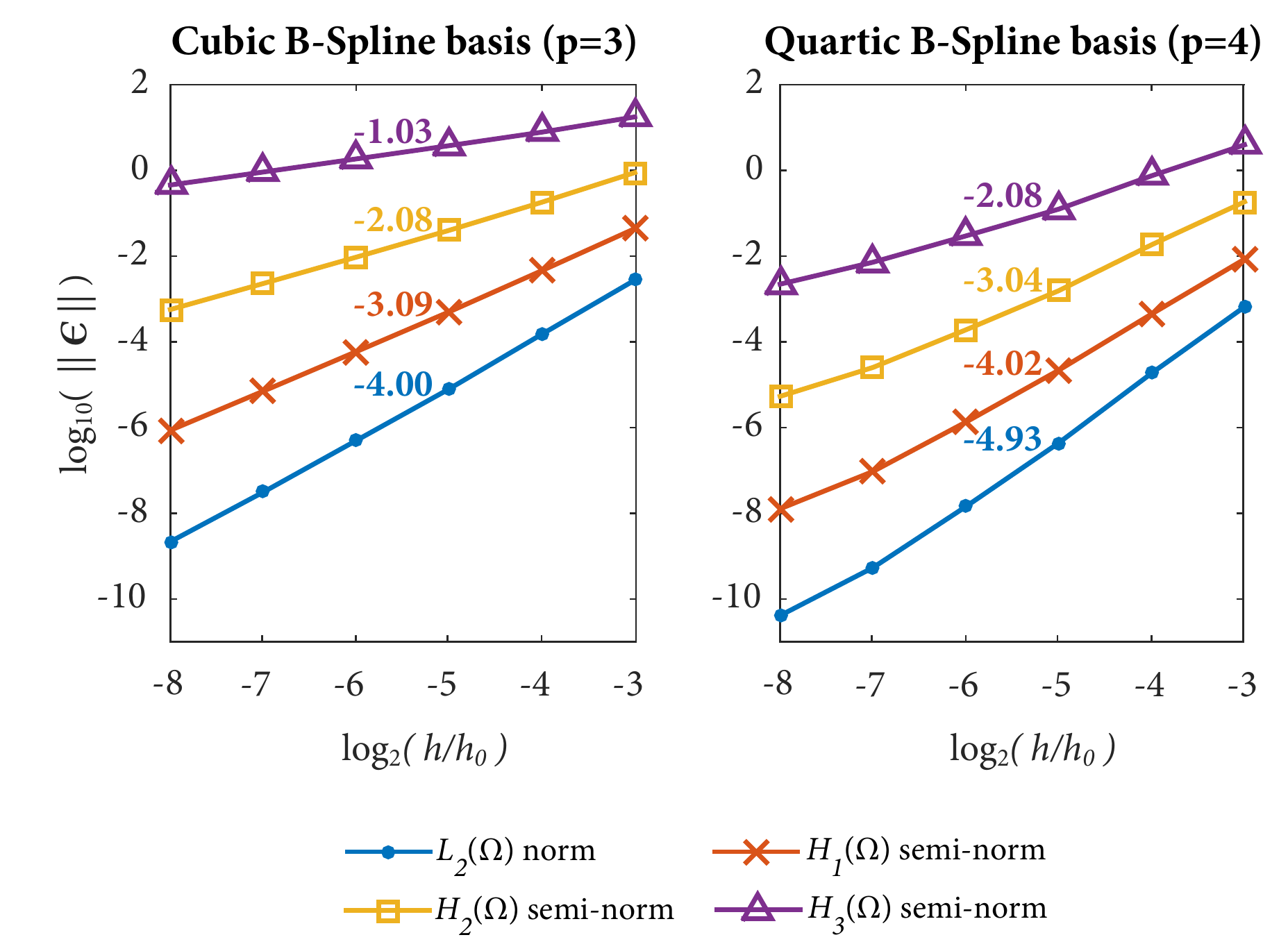}
		\label{fig_ex_optimalconvergence}
	\end{subfigure}
	\caption{Error convergence analysis. \emph{a)} Synthetic solution $\{\Displacement^*,\phi^*\}$ on the physical domain $\Omega$, and \emph{b)} convergence plots for B-spline basis of polynomial degree $p=\{3,4\}$. The numerical error $\toVect{\upepsilon}$ is measured in the $L_2$ norm and $H_1$, $H_2$ and $H_3$ semi-norms, and plotted as a function of the cell size $h$. The numbers on the curves denote the error convergence rates between the meshes $h/h_0=2^{-5}$ and $h/h_0=2^{-6}$, computed as the slope of the curves multiplied by a factor of $\log_2(10)$.}
	\label{fig_ex_convergence}
\end{figure}

In this example the convergence of the method is assessed for different high-order B-spline approximations, namely at $p=\{3,4\}$.
We consider the following synthetic solution $\{\Displacement^*,\phi^*\}$:
\begin{subequations}\begin{align}
		{u_1}^*(\x)&\coloneqq \sin\left(x_1\pi\right);\\
		{u_2}^*(\x)&\coloneqq \sin\left(x_2\pi\right);\\
		{\phi}^*(\x)&\coloneqq \sin\left(\frac{x_1}{10}\pi\right)+\sin\left(\frac{x_2}{10}\pi\right);
\end{align}\end{subequations}
depicted in \fig\ref{fig_ex_conv_analytical} on the domain $\Omega$, which is conformed by a circle of radius $\SI{1}{\micro\meter}$ with a square-shaped hole of size
$\SI{0.4 x 0.4}{\um}$
rotated 30 degrees with respect to the Cartesian coordinates; both geometries are centered at $\x=(0,0)$.

Dirichlet conditions are enforced at the boundaries and the corners of the domain. The boundary is exactly mapped by means of the NEFEM mapping, which makes the geometrical error vanish. The numerical integration is rich enough so that the integration error is negligible and does not pollute the convergence plots.

The material parameters are chosen as follows:
\begin{center}
$
\quad
E=\SI{100}{GPa};
\quad
\nu=0.37;
\quad
l=\SI{2}{\nm};
\quad
\dielec_L=\SI{11}{\nano\joule\per\square\volt\per\meter};
\quad
\toVect{d_\text{piezo}}=x_2;
\quad
\piezo_L=\SI{8.8}{\joule\per\volt\per\square\meter};
$

$
\piezo_T=\SI{-4.4}{\joule\per\volt\per\square\meter};
\quad
\piezo_S=\SI{1.1}{\joule\per\volt\per\square\meter};
\quad
\flexo_L=\SI{0.5}{\micro\joule\per\volt\per\meter};
\quad
\flexo_T=\SI{1}{\micro\joule\per\volt\per\meter};
\quad
\flexo_S=\SI{0.5}{\micro\joule\per\volt\per\meter}.
$\end{center}

Error convergence results for cubic $(p=3)$ and quartic $(p=4)$ B-spline bases are shown in \fig\ref{fig_ex_optimalconvergence}. The error $\toVect{\upepsilon}$ is measured in the $L_2$ norm and $H_1$, $H_2$ and $H_3$ semi-norms at six recursively-refined uniform meshes of cell size 
$h/h_0=\{\num{2e-3},\dots,\num{2e-8}\}$, where $h_0\coloneqq\SI{1}{\um}$ is a normalization factor.
As expected, optimal convergence rates are (asymptotically) achieved with both B-spline bases. The asymptotic behavior of the error convergence rate is expected, and is due to the extended B-spline stabilization performed on trimmed basis functions near the boundary $\partial\Omega$ (see remark in Section \ref{stab}). For relatively coarse meshes, a small additional error is introduced since the approximation space is coarsened at the boundary. However, for fine enough meshes, this effect is negligible and error convergence rates tend to optimality \cite{Hollig2001,Hollig2012}.

% % % % %
\subsection{Transversal transduction: cantilever beam under bending}
% % % % %
Bending a cantilever beam is a natural way to mobilize transversal strain-gradients and, consequently, to trigger the transverse flexoelectric effect. For this reason, electroactive beam bending has been extensively used in experiments for the characterization of the transversal flexoelectric coefficient \cite{Baskaran2012,Chu2012,Ma2002,Ma2005,Ma2006}. This setup has also been modeled with approximate analytical \cite{Majdoub2009} and numerical \cite{Abdollahi2014} models.

\begin{figure}[p]\centering
	\begin{subfigure}{\textwidth}\centering
		\vspace{1em}\captionof{figure}{}\vspace{-0.5em}
		\includegraphics[scale=1]{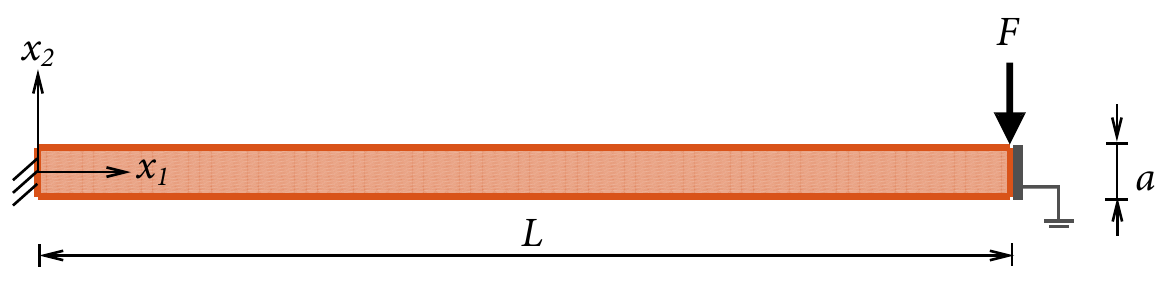}
		\label{fig_ex_beam}\\
	\end{subfigure}
	\begin{subfigure}{.450\textwidth}\centering
		\vspace{1em}\captionof{figure}{}\vspace{-0.5em}
		\includegraphics[width=\textwidth]{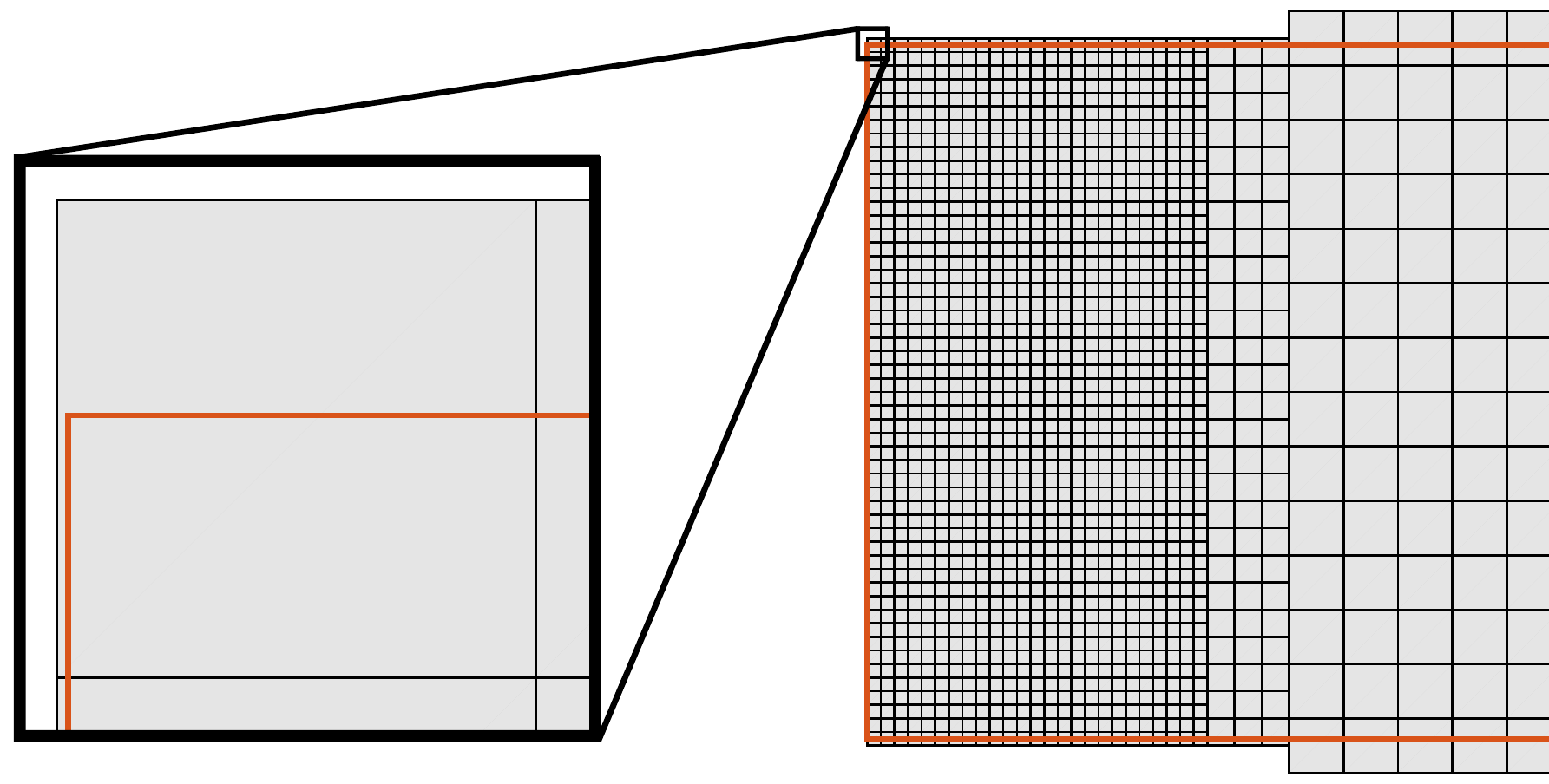}
		\label{fig_ex_beam_mesh1}
	\end{subfigure}\hfill
	\begin{subfigure}{.450\textwidth}\centering
		\vspace{1em}\captionof{figure}{}\vspace{-0.5em}
		\includegraphics[width=\textwidth]{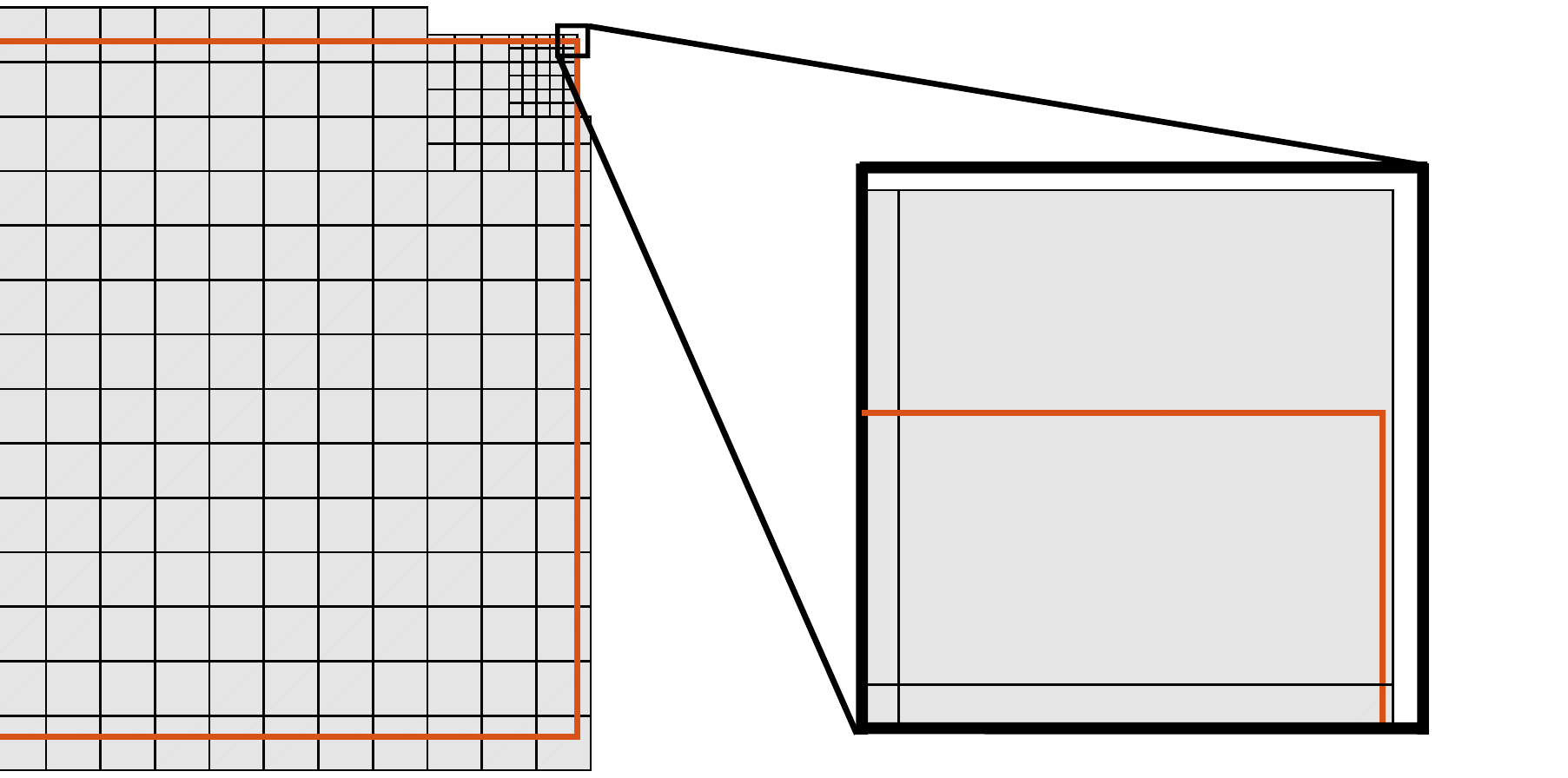}
		\label{fig_ex_beam_mesh2}
	\end{subfigure}
	\caption{Electroactive cantilever beam. \emph{a)} Geometrical 2D model, \emph{b)} detail of the hierarchical mesh at the left-end of the beam and \emph{c)} detail of the hierarchical mesh at the right-end of the beam. The mesh does not conform to the boundary of the beam at any level of refinement.}
	\label{fig_ex_BEAM}
\end{figure}

Figure \ref{fig_ex_beam} depicts the geometrical model of a cantilever beam.
The aspect ratio of the beam is fixed to $L/a=20$, where $L$ represents the length of the beam and $a$ its width. We consider different lengths $L$ in order to capture the size-dependent nature of the flexoelectric coupling.
Mechanically, the beam is clamped on its left-end and undergoes a point force $F$ on its top-right-corner.
Electrically, it is grounded on its right-end, and the other edges are considered charge-free.
The corresponding boundary conditions are:
\begin{subequations}\begin{align}
u_1=u_2=0 &\quad\text{at}\quad x_1=0, \\
(\partial^n\Displacement)_2=-u_{2,1}=0 &\quad\text{at}\quad x_1=0, \\
u_1=u_2=0 &\quad\text{at}\quad (0,-h/2)\cup(0,h/2), \\
j_2=-F &\quad\text{at}\quad (L,h/2), \\
\phi=0 &\quad\text{at}\quad x_1=L.
\end{align}\end{subequations}

Strain-gradient elasticity is neglected to isolate the effect of piezoelectricity and flexoelectricity couplings. Three different cases are considered for the electromechanics: \emph{i)} piezoelectricity only, \emph{ii)} flexoelectricity only and \emph{iii)} combined piezoelectricity and flexoelectricity.

The beam accommodates the mechanical load by bending, which produces a linear distribution of the axial strain $\strain_{11}$ along the transversal ($x_2-$) direction which is well known in classical elasticity. Namely, both \emph{i)} strain and \emph{ii)} strain gradients are generated, which are the triggers for direct piezoelectric and flexoelectric effects, respectively. Therefore, a non-zero electric field is generated on the sample as a consequence of the mechanical loading.
For a piezoelectric beam, we expect the electromechanical response to be the same regardless of the size of the sample. However, for a flexoelectric or flexo-piezoelectric beam it should grow inversely to the scale due to the size-dependent nature of the strain gradient field and thus of flexoelectric coupling \cite{Abdollahi2014}.

\begin{figure}[p]\centering
	\begin{subfigure}{.450\textwidth}\centering
		\vspace{1em}\captionof{figure}{}\vspace{-0.5em}
		\includegraphics[width=\textwidth]{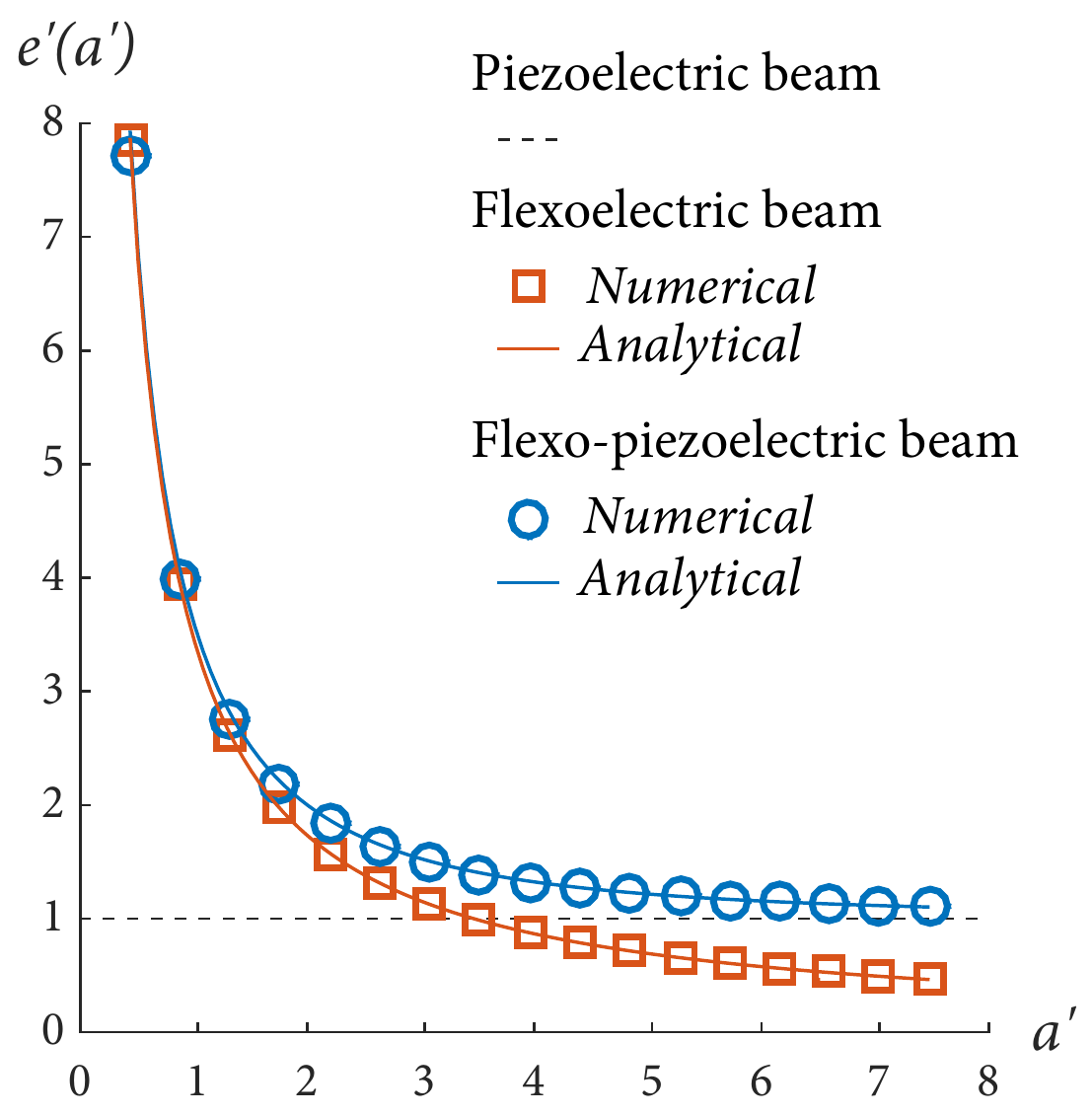}
		\label{fig_PAPER01_Amir2D_ECF}
	\end{subfigure}\hfill
	\begin{subfigure}{.450\textwidth}\centering
		\vspace{-0.6em}\captionof{figure}{Piezoelectric beam}\vspace{-0.5em}
		\includegraphics[width=\textwidth]{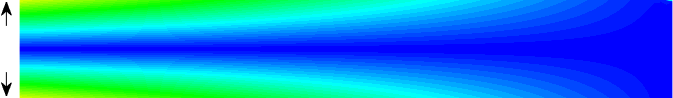}
		\label{fig_PAPER01_Amir2D_E_Piezo}
		\vspace{1em}\captionof{figure}{Flexoelectric beam}\vspace{-0.5em}
		\includegraphics[width=\textwidth]{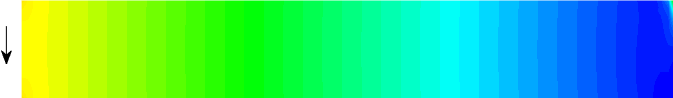}
		\label{fig_PAPER01_Amir2D_E_Flexo}
		\vspace{1em}\captionof{figure}{Flexo-piezoelectric beam}\vspace{-0.5em}
		\includegraphics[width=\textwidth]{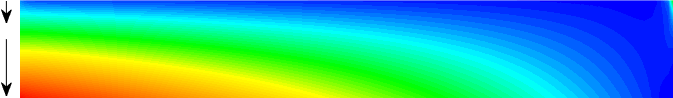}
		\\[2em]
		\includegraphics[width=0.6\textwidth]{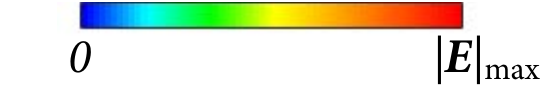}
		\label{fig_PAPER01_Amir2D_E_Flexopiezo}
	\end{subfigure}
	\caption{Numerical results of the electroactive cantilever beam. \emph{a)} Normalized effective piezoelectric constant $e'$ against normalized beam thickness $a'$ for the cases of flexoelectric and flexo-piezoelectric beams. Numerical results match the analytical estimations in \eq\eqref{eq_e}. \emph{b)-d)} Qualitative spatial distribution of the electric field modulus $|\E|$ at scale $a'=1.76$ for the cases of piezoelectric, flexoelectric and flexo-piezoelectric beams. The electric field direction is represented by arrows.}
	\label{fig_PAPER01_Amir2D}
\end{figure}

In order to quantify the energy conversion, we define the \emph{electromechanical coupling factor} $k_\text{eff}$ as
\begin{equation}
k_\text{eff}:=\left(\frac
{\frac{1}{2}\int_{\Omega}\E\cdot\Dielec\cdot\E~\mathrm{d}\Omega}
{\frac{1}{2}\int_{\Omega}\Strain\colon\Elast\colon\Strain~\mathrm{d}\Omega}\right)^\frac{1}{2},
\end{equation}
which is a positive, dimensionless scalar that indicates the relationship between dielectric and mechanic energies required to accommodate the external mechanical load. 

An analytical estimation of $k_\text{eff}$ for the flexo-piezoelectric beam can be found in the literature \cite{Majdoub2009} as
\begin{equation}\label{eq_keff}
k_\text{eff}\approx
\frac{\kappa_L(\kappa_L-\varepsilon_0)}{E^2}\sqrt{{\piezo_T}^2+12\frac{{\flexo_T}^2}{h^2}},
\end{equation}
where $\varepsilon_0\approx\SI{8.854e-12}{\coulomb\per\volt\per\meter}$ is the vacuum permittivity constant, and it has been assumed that the material parameters $\nu=l=\mu_L=\mu_S=\piezo_L=\piezo_S=0$ (see \ref{sec_app03}).

One can also define the \emph{normalized effective piezoelectric constant} $e'$ \cite{Majdoub2009} as 
\begin{equation}\label{eq_eprime}
e'\coloneqq\frac{k_\text{eff}}{k_\text{eff}\big|_{\mu_T=0}},
\end{equation}
which indicates the ratio between the current $k_\text{eff}$ and the $k_\text{eff}$ that would be obtained if the beam was purely piezoelectric. By combining \eq\eqref{eq_keff} and \eqref{eq_eprime}, analytical estimations of $e'$ are obtained for flexoelectric and flexo-piezoelectric beams as
\begin{align}\label{eq_e}
e'\big|_\text{flexo}(a')\approx\sqrt{\frac{12}{{a'}^2}},&&e'\big|_\text{flexo-piezo}(a')\approx\sqrt{1+\frac{12}{{a'}^2}};
\end{align}
where $a'\coloneqq-a\piezo_T\flexo_T^{-1}$ is the \emph{normalized beam thickness}.

For the numerical approximation of the solution let us consider a cubic ($p=3$) hierarchical B-spline basis on a mesh with two levels of refinement around the left-end and top-right-corner of the beam, as depicted in \fig\ref{fig_ex_beam_mesh1} and \ref{fig_ex_beam_mesh2}. The material parameters are the following:
\begin{center}$
\nu=l=\flexo_L=\flexo_S=\piezo_L=\piezo_S=0;
\quad
E=\SI{100}{GPa};
\quad
\dielec_L=\SI{11}{\nano\joule\per\square\volt\per\meter};
\quad
\piezo_T=\SI{-4.4}{\joule\per\volt\per\square\meter};
\quad
\flexo_T=\SI{1}{\micro\joule\per\volt\per\meter}.
$\end{center}

The normalized effective piezoelectric constant $e'(a')$ is depicted in \fig\ref{fig_PAPER01_Amir2D_ECF}, which shows very good agreement between our numerical results and the analytical estimations in \eq\eqref{eq_e}. The flexoelectric effect is negligible at large scales; as a consequence, for large $a'$ the electromechanical response tends to vanish for the flexoelectric beam, whereas the flexo-piezoelectric one tends to behave as purely piezoelectric. On the contrary, at smaller scales the flexoelectric effect is much more relevant and leads to an enhanced electromechanical transduction in both cases. 

Figures \ref{fig_PAPER01_Amir2D_E_Piezo}-\ref{fig_PAPER01_Amir2D_E_Flexopiezo} depict the spatial distribution of the mechanically-induced electric field $\E$ at scale $a'=1.76$ for the cases of piezoelectric, flexoelectric and flexo-piezoelectric beams.
In the piezoelectric case, the distribution is skew-symmetric (divergent) with respect to the neutral axis of the beam, in accordance with the axial strain $\strain_{11}$ distribution. Highest values appear close the left-end and away from the neutral axis.
The other two cases have a similar $k_\text{eff}$ but very different distributions of the electric field $\E$. The electric field on the flexoelectric beam points downwards, and remains almost constant along the transversal direction while increases close to the clamped end. The flexo-piezoelectric beam, however, presents an inhomogeneous distribution which can be thought as a combination of the two previous ones. Depending on the scale, the piezoelectric effect dominates the flexoelectric one or viceversa. Here, $a'=1.76$ is intentionally chosen since it leads to comparable electromechanical effects.

% % % % %
\subsection{Longitudinal transduction: truncated pyramid under compression}
% % % % %
Another frequent experimental setup is the truncated pyramid compression, widely used by experimentalists to characterize the longitudinal flexoelectric coefficient \cite{Cross2006,Lu2016,Baskaran2011a,Baskaran2011b}. Although analytical expressions are not available for this more complex setup, numerical solutions have been developed by our group \cite{Abdollahi2014}.

\begin{figure}[t!]
	\begin{subfigure}[t]{.38\textwidth}\centering
		\captionof{figure}{}\vspace{4em}
		\includegraphics[width=\textwidth]{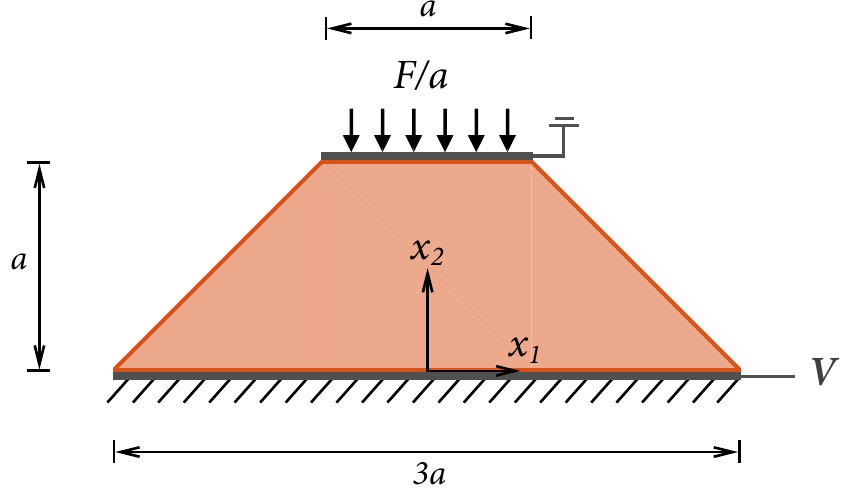}
		\label{fig_ex_trunc}
	\end{subfigure}%
	\begin{subfigure}[t]{.32\textwidth}\centering
		\captionof{figure}{}
		\includegraphics[width=.95\textwidth]{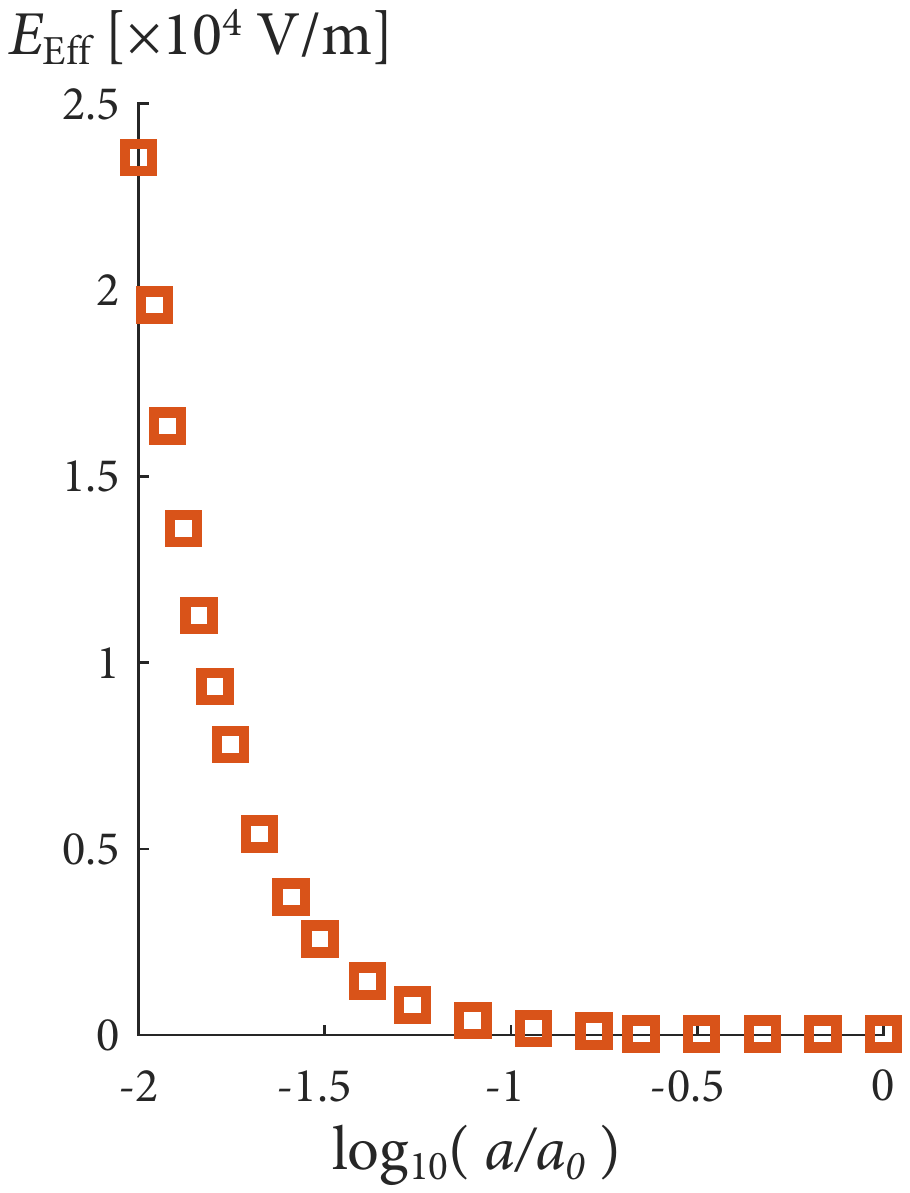}
		\label{fig_ex_trunc_results}
	\end{subfigure}%
	\begin{subfigure}[t]{.3\textwidth}\centering
		\captionof{figure}{}
		\includegraphics[width=0.9\textwidth]{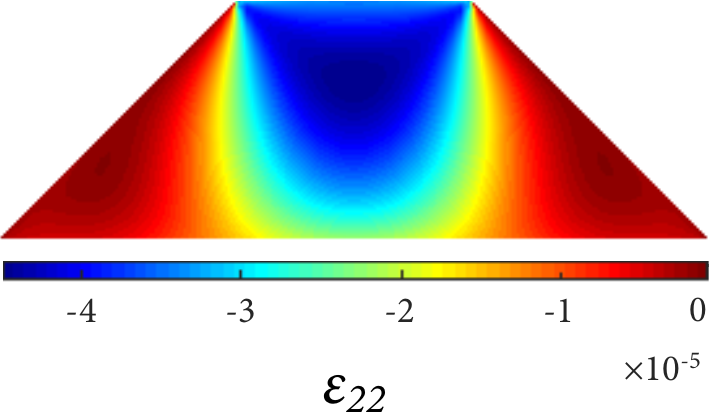}
		\label{fig_ex_trunc_strain}
		\\\vspace{1em}
		\captionof{figure}{}
		\includegraphics[width=0.9\textwidth]{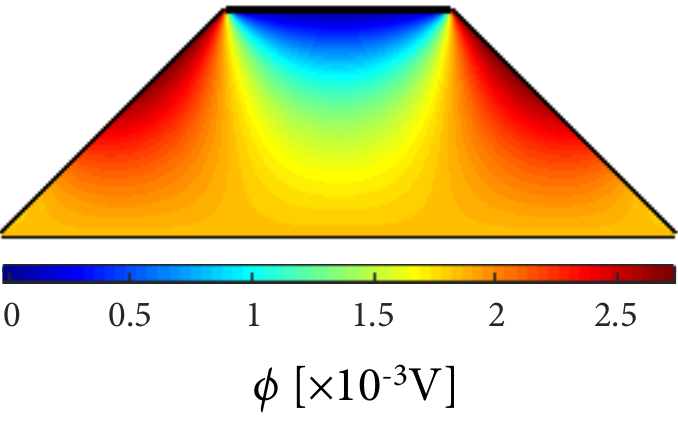}
		\label{fig_ex_trunc_potential}
	\end{subfigure}
	\caption{Electroactive truncated pyramid. \emph{a)} Geometrical 2D model. \emph{b)} Effective electric field $E_\text{Eff}$ as a function of the length scale, where $a_0=\SI{750}{\um}$ is a normalizing factor. \emph{c)} Strain $\strain_{22}$ distribution and \emph{d)} electric potential distribution.}
	\label{fig_trunc}
\end{figure}

Figure \ref{fig_ex_trunc} depicts the geometrical model of a flexoelectric truncated pyramid of height $a$ and bases $a$ (top) and $3a$ (bottom). The angle between bases and lateral boundaries is $\pi/4$.
Mechanically, the bottom basis is fixed and a compressive force $F$ is uniformly distributed on the top one.
Electrically, it is grounded on the top basis, and a sensing electrode is placed at the bottom.
The corresponding boundary conditions are:
\begin{subequations}\begin{align}
	u_1=u_2=0 &\quad\text{at}\quad x_2=0, \\
	u_1=u_2=0 &\quad\text{at}\quad (-3a/2,0)\cup(3a/2,0), \\
	t_2=-F/a &\quad\text{at}\quad x_2=a, \\
	\phi=0 &\quad\text{at}\quad x_2=a,\\
	\phi=V &\quad\text{at}\quad x_2=0;
\end{align}\end{subequations}
where $V$ is \emph{a priori} unknown but constant.

Since the top and bottom bases have different sizes, they undergo different compressive tractions and a longitudinal strain gradient in the $x_2$-direction arises as a result, which triggers the longitudinal flexoelectric effect. Another source of strain gradient is the bottom layer being fixed, which results in an inhomogeneous distribution of the traction along the bottom boundary.
Therefore, a non-zero electric field is generated on the sample as a consequence of the mechanical loading. A direct measure of the electromechanical transduction is the value $V$ of the sensing electrode at the bottom of the truncated pyramid. More interestingly, one can compute the \emph{effective electric field} $E_\text{Eff}$ measured as the voltage difference between electrodes over the height of the pyramid, \ie
\begin{equation}
E_\text{Eff}\coloneqq\frac{V}{a}.
\end{equation}

Numerical results are obtained with a cubic ($p=3$) B-spline basis on a uniform unfitted mesh of cell size $h\approx0.01195a$, and shown in \fig\ref{fig_trunc} for $F=\SI{4.5}{\newton\per\milli\metre}$, with the following material parameters:
\begin{center}$
l=\flexo_S=\piezo_L=\piezo_T=\piezo_S=0;
\quad
E=\SI{100}{\GPa};
\quad
\nu=0.37;
\quad
\dielec_L=\SI{11}{\nano\joule\per\square\volt\per\meter};
\quad
\flexo_L=\flexo_T=\SI{1}{\micro\joule\per\volt\per\meter}.
$\end{center}
Very good agreement with previous works in the literature \cite{Abdollahi2014} is reported.
The size-dependent nature of the flexoelectric coupling is evidenced in \fig\ref{fig_ex_trunc_results}, which shows the effective electric field $E_\text{Eff}$ as a function of the size of the truncated pyramid. Electromechanical transduction of the device takes place mainly at the micro- and nanoscale, whereas it is not relevant at larger scales.

In order to illustrate the complexity of the physics, the vertical strain $\strain_{22}$ and electric potential $\phi$ distributions at scale $a=7.5\um$ are depicted in \fig\ref{fig_ex_trunc_strain} and \ref{fig_ex_trunc_potential}, respectively. A highly inhomogeneous distribution of the strain takes place, specially near the corners of the grounded electrode on top of the device, which causes large strain-gradients triggering the flexoelectric effect. As a consequence, the electric potential distribution is also inhomogeneous within the domain. For this reason, simplified 1D models are not reliable to simulate the flexoelectric truncated pyramid, and numerical simulations are required \cite{Abdollahi2014}.

% % % % %
\subsection{Shear transduction: conical semicircular rod under torsion}
% % % % %
Unlike the longitudinal and transversal flexoelectric coefficients, the shear coefficient has been scarcely characterized experimentally. One reason is that, in many setups such as in the cylindrical rod torsion, shear strain gradients are effectively mobilized but the overall net polarization vanishes, and therefore no flexoelectric measurement can be effectively done. 
An alternative setup proposed by Mocci \etal\cite{MocciShearUnpublished} consists on a conical rod with semicircular cross section under torsion, where a net angular polarization arises thanks to the longitudinal variation of the cross section.

\begin{figure}[p]\centering
	\begin{subfigure}[t]{.4\textwidth}\centering
		\captionof{figure}{}
		\includegraphics[width=\textwidth]{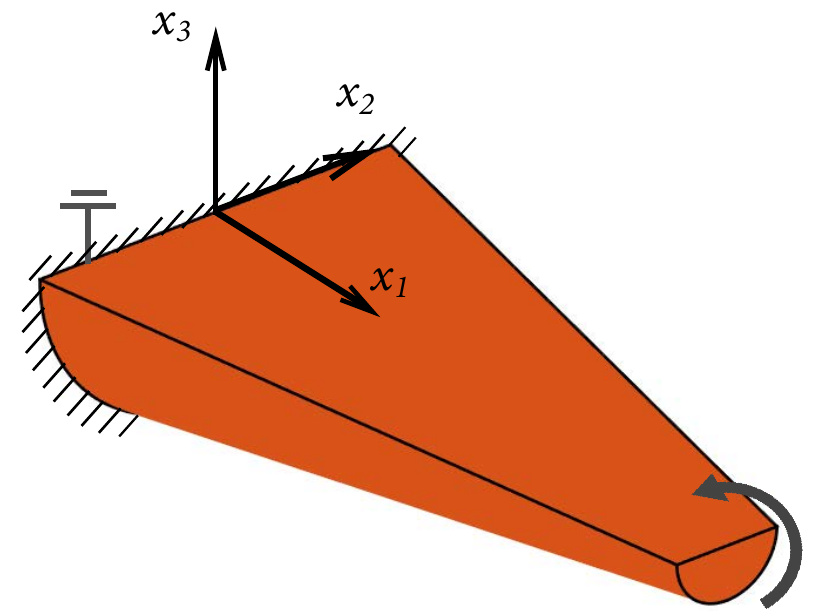}
		\label{fig_torsion_geom}
	\end{subfigure}\hspace{2em}
	\begin{subfigure}[t]{.4\textwidth}\centering
		\captionof{figure}{}\vspace{3.2em}
		\includegraphics[width=\textwidth]{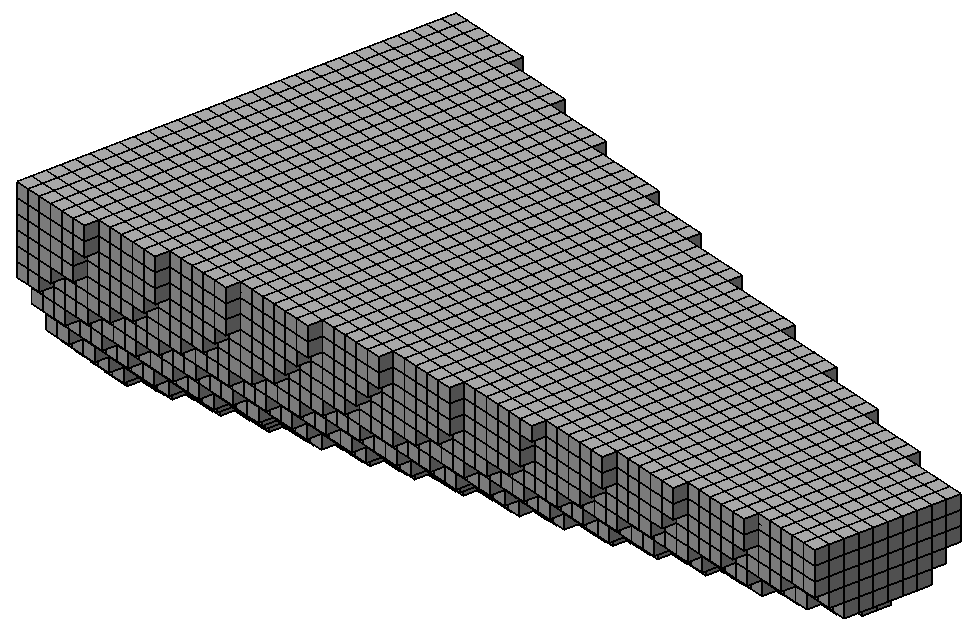}
		\label{fig_torsion_mesh}
	\end{subfigure}%
	\caption{Conical semicircular rod. \emph{a)} Geometrical 3D model, \emph{b)} Unfitted uniform B-spline mesh.}
	\label{fig_torsion}
\end{figure}
\begin{figure}[p]\centering
	\begin{subfigure}[c]{.35\textwidth}\centering
		\captionof{figure}{Shear flexoelectric coupling}
		\includegraphics[width=\textwidth]{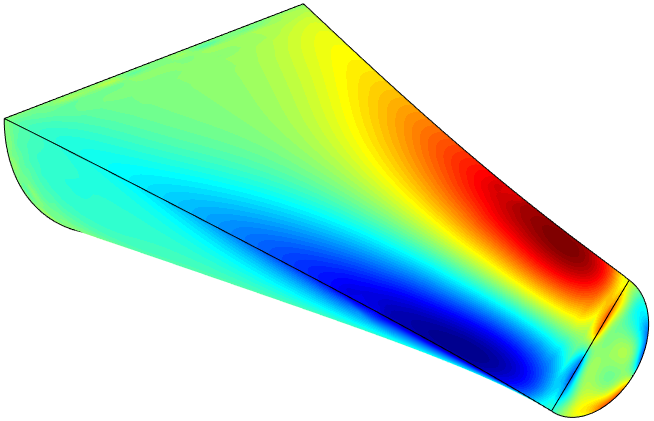}
		\label{fig_torsion_SHEAR}
	\end{subfigure}\hspace{2em}
	\begin{subfigure}[c]{.03\textwidth}\centering\vspace{-1em}
		\includegraphics[height=12em]{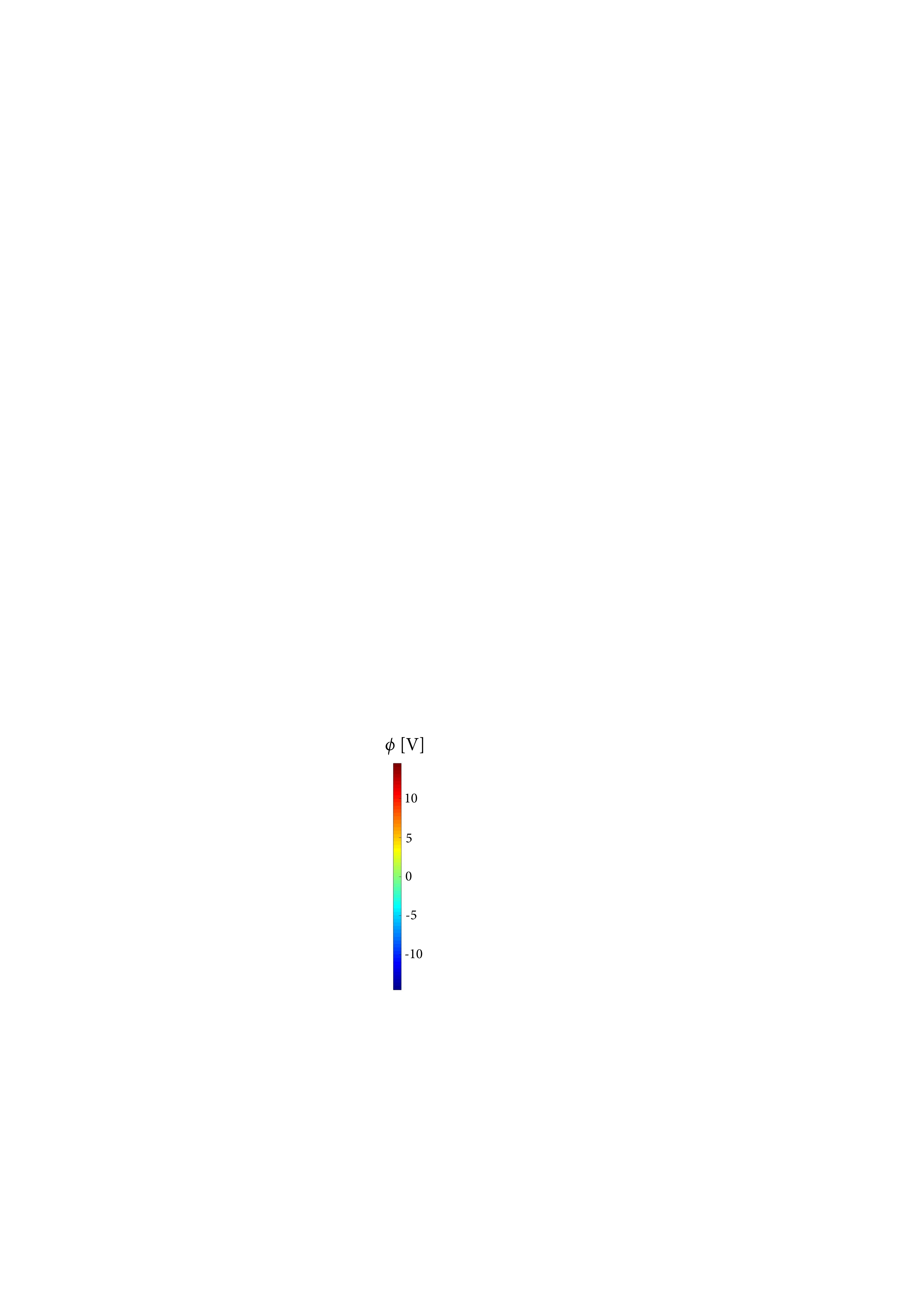}
	\end{subfigure}\hspace{2em}
	\begin{subfigure}[c]{.35\textwidth}\centering
		\captionof{figure}{Full flexoelectric coupling}
		\includegraphics[width=\textwidth]{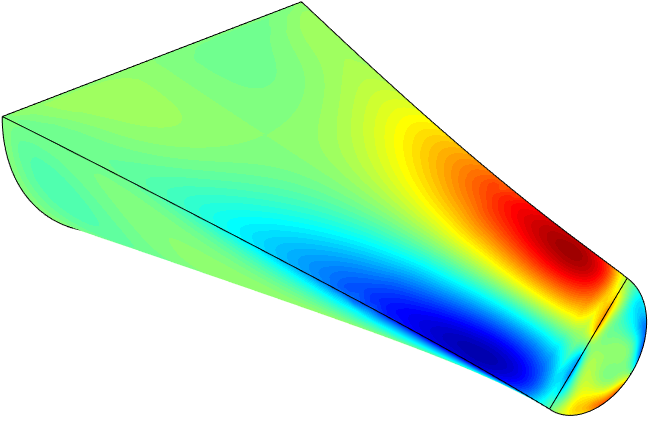}
		\label{fig_torsion_ALL}
	\end{subfigure}
	\caption{Numerical results of the conical semicircular rod. Electric potential over deformed shape of the rod ($\times10$ magnification) after torsion. \emph{a)} Shear flexoelectric coupling considered, \emph{b)} full flexoelectric coupling considered.}
	\label{fig_torsion_Result}
\end{figure}
\begin{figure}[p]\centering
	\begin{subfigure}[c]{.35\textwidth}\centering
		\captionof{figure}{Shear flexoelectric coupling}
		\includegraphics[width=\textwidth]{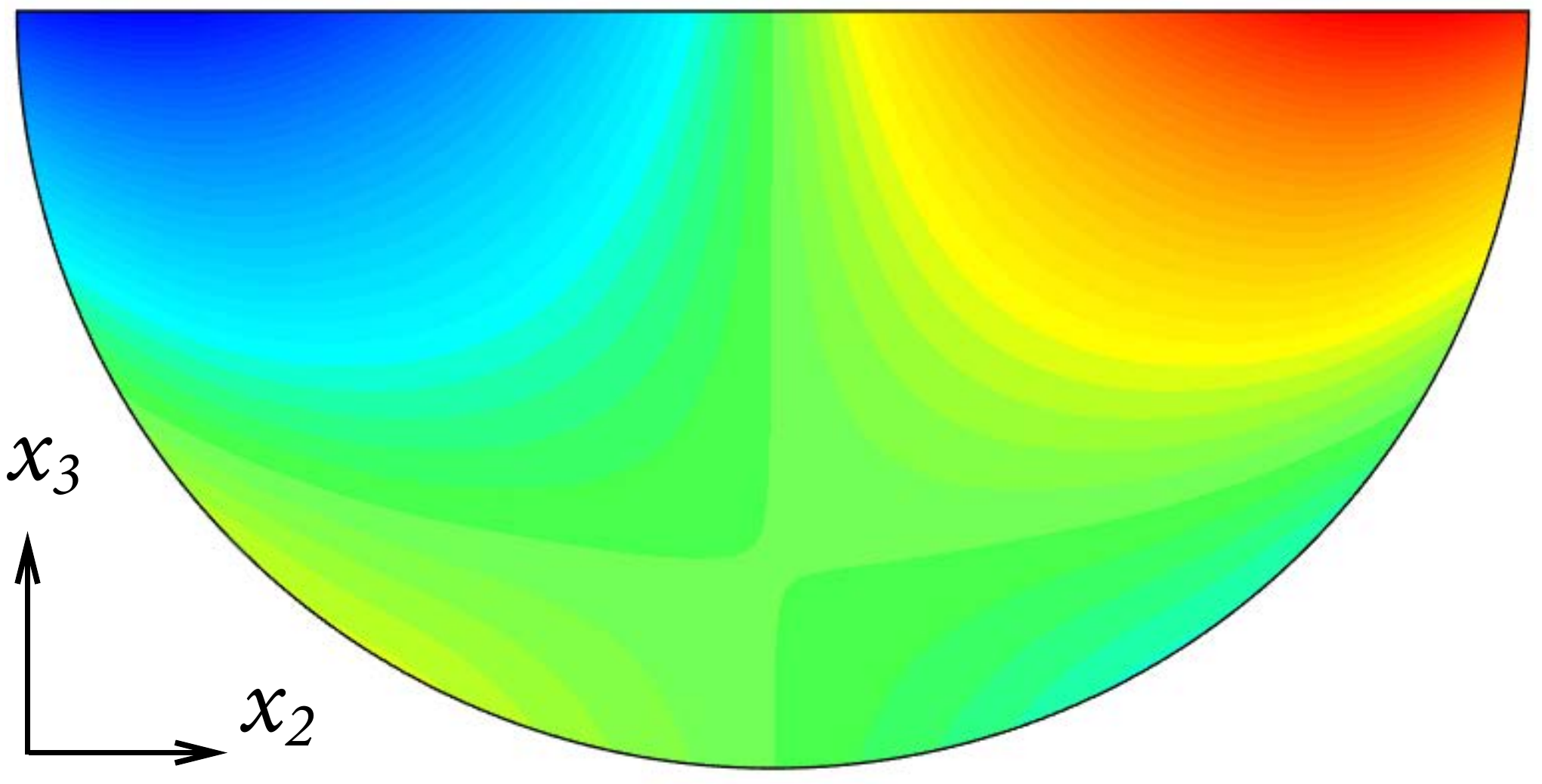}
		\label{fig_torsion_CrossSec_SHEAR}
	\end{subfigure}\hspace{2em}
	\begin{subfigure}[c]{.03\textwidth}\centering\vspace{-1em}
		\includegraphics[height=10em]{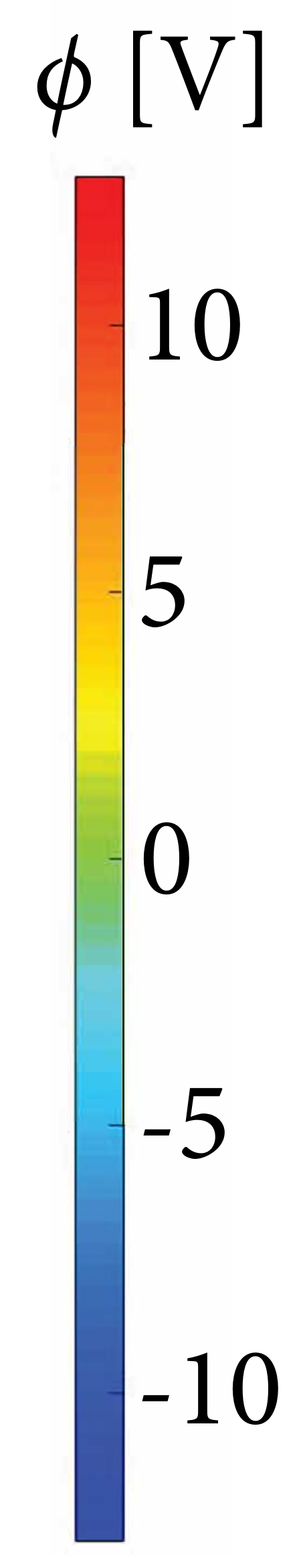}
	\end{subfigure}\hspace{2em}
	\begin{subfigure}[c]{.35\textwidth}\centering
		\captionof{figure}{Full flexoelectric coupling}
		\includegraphics[width=\textwidth]{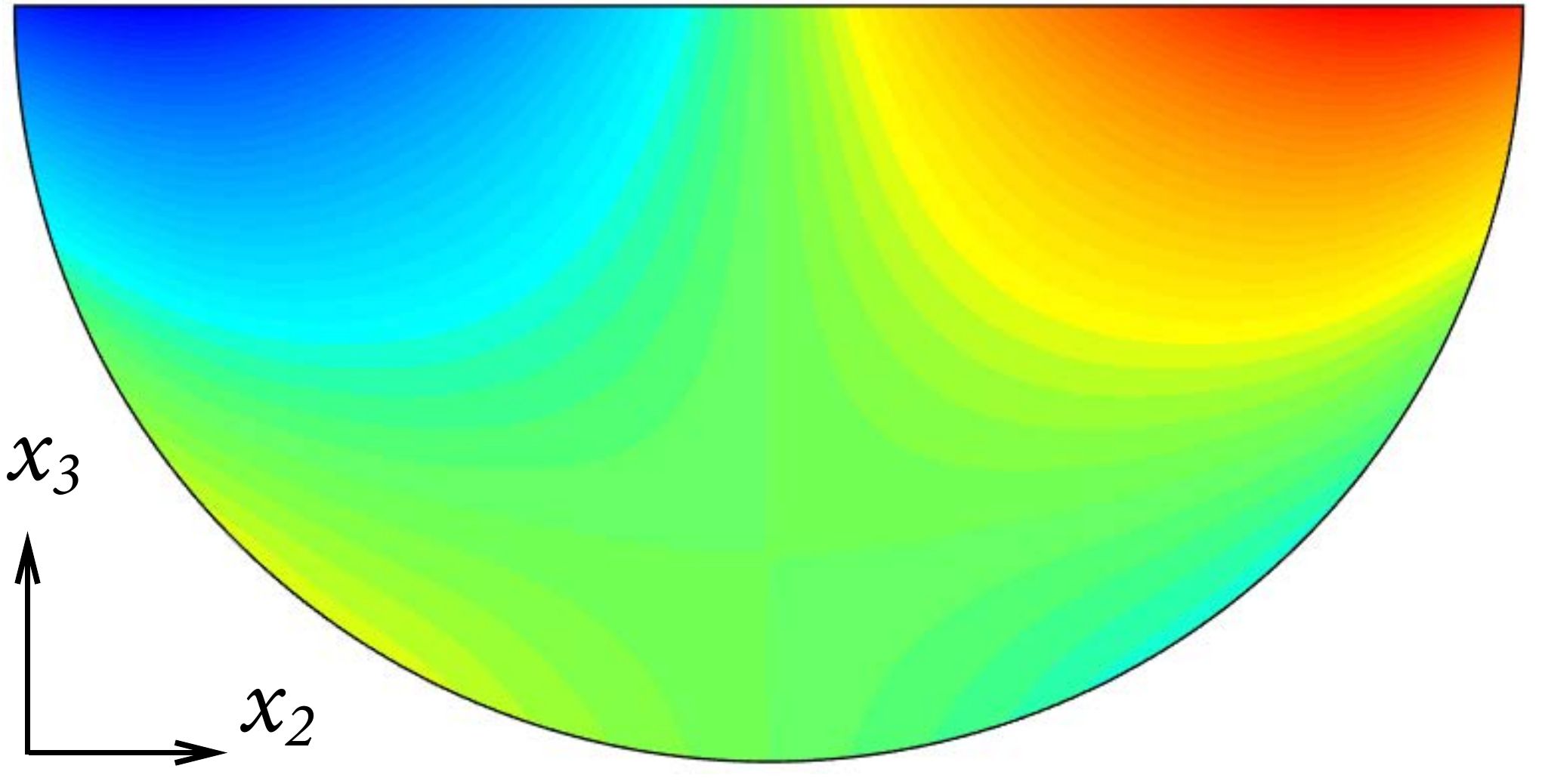}
		\label{fig_torsion_CrossSec_ALL}
	\end{subfigure}
	\caption{Numerical results of the conical semicircular rod. Electric potential distribution in the cross section $x_1=\SI{90.5}{\um}$. \emph{a)} Shear flexoelectric coupling considered, \emph{b)} full flexoelectric coupling considered.}
	\label{fig_torsion_CrossSec}
\end{figure}
Figure \ref{fig_torsion_geom} shows the geometrical model of the conical semicircular rod, with a length of $100\um$. The radii of the semicircular bases are $26.3\um$ and $7.5\um$, and their centers are located at $\x_O=(0,0,0)\um$ and $\x_o=(100,0,0)\um$.

The larger semicircular basis is clamped and grounded, and
torsion is enforced at the opposite basis by prescribing the displacement field. 
The corresponding boundary conditions are:
\begin{subequations}\begin{align}
	u_1=u_2=u_3=0 &\quad\text{at}\quad x_1=0&\text{(larger basis and its perimeter)}, \\
	u_2=-\alpha x_3 &\quad\text{at}\quad x_1=100\um&\text{(smaller basis and its perimeter)}, \\
	u_3=\alpha x_2 &\quad\text{at}\quad x_1=100\um&\text{(smaller basis and its perimeter)}, \\
	\phi=0 &\quad\text{at}\quad x_1=0;
\end{align}\end{subequations}
where $\alpha$ is the tangent of the prescribed torsional angle.

The mechanical response of the rod is composed by several effects, including non-constant twisting (in-plane rotation) and warping (out-of-plane displacement). Without going into the details, one can think of the rod undergoing $\strain_{12}$ and $\strain_{13}$ shear strains varying along the $x_1$ direction, hence triggering the shear flexoelectric effect along the $x_2-x_3$ planes.

Numerical simulations are performed with a cubic $(p=3)$ trivariate B-spline basis on an unfitted uniform mesh of cell size $1.778\um$ (see \fig\ref{fig_torsion_mesh}). The prescribed torsion is set to $\alpha=0.1$, which corresponds to a counterclockwise torsion of about \SI{5.7}{\degree}. The material constants are set to match those of barium strontium titanate, a strongly flexoelectric ceramic, in its paraelectric phase:
\begin{align*}
	\piezo_L=\piezo_T=\piezo_S=0;
	\quad
	E=\SI{152}{\GPa};
	\quad
	\nu=0.33;
	\quad
	l=10\um;
	\\
	\dielec_L=\SI{11}{\nano\joule\per\square\volt\per\meter};
	\quad
	\flexo_L=\flexo_T=\flexo_S=\SI{121}{\micro\joule\per\volt\per\meter}.
\end{align*}
In order to isolate the shear component of the flexoelectric effect, two simulations are performed. In the first one, only shear flexoelectricity is taken into account, namely $\flexo_L=\flexo_T=0$. In the second one, the complete flexoelectricity tensor $\Flexo$ is considered.

Numerical results are shown in \fig\ref{fig_torsion_Result}. The electric potential takes positive values for $x_2>0$ and negative values otherwise, being more prominent near the free end. An effective electric field arises in the \emph{polar} direction contained in the $x_2-x_3$ plane \cite{MocciShearUnpublished}, which can be readily seen by plotting the electric potential in a cross section of the rod (see \fig\ref{fig_torsion_CrossSec}). This distribution allows us to measure the electric potential difference between both sides of the rod, and therefore can be used to measure the shear flexoelectric coefficient \cite{MocciShearUnpublished}.

Results do not vary much by considering or disregarding the longitudinal and transversal coefficients of the flexoelectric tensor (see Subfig. (a) and (b) in \fig\ref{fig_torsion_Result} and \ref{fig_torsion_CrossSec}).
In order to quantify it, the voltage difference at the corners of the cross sections in \fig\ref{fig_torsion_CrossSec}, namely at $\x_{+}=(90.5,9.265,0)\um$ and $\x_{-}=(90.5,-9.265,0)\um$, is evaluated for the two cases, yielding:
\begin{align*}
\phi^\text{Shear}(\x_{+})-\phi^\text{Shear}(\x_{-})=11.45\si{\volt}-(-11.46)\si{\volt}&=22.92\si{\volt},
\\
\phi^\text{Full}(\x_{+})-\phi^\text{Full}(\x_{-})=11.13\si{\volt}-(-11.13)\si{\volt}&=22.28\si{\volt};
\end{align*}
which shows that considering the longitudinal and transversal coefficients of the flexoelectric tensor affects the voltage difference only by 2.87\%.
Therefore, it is apparent that the flexoelectric behavior of this setup is mainly controlled by the shear flexoelectric coefficient $\mu_S$.
\section{Concluding remarks}\label{sec_06}
We have developed a computational approach with unfitted meshes to simulate the electromechanical response of small scale dielectrics at infinitesimal strains, including the piezoelectric and flexoelectric electromechanical couplings. The high-order nature of the equations is addressed by smooth B-spline basis functions on a background Cartesian mesh, which can be hierarchically refined to resolve local features. The unfitted nature of the method allows B-spline-based simulations on arbitrary domain shapes, which can be explicitly represented by means of NURBS surfaces (3D) or curves (2D), and exactly integrated by means of the NEFEM mapping. Therefore, the method is suitable for the rational design and optimization of nanoscale electromechanical devices with no geometrical limitations. The formulation has been detailed for infinitesimal deformations. The ideas described here are currently being extended to a finite deformation setting suitable for the study of flexoelectricity in soft materials as well \cite{CodonyNonlinearFlexoUnpublished}.

Our work highlights two features of the flexoelectric formulation that have been scarcely commented in the literature, namely \emph{a)} the correct way of considering non-smooth boundaries (\ie corners in 2D) and their corresponding non-local boundary conditions, and \emph{b)} the role that the curvature of the boundary plays in physical quantities such as the mechanical tractions.

The Nitsche's method has been particularized to the flexoelectric problem to weakly enforce essential boundary conditions and sensing electrode conditions. Optimal high-order error convergence rates are reported on non-trivial geometries featuring both curved boundaries and corners.

We have simulated several electromechanical setups that are traditionally used to quantify the longitudinal and transversal components of the cubic flexoelectric tensor and are standard benchmarks in the literature of linear flexoelectricity, such as the bending of a cantilever beam and the compression of a truncated pyramid. In all cases, our results match the ones in the literature. Additionally, we perform a 3D simulation of the torsion of a conical semicircular rod, to illustrate the ability of the proposed method to deal accurately with complex geometries. The results are in excellent agreement with \cite{MocciShearUnpublished}.  

The proposed computational framework can assist the design and optimization of a new generation of nanoscale electromechanical devices such as actuators, sensors and energy harvesters, allowing any complexity on the domain shape.
	
\section*{Acknowledgments}
This work was supported by the Generalitat de Catalunya  (``ICREA Academia'' award for excellence in research to I.A., and Grant No.~2017-SGR-1278), and the European Research Council (StG-679451 to I.A.).

\appendix
\section{Particular expression for the traction}\label{sec_app01}
Equation \eqref{eq_Traction} in Section \ref{sec_02} presents a particular expression for the traction $\Traction(\Displacement,\phi)$, which is not the usual one in strain-gradient elasticity or flexoelectricity formulations in the literature \cite{Mindlin1964,Mindlin1968a,Mindlin1968b,Abdollahi2014,Amanatidou2002,Mao2014,Mao2016,Sharma2010,Shen2010}, which is:
\begin{equation}\label{eq_app1}
t_i(\Displacement,\phi)= \left(\cauchyStress_{ij}(\Displacement,\phi)-\hyperStress_{ijk,k}(\Displacement,\phi)+\surfaceGradient_l(n_l)\hyperStress_{ijk}(\Displacement,\phi)n_k\right)n_j-\surfaceGradient_{j}\left(\hyperStress_{ijk}(\Displacement,\phi)n_k\right).
\end{equation}
We obtain \eq\eqref{eq_Traction} by expanding and rearranging terms in \eq\eqref{eq_app1} in the following way:
\begin{align}\label{eq_app2}
t_i(\Displacement,\phi) =&\ \left(\cauchyStress_{ij}(\Displacement,\phi)-\hyperStress_{ijk,k}(\Displacement,\phi)\right)n_j+\surfaceGradient_l(n_l)\hyperStress_{ijk}(\Displacement,\phi)n_jn_k\nonumber\\&-\surfaceGradient_{j}\left(\hyperStress_{ijk}(\Displacement,\phi)\right)n_k-\hyperStress_{ijk}(\Displacement,\phi)\surfaceGradient_{j}(n_k)
\nonumber\\{}=&\
\left(\cauchyStress_{ij}(\Displacement,\phi)-\hyperStress_{ijk,k}(\Displacement,\phi)\right)n_j+\hyperStress_{ijk}(\Displacement,\phi)\left(\surfaceGradient_l(n_l)n_jn_k-\surfaceGradient_{j}(n_k)\right)\nonumber\\&-\surfaceGradient_{k}\left(\hyperStress_{ikj}(\Displacement,\phi)\right)n_j
\nonumber\\{}=&\
\left(\cauchyStress_{ij}(\Displacement,\phi)-\hyperStress_{ijk,k}(\Displacement,\phi)-\surfaceGradient_{k}\hyperStress_{ikj}(\Displacement,\phi)\right)n_j+\hyperStress_{ijk}(\Displacement,\phi)\left(\shapeOperator_{jk}-2Hn_jn_k\right)
\nonumber\\{}=&\
\left(\cauchyStress_{ij}(\Displacement,\phi)-\hyperStress_{ijk,k}(\Displacement,\phi)-\surfaceGradient_{k}\hyperStress_{ikj}(\Displacement,\phi)\right)n_j+\hyperStress_{ijk}(\Displacement,\phi)\curvatureProjector_{jk};
\end{align}
where $\ShapeOperator$ is the \emph{shape operator} of a surface and $H=\frac{1}{2}\trace{\ShapeOperator}$ its \emph{mean curvature} \cite{Block1997}, as defined in Section \ref{sec_02}.

Equation \eqref{eq_app2} reveals that the traction $\Traction(\Displacement,\phi)$ has two contributions: the former involves stress measures, \ie $\CauchyStress(\Displacement,\phi)$ and $\divergence\HyperStress(\Displacement,\phi)$, dotted with the normal vector $\n$ (a first-order measure of the geometry), whereas the latter involves the double stress measure $\HyperStress(\Displacement,\phi)$ dotted with second-order geometry measures, namely the \emph{second-order geometry tensor} defined as $\CurvatureProjector:=\ShapeOperator-2H\n\otimes\n$. Thus, it is clear that high-order physics are intrinsically linked to high-order geometrical measures of the domain, such as the curvature of its boundary.
\section{Material tensors}\label{sec_app03}
In the following, $n_d$ refers to the number of dimensions of the physical space $\Omega$ (either 2 or 3).
Material tensors are described component-wise, and only the non-zero components are specified.
% % % %
%\subsection{Elasticity}
% % % % 

Isotropic elasticity is represented by the fourth-order tensor $\Elast$, which depends on the \emph{Young modulus} $E$ and the \emph{Poisson ratio} $\nu$ as
\begin{align}
\elast_{iiii}&=C_L,&&& i&=1,\dots,n_d;\nonumber\\
\elast_{iijj}&=C_T,&&& i,j&=1,\dots,n_d\quad\emph{such that}\quad i\neq j;\nonumber\\
\elast_{ijij}=\elast_{ijji}&=C_S,&&& i,j&=1,\dots,n_d\quad\emph{such that}\quad i\neq j,
\end{align}
where the parameters $C_L$, $C_S$ and $C_T$ are 
\begin{align}\label{eq_c3}
C_L\coloneqq\frac{E\left(1-\nu\right)}{(1+\nu)(1-2\nu)},&&
C_T\coloneqq\frac{E\nu}{(1+\nu)(1-2\nu)},&&
C_S\coloneqq\frac{E}{2(1+\nu)}
\end{align}
in the \emph{3D} and \emph{plane strain 2D} cases, and
\begin{align}\label{eq_c2}
C_L\coloneqq\frac{E}{1-\nu^2},&&
C_T\coloneqq\frac{E\nu}{1-\nu^2},&&
C_S\coloneqq\frac{E}{2(1+\nu)}
\end{align}
for the \emph{plane stress 2D} case.

% % % %
%\subsection{Strain gradient elasticity}
% % % % 
We consider an isotropic simplified strain gradient elasticity model \cite{Altan1997}, which is a particular case of the general model of strain gradient elasticity in \cite{Mindlin1968a}.
The strain gradient elasticity tensor is represented by the sixth-order tensor $\StrGr$, which depends on the \emph{Young modulus} $E$, the \emph{Poisson ratio} $\nu$ and a single internal length scale $l$ in the following form:
\begin{align}
	\strGr_{iikiik}&=l^2C_L,&&& i,k&=1,\dots,n_d;\nonumber\\
	\strGr_{iikjjk}&=l^2C_T,&&& i,j,k&=1,\dots,n_d\quad\emph{such that}\quad i\neq j;\nonumber\\
	\strGr_{ijkijk}=\strGr_{ijkjik}&=l^2C_S,&&& i,j,k&=1,\dots,n_d\quad\emph{such that}\quad i\neq j
\end{align}
with parameters $C_L$, $C_S$ and $C_T$ defined in \eq\eqref{eq_c3} and \eqref{eq_c2}.

% % % %
%\subsection{Dielectricity}
% % % % 
Isotropic dielectricity is represented by the second-order tensor $\Dielec$, which depends on the parameter $\dielec_L$ as 
\begin{align}
	\dielec_{ii}&=\dielec_L,&&& i&=1,\dots,n_d.
\end{align}

% % % %
%\subsection{Piezoelectricity}
% % % % 
Piezoelectricity is represented by the third-order tensor $\Piezo$.
%Two symmetry cases are considered, i.e. cubic and tetragonal symmetry.
%
%Cubic symmetry involves an isotropic shear coupling represented by the parameter $\piezo_S$. It is defined only for $n_d=3$ as
%\begin{align}
%	\piezo_{ijk}&=\piezo_S,&&& i,j,k&=1,\dots,3\quad\emph{such that}\quad i\neq j\neq k\neq i.
%\end{align}
%
%A more general symmetry is the tetragonal one, 
Tetragonal symmetry is considered,
which has a principal direction and involves longitudinal, transversal and shear couplings represented by the parameters $\piezo_L$, $\piezo_T$ and $\piezo_S$, respectively. For a material with principal direction $\x_1$, the piezoelectric tensor $\Piezo^{<\x_1>}$ reads
\begin{align}
	{\piezo^{<\x_1>}}_{111}&=\piezo_L;\nonumber\\
	{\piezo^{<\x_1>}}_{1jj}&=\piezo_T,&&& j&=2,\dots,n_d;\nonumber\\
	{\piezo^{<\x_1>}}_{j1j}={\piezo^{<\x_1>}}_{jj1}&=\piezo_S,&&& j&=2,\dots,n_d.
\end{align}
The piezoelectric tensor $\Piezo$ oriented in an arbitrary direction $\toVect{d}$ is obtained by rotating $\Piezo^{<\x_1>}$ as
\begin{equation}
\piezo_{lij}=R_{lL}R_{iI}R_{jJ}{\piezo^{<\x_1>}}_{LIJ},
\end{equation}
where $\toMat{R}$ is a rotation matrix that rotates $\x_1$ onto $\toVect{d}$.

% % % %
%\subsection{Flexoelectricity}
% % % % 
Flexoelectricity is represented by the fourth-order tensor $\Flexo$. Cubic symmetry is considered, which leads to a flexoelectric tensor involving longitudinal, transversal and shear couplings represented by the parameters $\flexo_L$, $\flexo_T$ and $\flexo_S$, respectively. We refer to \cite{LeQuang2011} for an extensive analysis of other possible symmetries for the flexoelectric tensor. The components of the flexoelectric tensor $\Flexo^{<{\x}>}$ of a material oriented in the Cartesian axes are the following:
\begin{align}
	{\flexo^{<{\x}>}}_{iiii}&=\flexo_L,&&& i&=1,\dots,n_d\nonumber;\\
	{\flexo^{<{\x}>}}_{ijji}&=\flexo_T,&&& i,j&=1,\dots,n_d\quad\emph{such that}\quad i\neq j;\nonumber\\
    {\flexo^{<{\x}>}}_{iijj}={\flexo^{<{\x}>}}_{ijij}&=\flexo_S,&&& i,j&=1,\dots,n_d\quad\emph{such that}\quad i\neq j.
\end{align}
The flexoelectric tensor $\Flexo$ oriented in an arbitrary orthonormal basis is obtained by rotating $\Flexo^{<\x>}$ as
\begin{equation}
\flexo_{lijk}=R_{lL}R_{iI}R_{jJ}R_{kK}{\flexo^{<\x>}}_{LIJK},
\end{equation}
where $\toMat{R}$ is the rotation matrix that rotates the Cartesian basis to the desired orthonormal basis.

Assuming the material models presented in this Appendix, the restrictions on material tensors in \eq\eqref{restrictions} simplify to restrictions on material coefficients as follows:
\begin{align}
 \dielec_L, C_L, C_S >  0                        , \quad
 l \geq 0
 %,                                       \quad
%|\piezo_L|,|\piezo_T|   <  (C_L\dielec_L)^{1/2}, \quad
%|\piezo_S|              <  (C_S\dielec_L)^{1/2}, \quad
%|\flexo_L|,|\flexo_T|   < l(C_L\dielec_L)^{1/2}, \quad
%|\flexo_S|              < l(C_S\dielec_L)^{1/2}
.
\end{align}
\section{Computation of penalty parameters}\label{sec_app02}
Nitsche's method involves numerical penalty-like parameters to enforce Dirichlet boundary conditions. The stability of the formulation is guaranteed for large values of the parameters. Too large values would lead to ill-conditioning of the system matrix but, differently to penalty methods, usually moderate values provide good results.
Lower bounds of the penalty parameters can be assessed globally (constant for the whole mesh) \cite{Embar2010} or locally (cell-wise) \cite{Griebel2003}, being the latter more appealing due to \emph{i)} lower condition number of the resulting algebraic system and \emph{ii)} lower computational cost \cite{dePrenter2016}.

Typically, for elliptic PDE, lower bounds are found by studying the coercivity of the formulation (see for instance \cite{Embar2010}). In the case of flexoelectricity, it is a saddle point problem corresponding to the coupling between mechanical (positive definite) and electrical (negative definite) problems. Therefore, coercivity is met by checking the positivity and negativity of the second variations of the energy functional with respect to the mechanical and electrical unknowns, respectively, as stated in \eq\eqref{eq_2var}:
\begin{equation}\label{eq_coerc}
\vvar[\Displacement]{\Pi}[\delta\Displacement]>0,~\forall\delta\Displacement\in\mathcal{U}
;\qquad
\vvar[\phi]{\Pi}[\delta\phi]<0,~\forall\delta\phi\in\mathcal{P}
.
\end{equation}
From \eq\eqref{eq_coerc}, one can readily see that the second variations of the energy functional depend solely on $\delta\Displacement$ \emph{or} $\delta\phi$. Therefore, the coupling physics (\ie piezoelectricity and flexoelectricity) do not play any role in determining the coercivity of the formulation. As a consequence, the coercivity analysis can be performed for the uncoupled mechanical and electrical problems independently, namely
\begin{equation}\label{eq_coerc_sep}
\vvar[\Displacement]{\Pi^\text{Mechanical}}[\delta\Displacement]>0,~\forall\delta\Displacement\in\mathcal{U}
;\qquad
\vvar[\phi]{\Pi^\text{Electrical}}[\delta\phi]<0,~\forall\delta\phi\in\mathcal{P}
,
\end{equation}
being $\Pi^\text{Mechanical}\coloneqq\Pi|_{\delta\phi=0}$ and $\Pi^\text{Electrical}\coloneqq\Pi|_{\delta\Displacement=\toVect{0}}$.

The mechanical part of the problem corresponds to a strain-gradient elasticity formulation, for which we derive next the lower bounds arising from the coercivity analysis.
The electrical part of the problem corresponds to a (negative) Poisson equation, which has been largely studied in the literature \cite{Embar2010} and whose corresponding lower bounds are well known and not presented here.
%%%
\subsection{Deriving conditions for coercivity of the strain-gradient elasticity formulation}
%%%
The bilinear form of the strain-gradient elasticity formulation corresponds to the mechanical part of the flexoelectric bilinear form in \eq\eqref{bil}:
\begin{multline}
\mathcal{B}[\Displacement,\delta\Displacement]=
	\int_\Omega\Big(\cauchyStress_{ij}(\Displacement)\strain_{ij}(\delta\Displacement)+\hyperStress_{ijk}(\Displacement)\strain_{ij,k}(\delta\Displacement)\Big)\mathrm{d}\Omega
	+\int_{\partial\Omega_u}\Big(\Big(\beta_uu_{i}-t_i(\Displacement)\Big)\delta u_{i}-u_{i}t_i(\delta \Displacement)\Big)\dd\Gamma
	+{}\\{}+
	\int_{\partial\Omega_v}\Big(\Big(\beta_v\partial^n(u_i)-r_i(\Displacement)\Big)\partial^n(\delta u_i)-\partial^n(u_i)r_i(\delta \Displacement)\Big)\dd\Gamma
	+\int_{C_u}\Big(\Big(\beta_{C_u}u_i-j_i(\Displacement)\Big)\delta u_i-u_i j_i(\delta\Displacement) \Big)\dd\ds,
	\qquad
	\\[-0.5em]\text{ with }\quad
	\beta_u,\beta_v,\beta_{C_u},\in\mathds{R}^+.
\end{multline}
It is easy to verify that the second variation $\vvar[\Displacement]{\Pi^\text{Mechanical}}[\delta\Displacement]$ corresponds to $\mathcal{B}[\delta\Displacement,\delta\Displacement]$.
Since $\mathcal{B}[\Displacement,\delta\Displacement]$ can be expressed as $\mathcal{B}[\Displacement,\delta\Displacement]=\sum_c\mathcal{B}^{\Omega_\square^c\cap\Omega}[\Displacement,\delta\Displacement]$, coercivity of $\mathcal{B}[\Displacement,\delta\Displacement]$ is cell-wise met by proving coercivity of $\mathcal{B}^{\Omega_\square^c\cap\Omega}[\Displacement,\delta\Displacement]$, $\forall\Omega_\square^c\in\mathcal{I}\cup\mathcal{C}$.
Therefore, a sufficient cell-wise coercivity condition is:
\begin{equation}\label{eq_coer}
\mathcal{B}^{\Omega_\square^c\cap\Omega}[\delta\Displacement,\delta\Displacement]>0
%\gamma^c\int_{\Omega_\square^c\cap\Omega}\left(\delta\cauchyStress_{ij}\delta\strain_{ij}
%+
%\delta\hyperStress_{ijk}\delta\strain_{ij,k}\right)\mathrm{d}\Omega
,\qquad \forall\delta\Displacement\neq 0,\quad \forall\Omega_\square^c\in\mathcal{I}\cup\mathcal{C}.
\end{equation}
Since \eq\eqref{eq_coer} is trivially fulfilled on inner cells, the coercivity analysis is carried out only on cut cells:
\begin{multline}
\mathcal{B}^{\Omega_\square^c\cap\Omega}[\delta\Displacement,\delta\Displacement]=
\int_{\Omega_\square^c\cap\Omega}\Big(\cauchyStress_{ij}(\delta\Displacement)\strain_{ij}(\delta\Displacement)+\hyperStress_{ijk}(\delta\Displacement)\strain_{ij,k}(\delta\Displacement)\Big)\mathrm{d}\Omega
+{}\\{}+
\int_{\Omega_\square^c\cap\partial\Omega_u}\Big(\beta^c_u[\delta u_{i}]^2-2 t_i(\delta\Displacement) \delta u_{i}\Big)\dd\Gamma
+{}\\{}+
\int_{\Omega_\square^c\cap\partial\Omega_v}\Big(\beta^c_v[\partial^n( \delta u_i)]^2-2 r_i(\delta\Displacement)\partial^n(\delta u_i)\Big)\dd\Gamma
{}\\{}+
\int_{\Omega_\square^c\cap C_u}\Big(\beta^c_{C_u}[ \delta u_i]^2-2 j_i(\delta\Displacement) \delta u_i\Big)\dd\ds>0
%\qquad \text{ with }
%\beta^c_u,\beta^c_v,\beta^c_{C_u},\in\mathds{R}^+;
\end{multline}
where the penalty parameters $\beta^c_u,\beta^c_v,\beta^c_{C_u}\in\mathds{R}^+$ correspond to each cell $\Omega_\square^c$.

Applying the \emph{Cauchy-Schwarz inequality} and the \emph{Young's inequality} leads to:
\begin{multline}\label{eq_schwarz}
	\mathcal{B}^{\Omega_\square^c\cap\Omega}[\delta\Displacement,\delta\Displacement]
	\geq
	\int_{\Omega_\square^c\cap\Omega}\left(\cauchyStress_{ij}(\delta\Displacement)\strain_{ij}(\delta\Displacement)+\hyperStress_{ijk}(\delta\Displacement)\strain_{ij,k}(\delta\Displacement)\right)\mathrm{d}\Omega
	-{}\\{}-
	\normLL[\Omega_\square^c\cap\partial\Omega_u]{ \Traction(\delta\Displacement)}\normLL[\Omega_\square^c\cap\partial\Omega_u]{\delta \Displacement}
	-{}\\{}-
	\normLL[\Omega_\square^c\cap\partial\Omega_v]{ \HighOrderTraction(\delta\Displacement)}\normLL[\Omega_\square^c\cap\partial\Omega_v]{\partial^n(\delta\Displacement)}
	-
	\normLL[\Omega_\square^c\cap\partial\Omega_u]{ \HighOrderJump(\delta\Displacement)}\normLL[\Omega_\square^c\cap C_u]{\delta \Displacement}
	+{}\\{}+
	\beta^c_u\normLL[\Omega_\square^c\cap\partial\Omega_u]{\delta \Displacement}^2
	+
	\beta^c_v\normLL[\Omega_\square^c\cap\partial\Omega_v]{\partial^n(\delta\Displacement)}^2\mathrm{d}\Gamma
	+
	\beta^c_{C_u}\normLL[\Omega_\square^c\cap C_u]{\delta \Displacement}^2
	\geq{}\\[1.5em]{}\geq
	\int_{\Omega_\square^c\cap\Omega}\left(\cauchyStress_{ij}(\delta\Displacement)\strain_{ij}(\delta\Displacement)
	+
	\hyperStress_{ijk}(\delta\Displacement)\strain_{ij,k}(\delta\Displacement)\right)\mathrm{d}\Omega
	-
	\frac{1}{\epsilon_u}\normLL[\Omega_\square^c\cap\partial\Omega_u]{ \Traction(\delta\Displacement)}^2
	-
	\epsilon_u\normLL[\Omega_\square^c\cap\partial\Omega_u]{\delta \Displacement}^2
	-{}\\{}-
	\frac{1}{\epsilon_v}\normLL[\Omega_\square^c\cap\partial\Omega_v]{\HighOrderTraction(\delta\Displacement)}^2
	-
	\epsilon_v\normLL[\Omega_\square^c\cap\partial\Omega_v]{\partial^n(\delta\Displacement)}^2
	-
	\frac{1}{\epsilon_{C_u}}\normLL[\Omega_\square^c\cap C_u]{\HighOrderJump(\delta\Displacement)}^2
	-
	\epsilon_{C_u}\normLL[\Omega_\square^c\cap C_u]{\delta \Displacement}^2
	+{}\\{}+
	\beta^c_u\normLL[\Omega_\square^c\cap\partial\Omega_u]{\delta \Displacement}^2
	+
	\beta^c_v\normLL[\Omega_\square^c\cap\partial\Omega_v]{\partial^n(\delta\Displacement)}^2
	+
	\beta^c_{C_u}\normLL[\Omega_\square^c\cap C_u]{\delta \Displacement}^2
	={}\\[1.5em]{}=
	\int_{\Omega_\square^c\cap\Omega}\left(\cauchyStress_{ij}(\delta\Displacement)\strain_{ij}(\delta\Displacement)
	+
	\hyperStress_{ijk}(\delta\Displacement)\strain_{ij,k}(\delta\Displacement)\right)\mathrm{d}\Omega
	-{}\\{}-
	\frac{1}{\epsilon_u}\normLL[\Omega_\square^c\cap\partial\Omega_u]{\Traction(\delta\Displacement)}^2
	-
	\frac{1}{\epsilon_v}\normLL[\Omega_\square^c\cap\partial\Omega_v]{\HighOrderTraction(\delta\Displacement)}^2
	-
	\frac{1}{\epsilon_{C_u}}\normLL[\Omega_\square^c\cap C_u]{\HighOrderJump(\delta\Displacement)}^2
	{}\\{}
	+
	(\beta^c_u-\epsilon_u)\normLL[\Omega_\square^c\cap\partial\Omega_u]{\delta \Displacement}^2
	+
	(\beta^c_v-\epsilon_v)\normLL[\Omega_\square^c\cap\partial\Omega_v]{\partial^n(\delta\Displacement)}^2
	+
	(\beta^c_{C_u}-\epsilon_{C_u})\normLL[\Omega_\square^c\cap C_u]{\delta \Displacement}^2>0;
\end{multline}
which holds $\forall\epsilon_u,\epsilon_v,\epsilon_{C_u}\in\mathds{R}^+$.
Let us consider now the mesh-dependent constants $K_u,K_v,K_{C_u}\in\mathds{R}^+$ such that
\begin{subequations}\label{cucv}\begin{align}
		\normLL[\Omega_\square^c\cap\partial\Omega_u]{\Traction(\delta\Displacement)}^2\leq&~ K_u\int_{\Omega_\square^c\cap\Omega}\left(\cauchyStress_{ij}(\delta\Displacement)\strain_{ij}(\delta\Displacement)
		+\hyperStress_{ijk}(\delta\Displacement)\strain_{ij,k}(\delta\Displacement)\right)\mathrm{d}\Omega,
		\\
		\normLL[\Omega_\square^c\cap\partial\Omega_v]{\HighOrderTraction(\delta\Displacement)}^2\leq&~ K_v\int_{\Omega_\square^c\cap\Omega}\left(\cauchyStress_{ij}(\delta\Displacement)\strain_{ij}(\delta\Displacement)
		+\hyperStress_{ijk}(\delta\Displacement)\strain_{ij,k}(\delta\Displacement)\right)\mathrm{d}\Omega,
		\\
		\normLL[\Omega_\square^c\cap{C_u}]{\HighOrderJump(\delta\Displacement)}^2\leq&~ K_{C_u}\int_{\Omega_\square^c\cap\Omega}\left(\cauchyStress_{ij}(\delta\Displacement)\strain_{ij}(\delta\Displacement)
		+\hyperStress_{ijk}(\delta\Displacement)\strain_{ij,k}(\delta\Displacement)\right)\mathrm{d}\Omega;
\end{align}\end{subequations}
% Note that equations \eqref{cucv} do \emph{not} have the usual form of fourth-order problems as in \cite{Embar2010}. Here, both inequalities are compared at the right hand side to the \emph{same} domain integral, which incorporates both mechanic and hyper-mechanic terms. This procedure provides more flexibility on the choice of $\beta_u$ and $\beta_v$, as we will see in next lines.
for all admissible $\delta\Displacement$, which can be computed as detailed in \ref{secKuKv}.
Then, \eq\eqref{eq_schwarz} leads to:
\begin{multline}\label{eq_schwarz2}
	\mathcal{B}^{\Omega_\square^c\cap\Omega}[\delta\Displacement,\delta\Displacement]
	\geq
	\left(1-\frac{K_u}{\epsilon_u}-\frac{K_v}{\epsilon_v}-\frac{K_{C_u}}{\epsilon_{C_u}}\right)
	\int_{\Omega_\square^c\cap\Omega}\left(\cauchyStress_{ij}(\delta\Displacement)\strain_{ij}(\delta\Displacement)
	+
	\hyperStress_{ijk}(\delta\Displacement)\strain_{ij,k}(\delta\Displacement)\right)\mathrm{d}\Omega
	+{}\\{}+
	(\beta^c_u-\epsilon_u)\normLL[\Omega_\square^c\cap\partial\Omega_u]{\delta \Displacement}^2
	+
	(\beta^c_v-\epsilon_v)\normLL[\Omega_\square^c\cap\partial\Omega_v]{\partial^n(\delta\Displacement)}^2
	+{}\\{}+
	(\beta^c_{C_u}-\epsilon_{C_u})\normLL[\Omega_\square^c\cap{C_u}]{\delta \Displacement}^2>0
	;
	\qquad\forall\epsilon_u,\epsilon_v,\epsilon_{C_u}\in\mathds{R}^+.
\end{multline}
Thus, the conditions 
\begin{align}\label{coercivitycond}
	\frac{K_u}{\epsilon_u}+\frac{K_v}{\epsilon_v}+\frac{K_{C_u}}{\epsilon_{C_u}}<1,&&
	\beta^c_u>\epsilon_u,&&
	\beta^c_v>\epsilon_v,&&
	\beta^c_{C_u}>\epsilon_{C_u};
\end{align}
are sufficient conditions for the coercivity of $\mathcal{B}^{\Omega_\square^c\cap\Omega}[\delta\Displacement,\delta\Displacement]$,  for any $\epsilon_u,\epsilon_v,\epsilon_{C_u}\in\mathds{R}^+$.
%\\[1em]\emph{\textbf{Note:\quad} This analysis can be performed \textbf{globally} \cite{Griebel2003} or \textbf{element-wise} \cite{Embar2010}, considering the elemental discretization of  $\mathcal{B}(\Displacement,\delta\Displacement)=\sum_e\mathcal{B}^e(\Displacement,\delta\Displacement)$. As said in \cite{dePrenter2016} and many other publications, global stabilization negatively affects the conditioning of the method and can degenerate it to a penalty method. Therefore, we use and recommend the element-wise analysis, where condition \eqref{coercivitycond} is to be enforced on every single element, leading to elemental $C_u^e,C_v^e,\beta_u^e$ and $\beta_v^e$ parameters. We just show the global analysis here for the sake of simplicity.}
%%%
\subsection{Conditions on the penalty parameters}
%%%
In order to get explicit bounds of the penalty parameters, let us define $\alpha_1\coloneqq\epsilon_u/\epsilon_v\in\mathds{R}^+$ and $\alpha_2\coloneqq\epsilon_u/\epsilon_{C_u}\in\mathds{R}^+$.
The first coercivity condition in \eq\eqref{coercivitycond} is rewritten as
\begin{equation}\label{coercivitycond2}
\epsilon_u>K_u+\alpha_1 K_v+\alpha_2 K_{C_u};\qquad\forall\epsilon_u,\alpha_1,\alpha_2\in\mathds{R}^+.
\end{equation}
By means of \eq\eqref{coercivitycond2}, we can rewrite the remaining conditions in \eq\eqref{coercivitycond} as
\begin{subequations}\label{coercivitycond3}\begin{align}
	\beta^c_u>\epsilon_u&>C_u+\alpha_1 K_v+\alpha_1 K_{C_u},\\
	\beta^c_v>\epsilon_v=\frac{\epsilon_u}{\alpha_1}&>\frac{1}{\alpha_1}K_u+K_v+\frac{\alpha_2}{\alpha_1}K_{C_u},\\
	\beta^c_{C_u}>\epsilon_{C_u}=\frac{\epsilon_u}{\alpha_2}&>\frac{1}{\alpha_2}K_u+\frac{\alpha_1}{\alpha_2}K_v+K_{C_u}.
\end{align}\end{subequations}
Since \eq\eqref{coercivitycond3} ensure coercivity for any $\epsilon_u\in\mathds{R}^+$, we can express the bounds of the penalty parameters as a function of $\alpha_1$ and $\alpha_2$ only.

The bilinear form $\mathcal{B}^{\Omega_\square^c\cap\Omega}$ is coercive if $\beta_u^c,\beta_v^c,\beta_{C_u}^c$ satisfy
%\begin{subequations}\label{coercivitycond4}\begin{align}
%		\beta^c_u(\alpha_1,\alpha_2)&>C_u+\alpha_1 K_v+\alpha_1 K_{C_u};\\
%		\beta^c_v(\alpha_1,\alpha_2)&>\frac{1}{\alpha_1}K_u+K_v+\frac{\alpha_2}{\alpha_1}K_{C_u};\\
%		\beta^c_{C_u}(\alpha_1,\alpha_2)&>\frac{1}{\alpha_2}K_u+\frac{\alpha_1}{\alpha_2}K_v+K_{C_u};
%\end{align}\end{subequations}
\begin{equation}\label{coercivitycond5}
\begin{bmatrix}\beta^c_u\\\beta^c_v\\\beta^c_{C_u}\end{bmatrix}
>
\begin{bmatrix}
1 & \alpha_1 & \alpha_2 \\ 1/\alpha_1 & 1 & \alpha_2/\alpha_1\\ 1/\alpha_2 & \alpha_1/\alpha_2 & 1
\end{bmatrix}\cdot
\begin{bmatrix}K_u\\K_v\\K_{C_u}\end{bmatrix},
\end{equation}
for any $\alpha_1,\alpha_2\in\mathds{R}^+$.

At this point, it is worth mentioning that \eq\eqref{coercivitycond5} correspond to a particular cell $\Omega_\square^c\in\mathcal{C}$ such that 
\emph{i)} $\Omega_\square^c\cap\partial\Omega_u\neq\emptyset$ (classical Dirichlet boundary), 
\emph{ii)} $\Omega_\square^c\cap\partial\Omega_v\neq\emptyset$ (non-local Dirichlet boundary) and
\emph{iii)} $\Omega_\square^c\cap{C_u}\neq\emptyset$ (non-local Dirichlet edges),
\ie for cut cells with the complete set of Dirichlet boundary conditions.
In the case of a cell with a Neumann boundary condition, it is easy to verify that \eq\eqref{coercivitycond5} still holds if we just consider the bounds and constants corresponding to the Dirichlet boundary conditions. For instance, in a cell where double tractions $\HighOrderTraction$ are prescribed on the boundary, $\Omega_\square^c\cap\partial\Omega_v=\emptyset$ and we obtain the following bounds:

\begin{equation}\label{coercivitycond6}
\begin{bmatrix}\beta^c_u\\\beta^c_{C_u}\end{bmatrix}
>
\begin{bmatrix}
1 & \alpha_2 \\ 1/\alpha_2 & 1
\end{bmatrix}\cdot
\begin{bmatrix}K_u\\K_{C_u}\end{bmatrix},
\qquad\qquad
\forall\alpha_2\in\mathds{R}^+.
\end{equation}
Any choice of $\alpha_1,\alpha_2\in\mathds{R}^+$ leads to a coercive formulation; however, the condition number of the resulting linear system might be affected. Further investigation is required to assess suitable values for $\{\alpha_1,\alpha_2\}$ for providing condition numbers of the system as low as possible.
\subsection{Computation of the mesh-dependent constants $K_u$, $K_v$ and $K_{C_u}$}\label{secKuKv}
%%%
Mesh-dependent constants $K_u$, $K_v$ and $K_{C_u}$ satisfying \eq\eqref{cucv} can be taken as the largest eigenvalues of the following generalized eigenvalue problems that arise from the spatial discretization of \eq\eqref{cucv}:
\begin{subequations}\label{cucv1}\begin{align}
    \toMat{B_u}\cdot\toVect{x_u}&=~ \lambda_u\toMat{V}\cdot\toVect{x_u},\\
	\toMat{B_v}\cdot\toVect{x_v}&=~ \lambda_v\toMat{V}\cdot\toVect{x_v},\\
	\toMat{B_{C_u}}\cdot\toVect{x_{C_u}}&=~ \lambda_{C_u}\toMat{V}\cdot\toVect{x_{C_u}}.	
\end{align}\end{subequations}
The matrices $\mathbf{B_u},\mathbf{B_v},\mathbf{B_{C_u}}$ and $\mathbf{V}$ are defined at each cell $\Omega_\square^c$ as the discretization of the weak forms $\mathcal{B}_u^{\Omega_\square^c\cap\Omega}$, $\mathcal{B}_v^{\Omega_\square^c\cap\Omega}$, $\mathcal{B}_{C_u}^{\Omega_\square^c\cap\Omega}$ and $\mathcal{B}_\Omega^{\Omega_\square^c\cap\Omega}$ in \eq\eqref{cucv}, respectively, with
\begin{subequations}\label{cucv2}\begin{align}
		\mathcal{B}_u^{\Omega_\square^c\cap\Omega}[\Displacement,\displacement]\coloneqq& \int_{\Omega_\square^c\cap\partial\Omega_u} t_{i}(\Displacement)t_{i}(\displacement)\mathrm{d}\Gamma,\\
		\mathcal{B}_v^{\Omega_\square^c\cap\Omega}[\Displacement,\displacement]\coloneqq& \int_{\Omega_\square^c\cap\partial\Omega_v}r_{i}(\Displacement) r_{i}(\displacement)\mathrm{d}\Gamma,\\
		\mathcal{B}_{C_u}^{\Omega_\square^c\cap\Omega}[\Displacement,\displacement]\coloneqq& \int_{\Omega_\square^c\cap\partial\Omega_{C_u}}j_{i}(\Displacement)j_{i}(\displacement)\mathrm{d}s,\\
		\mathcal{B}_\Omega^{\Omega_\square^c\cap\Omega}[\Displacement,\displacement]\coloneqq& \int_{\Omega_\square^c\cap\Omega}\left(\cauchyStress_{ij}(\Displacement)\strain_{ij}(\displacement)+\hyperStress_{ijk}(\Displacement)\strain_{ij,k}(\displacement)\right)\mathrm{d}\Omega.
\end{align}\end{subequations}

The numerical solution of the generalized eigenvalue problems requires careful consideration \cite{dePrenter2016}, since \emph{a)} the matrix $\mathbf{V}$ is always singular (it is based solely on derivative quantities of the test functions) and \emph{b)} the matrix $\mathbf{V}$ can be bad-conditioned if the corresponding volume fraction $\chi_c$ of the cell is very small. We refer to \cite{dePrenter2016} for further details.

\bibliography{References}

\begin{thebibliography}{10}

\bibitem{Mashkevich1957}
V.~S. Mashkevich and K.~B. Tolpygo.
\newblock Electrical, optical and elastic properties of diamond type cristals.
  1.
\newblock {\em Soviet Physics}, 5(3):435--439, 1957.

\bibitem{Tolpygo1963}
KB~Tolpygo.
\newblock Long wavelength oscillations of diamond-type crystals including long
  range forces.
\newblock {\em Soviet Physics-Solid State}, 4(7):1297--1305, 1963.

\bibitem{Kogan1964}
S.M. Kogan.
\newblock Piezoelectric effect during inhomogeneous deformation and acoustic
  scattering of carriers in crystals.
\newblock {\em Sov. Phys. Solid State}, 5(10), 1964.
\newblock cited By 216.

\bibitem{bursian1968}
J~M Bursian and O~I Zaikovskii.
\newblock {Changes in curvature of a ferroelectric film due to polarization}.
\newblock {\em Soviet Physics Solid State}, 10(5):1121--1124, 1968.

\bibitem{Indenbom1981}
VL~Indenbom, EB~Loginov, and MA~Osipov.
\newblock Flexoelectric effect and crystal-structure.
\newblock {\em Kristallografiya}, 26(6):1157--1162, 1981.

\bibitem{Tagantsev1986}
A.~K. Tagantsev.
\newblock Piezoelectricity and flexoelectricity in crystalline dielectrics.
\newblock {\em Phys. Rev. B}, 34:5883--5889, Oct 1986.

\bibitem{Tagantsev1991}
Alexander~K Tagantsev.
\newblock Electric polarization in crystals and its response to thermal and
  elastic perturbations.
\newblock {\em Phase Transitions: A Multinational Journal}, 35(3-4):119--203,
  1991.

\bibitem{Nguyen2013}
Thanh~D. Nguyen, Sheng Mao, Yao-Wen Yeh, Prashant~K. Purohit, and Michael~C.
  McAlpine.
\newblock Nanoscale flexoelectricity.
\newblock {\em Advanced Materials}, 25(7):946--974, 2013.

\bibitem{Maranganti2006}
R~Maranganti, ND~Sharma, and P~Sharma.
\newblock Electromechanical coupling in nonpiezoelectric materials due to
  nanoscale nonlocal size effects: Green’s function solutions and embedded
  inclusions.
\newblock {\em Physical Review B}, 74(1):014110, 2006.

\bibitem{Majdoub2008}
MS~Majdoub, P~Sharma, and T~\ifmmode \mbox{\c{C}}\else
  \c{C}\fi{}a\ifmmode~\breve{g}\else \u{g}\fi{}in.
\newblock Enhanced size-dependent piezoelectricity and elasticity in
  nanostructures due to the flexoelectric effect.
\newblock {\em Physical Review B}, 77(12):125424, 2008.

\bibitem{Sharma2010}
N.D. Sharma, C.M. Landis, and P.~Sharma.
\newblock Piezoelectric thin-film superlattices without using piezoelectric
  materials.
\newblock {\em Journal of Applied Physics}, 108(2):1--25, 2010.

\bibitem{Cross2006}
L.~Eric Cross.
\newblock Flexoelectric effects: Charge separation in insulating solids
  subjected to elastic strain gradients.
\newblock {\em Journal of Materials Science}, 41(1):53--63, Jan 2006.

\bibitem{Catalan2004}
G~Catalan, LJ~Sinnamon, and JM~Gregg.
\newblock The effect of flexoelectricity on the dielectric properties of
  inhomogeneously strained ferroelectric thin films.
\newblock {\em Journal of Physics: Condensed Matter}, 16(13):2253, 2004.

\bibitem{Eliseev2009}
Eugene~A Eliseev, Anna~N Morozovska, Maya~D Glinchuk, and R~Blinc.
\newblock Spontaneous flexoelectric/flexomagnetic effect in nanoferroics.
\newblock {\em Physical Review B}, 79(16):165433, 2009.

\bibitem{Liu2014}
Liping Liu.
\newblock An energy formulation of continuum magneto-electro-elasticity with
  applications.
\newblock {\em Journal of the Mechanics and Physics of Solids}, 63:451 -- 480,
  2014.

\bibitem{Shen2010}
Shengping Shen and Shuling Hu.
\newblock A theory of flexoelectricity with surface effect for elastic
  dielectrics.
\newblock {\em Journal of the Mechanics and Physics of Solids}, 58(5):665 --
  677, 2010.

\bibitem{Hu2010}
ShuLing Hu and ShengPing Shen.
\newblock Variational principles and governing equations in nano-dielectrics
  with the flexoelectric effect.
\newblock {\em Science China Physics, Mechanics and Astronomy},
  53(8):1497--1504, Aug 2010.

\bibitem{Yudin2013}
P~V Yudin and A~K Tagantsev.
\newblock Fundamentals of flexoelectricity in solids.
\newblock {\em Nanotechnology}, 24(43):432001, 2013.

\bibitem{Zubko2013}
Pavlo Zubko, Gustau Catalan, and Alexander~K. Tagantsev.
\newblock Flexoelectric effect in solids.
\newblock {\em Annual Review of Materials Research}, 24(43):387--421, 2013.

\bibitem{krichen2016}
Sana Krichen and Pradeep Sharma.
\newblock Flexoelectricity: A perspective on an unusual electromechanical
  coupling.
\newblock {\em Journal of Applied Mechanics}, 83(3):030801, 2016.

\bibitem{Shu2011}
Longlong Shu, Xiaoyong Wei, Ting Pang, Xi~Yao, and Chunlei Wang.
\newblock Symmetry of flexoelectric coefficients in crystalline medium.
\newblock {\em Journal of Applied Physics}, 110(10):104106, 2011.

\bibitem{LeQuang2011}
H.~Le~Quang and Q.-C. He.
\newblock The number and types of all possible rotational symmetries for
  flexoelectric tensors.
\newblock {\em Proceedings of the Royal Society of London A: Mathematical,
  Physical and Engineering Sciences}, 467(2132):2369--2386, 2011.

\bibitem{Abdollahi2014}
Amir Abdollahi, Christian Peco, Daniel Mill\'an, Marino Arroyo, and Irene
  Arias.
\newblock Computational evaluation of the flexoelectric effect in dielectric
  solids.
\newblock {\em Journal of Applied Physics}, 116(9):093502, 2014.

\bibitem{Abdollahi2015a}
Amir Abdollahi, Daniel Mill\'an, Christian Peco, Marino Arroyo, and Irene
  Arias.
\newblock Revisiting pyramid compression to quantify flexoelectricity: A
  three-dimensional simulation study.
\newblock {\em Phys. Rev. B}, 91:104103, Mar 2015.

\bibitem{Abdollahi2015b}
Amir Abdollahi, Christian Peco, Daniel Mill\'an, Marino Arroyo, Gustau Catalan,
  and Irene Arias.
\newblock Fracture toughening and toughness asymmetry induced by
  flexoelectricity.
\newblock {\em Phys. Rev. B}, 92:094101, Sep 2015.

\bibitem{Ghasemi2017}
Hamid Ghasemi, Harold~S. Park, and Timon Rabczuk.
\newblock A level-set based iga formulation for topology optimization of
  flexoelectric materials.
\newblock {\em Computer Methods in Applied Mechanics and Engineering}, 313:239
  -- 258, 2017.

\bibitem{Nanthakumar2017}
S.S. Nanthakumar, Xiaoying Zhuang, Harold~S. Park, and Timon Rabczuk.
\newblock Topology optimization of flexoelectric structures.
\newblock {\em Journal of the Mechanics and Physics of Solids}, 105:217 -- 234,
  2017.

\bibitem{Ghasemi2018}
Hamid Ghasemi, Harold~S. Park, and Timon Rabczuk.
\newblock A multi-material level set-based topology optimization of
  flexoelectric composites.
\newblock {\em Computer Methods in Applied Mechanics and Engineering}, 332:47
  -- 62, 2018.

\bibitem{Yvonnet2017}
Julien Yvonnet and LP~Liu.
\newblock A numerical framework for modeling flexoelectricity and maxwell
  stress in soft dielectrics at finite strains.
\newblock {\em Computer Methods in Applied Mechanics and Engineering},
  313:450--482, 2017.

\bibitem{Mao2016}
Sheng Mao, Prashant~K. Purohit, and Nikolaos Aravas.
\newblock Mixed finite-element formulations in piezoelectricity and
  flexoelectricity.
\newblock {\em Proceedings of the Royal Society of London A: Mathematical,
  Physical and Engineering Sciences}, 472(2190), 2016.

\bibitem{Deng2017}
Feng Deng, Qian Deng, Wenshan Yu, and Shengping Shen.
\newblock Mixed finite elements for flexoelectric solids.
\newblock {\em Journal of Applied Mechanics}, 84(8):081004, 2017.

\bibitem{Abdollahi2015c}
Amir Abdollahi and Irene Arias.
\newblock Constructive and destructive interplay between piezoelectricity and
  flexoelectricity in flexural sensors and actuators.
\newblock {\em Journal of Applied Mechanics}, 82(12):121003, 2015.

\bibitem{Aravas2011}
Nikolaos Aravas.
\newblock Plane-strain problems for a class of gradient elasticity models - a
  stress function approach.
\newblock {\em Journal of Elasticity}, 104(1-2):45--70, 2011.

\bibitem{CodonyNonlinearFlexoUnpublished}
David Codony, Prakhar Gupta, Onofre Marco, and Irene Arias.
\newblock Modeling flexoelectricity in soft dielectrics at finite deformation,
  2020.

\bibitem{Sevilla2008}
Ruben Sevilla, Sonia Fern{\'a}ndez-M{\'e}ndez, and Antonio Huerta.
\newblock Nurbs-enhanced finite element method (nefem).
\newblock {\em International Journal for Numerical Methods in Engineering},
  76(1):56--83, 2008.

\bibitem{Sukumar2001}
N.~Sukumar, D.L. Chopp, Nicolas Mo{\"e}s, and Ted Belytschko.
\newblock Modeling holes and inclusions by level sets in the extended
  finite-element method.
\newblock {\em Computer Methods in Applied Mechanics and Engineering}, 190:6183
  -- 6200, 2001.

\bibitem{Belytschko2003}
Ted Belytschko, Chandu Parimi, Nicolas Mo{\"e}s, N.~Sukumar, and Shuji Usui.
\newblock Structured extended finite element methods for solids defined by
  implicit surfaces.
\newblock {\em International journal for Numerical Methods in Engineering},
  56:609 -- 635, 2003.

\bibitem{Bandara2016}
Kosala Bandara, Thomas R{\"u}berg, and Fehmi Cirak.
\newblock Shape optimisation with multiresolution subdivision surfaces and
  immersed finite elements.
\newblock {\em Computer Methods in Applied Mechanics and Engineering},
  300:510--539, 2016.

\bibitem{Scott2011}
M~A Scott, X~Li, T~W Sederberg, and T~J~R Hughes.
\newblock Local refinement of analysis-suitable t-splines, ices report 11-06.
\newblock {\em The Institute for Computational Engineering and Sciences, The
  University of Texas at Austin}, 2011.

\bibitem{Forsey1988}
David~R. Forsey and Richard~H. Bartels.
\newblock Hierarchical b-spline refinement.
\newblock {\em ACM Siggraph Computer Graphics}, 22(4):205--212, 1988.

\bibitem{Kraft1995}
Rainer Kraft.
\newblock Hierarchical b-splines.
\newblock 1995.

\bibitem{Kraft1997}
R.~Kraft.
\newblock {\em Adaptive and Linearly Independent Multilevel B-splines}.
\newblock Bericht. SFB 404, Gesch{\"a}ftsstelle, 1997.

\bibitem{Mindlin1964}
Raymond~David Mindlin.
\newblock Micro-structure in linear elasticity.
\newblock {\em Archive for Rational Mechanics and Analysis}, 16(1):51--78,
  1964.

\bibitem{Mindlin1968a}
R.~D. Mindlin and N.~N. Eshel.
\newblock On first strain-gradient theories in linear elasticity.
\newblock {\em International Journal of Solids and Structures}, 4(1):109--124,
  1968.

\bibitem{Mindlin1968b}
Raymond~David Mindlin.
\newblock Polarization gradient in elastic dielectrics.
\newblock {\em International Journal of Solids and Structures}, 4(6):637--642,
  1968.

\bibitem{Mao2014}
Sheng Mao and Prashant~K. Purohit.
\newblock Insights into flexoelectric solids from strain-gradient elasticity.
\newblock {\em ASME Journal of Applied Mechanics}, 81(8):1--10, 2014.

\bibitem{Block1997}
Ethan~D Bloch.
\newblock {\em A first course in geometric topology and differential geometry}.
\newblock Springer Science \& Business Media, 1997.

\bibitem{Majdoub2009}
M.~S. Majdoub, P.~Sharma, and T.~\ifmmode \mbox{\c{C}}\else
  \c{C}\fi{}a\ifmmode~\breve{g}\else \u{g}\fi{}in.
\newblock Erratum: Enhanced size-dependent piezoelectricity and elasticity in
  nanostructures due to the flexoelectric effect [phys. rev. b 77, 125424
  (2008)].
\newblock {\em Phys. Rev. B}, 79:119904, Mar 2009.

\bibitem{Landau2013}
Lev~Davidovich Landau and Evgenii~Mikhailovich Lifshitz.
\newblock {\em Course of theoretical physics}.
\newblock Elsevier, 2013.

\bibitem{Nitsche1971}
J.~Nitsche.
\newblock {\"U}ber ein {V}ariationsprinzip zur {L}{\"o}sung von
  {D}irichlet-{P}roblemen bei {V}erwendung von {T}eilr{\"a}umen, die keinen
  {R}andbedingungen unterworfen sind.
\newblock {\em Abhandlungen aus dem Mathematischen Seminar der Universit{\"a}t
  Hamburg}, 36(1):9--15, 1971.

\bibitem{Fernandez2004}
Sonia Fern{\'a}ndez-M{\'e}ndez and Antonio Huerta.
\newblock Imposing essential boundary conditions in mesh-free methods.
\newblock {\em Computer methods in applied mechanics and engineering},
  193(12):1257--1275, 2004.

\bibitem{Schillinger2016}
Dominik Schillinger, Isaac Harari, Ming-Chen Hsu, David Kamensky, Stein~KF
  Stoter, Yue Yu, and Ying Zhao.
\newblock The non-symmetric nitsche method for the parameter-free imposition of
  weak boundary and coupling conditions in immersed finite elements.
\newblock {\em Computer Methods in Applied Mechanics and Engineering},
  309:625--652, 2016.

\bibitem{badia2018}
Santiago Badia, Francesc Verdugo, and Alberto~F Mart{\'\i}n.
\newblock The aggregated unfitted finite element method for elliptic problems.
\newblock {\em Computer Methods in Applied Mechanics and Engineering},
  336:533--553, 2018.

\bibitem{Ruberg2016}
Thomas R{\"u}berg, Fehmi Cirak, and Jos{\'e}~Manuel Garc{\'i}a~Aznar.
\newblock An unstructured immersed finite element method for nonlinear solid
  mechanics.
\newblock {\em Advanced Modeling and Simulation in Engineering Sciences},
  3(1):22, 2016.

\bibitem{deBoor2001}
C.~de~Boor.
\newblock {\em A Practical Guide to Splines}.
\newblock Applied Mathematical Sciences. Springer New York, 2001.

\bibitem{Rogers2001}
D.F. Rogers.
\newblock {\em An Introduction to NURBS: With Historical Perspective}.
\newblock Morgan Kaufmann Series in Computer Graphics and Geometric Modeling.
  Morgan Kaufmann Publishers, 2001.

\bibitem{Piegl2012}
L.~Piegl and W.~Tiller.
\newblock {\em The NURBS Book}.
\newblock Monographs in Visual Communication. Springer Berlin Heidelberg, 2012.

\bibitem{Hughes2005}
Thomas J~R Hughes, John~A Cottrell, and Yuri Bazilevs.
\newblock Isogeometric analysis: Cad, finite elements, nurbs, exact geometry
  and mesh refinement.
\newblock {\em Computer methods in applied mechanics and engineering},
  194(39):4135--4195, 2005.

\bibitem{Peskin2002}
Charles~S. Peskin.
\newblock The immersed boundary method.
\newblock {\em Acta Numerica}, 11:479--517, 2002.

\bibitem{Mittal2005}
Rajat Mittal and Gianluca Iaccarino.
\newblock Immersed boundary methods.
\newblock {\em Annu. Rev. Fluid Mech.}, 37:239--261, 2005.

\bibitem{Fries2015}
Thomas-Peter Fries and Samir Omerovi{\'c}.
\newblock Higher-order accurate integration of implicit geometries.
\newblock {\em International Journal for Numerical Methods in Engineering},
  106(5):323--371, 2016.
\newblock nme.5121.

\bibitem{Fries2016}
Thomas-Peter Fries.
\newblock {\em Higher-Order Accurate Integration for Cut Elements with
  Chen-Babu{\v{s}}ka Nodes}, pages 245--269.
\newblock Springer International Publishing, 2016.

\bibitem{Kudela2016}
L{\'a}szl{\'o} Kudela, Nils Zander, Stefan Kollmannsberger, and Ernst Rank.
\newblock Smart octrees: Accurately integrating discontinuous functions in 3d.
\newblock {\em Computer Methods in Applied Mechanics and Engineering},
  306:406--426, 2016.

\bibitem{Legrain2012}
G.~Legrain, N.~Chevaugeon, and K.~Dr\'{e}au.
\newblock High order x-fem and levelsets for complex microstructures:
  Uncoupling geometry and approximation.
\newblock {\em Computer Methods in Applied Mechanics and Engineering}, 241:172
  -- 189, 2012.

\bibitem{Marco2015}
Onofre Marco, Ruben Sevilla, Yongjie Zhang, Juan~Jos{\'e} R{\'o}denas, and
  Manuel Tur.
\newblock Exact 3d boundary representation in finite element analysis based on
  cartesian grids independent of the geometry.
\newblock {\em International Journal for Numerical Methods in Engineering},
  103(6):445--468, 2015.

\bibitem{Marco2017}
Onofre Marco, Juan~Jos{\'e} R{\'o}denas, Jos{\'e}~Manuel Navarro-Jim{\'e}nez,
  and Manuel Tur.
\newblock Robust h-adaptive meshing strategy considering exact arbitrary cad
  geometries in a cartesian grid framework.
\newblock {\em Computers \& Structures}, 193:87--109, 2017.

\bibitem{Vincent2015}
F.D. Witherden and P.E. Vincent.
\newblock On the identification of symmetric quadrature rules for finite
  element methods.
\newblock {\em Computers \& Mathematics with Applications}, 69(10):1232 --
  1241, 2015.

\bibitem{Lorensen1987}
William~E. Lorensen and Harvey~E. Cline.
\newblock Marching cubes: A high resolution 3d surface construction algorithm.
\newblock {\em Computer Graphics}, 21(4):163--169, 1987.

\bibitem{Duster2008}
A.~D\"{u}ster, J.~Parvizian, Z.~Yang, and E.~Rank.
\newblock The finite cell method for three-dimensional problems of solid
  mechanics.
\newblock {\em Computer Methods in Applied Mechanics and Engineering},
  197(45):3768 -- 3782, 2008.

\bibitem{Schillinger2015}
Dominik Schillinger and Martin Ruess.
\newblock The finite cell method: A review in the context of higher-order
  structural analysis of cad and image-based geometric models.
\newblock {\em Archives of Computational Methods in Engineering},
  22(3):391--455, Jul 2015.

\bibitem{Sevilla2011Seamless}
Ruben Sevilla, Sonia Fern{\'a}ndez-M{\'e}ndez, and Antonio Huerta.
\newblock Nurbs-enhanced finite element method (nefem).
\newblock {\em Archives of Computational Methods in Engineering}, 18(4):441,
  2011.

\bibitem{Sevilla2011Numerical}
Ruben Sevilla and Sonia Fern{\'a}ndez-M{\'e}ndez.
\newblock Numerical integration over 2d nurbs-shaped domains with applications
  to nurbs-enhanced fem.
\newblock {\em Finite Elements in Analysis and Design}, 47(10):1209--1220,
  2011.

\bibitem{Sevilla20113D}
Ruben Sevilla, Sonia Fern{\'a}ndez-M{\'e}ndez, and Antonio Huerta.
\newblock 3d nurbs-enhanced finite element method (nefem).
\newblock {\em International Journal for Numerical Methods in Engineering},
  88(2):103--125, 2011.

\bibitem{Legrain2013}
Gr{\'e}gory Legrain.
\newblock A nurbs enhanced extended finite element approach for unfitted cad
  analysis.
\newblock {\em Computational Mechanics}, 52(4):913--929, 2013.

\bibitem{dePrenter2016}
F~de~Prenter, CV~Verhoosel, GJ~van Zwieten, and EH~van Brummelen.
\newblock Condition number analysis and preconditioning of the finite cell
  method.
\newblock {\em Computer Methods in Applied Mechanics and Engineering}, 2016.

\bibitem{Burman2010b}
Erik Burman.
\newblock Ghost penalty.
\newblock {\em Comptes Rendus Mathematique}, 348(21-22):1217--1220, 2010.

\bibitem{Hollig2001}
Klaus H{\"o}llig, Ulrich Reif, and Joachim Wipper.
\newblock Weighted extended b-spline approximation of dirichlet problems.
\newblock {\em SIAM Journal on Numerical Analysis}, 39(2):442--462, 2001.

\bibitem{Hollig2012}
Klaus H{\"o}llig, J{\"o}rg H{\"o}rner, and Axel Hoffacker.
\newblock Finite element analysis with b-splines: weighted and isogeometric
  methods.
\newblock In {\em International Conference on Curves and Surfaces}, pages
  330--350. Springer, 2012.

\bibitem{Ruberg2012}
T~R{\"u}berg and F~Cirak.
\newblock Subdivision-stabilised immersed b-spline finite elements for moving
  boundary flows.
\newblock {\em Computer Methods in Applied Mechanics and Engineering},
  209:266--283, 2012.

\bibitem{Zorin2000}
Denis Zorin.
\newblock Subdivision for modeling and animation.
\newblock {\em SIGGRAPH 2000 Course Notes}, pages 65--104, 2000.

\bibitem{Schillinger2012}
Dominik Schillinger, Luca Ded\`{e}, Michael~A. Scott, John~A. Evans, Michael~J.
  Borden, Ernst Rank, and Thomas~J.R. Hughes.
\newblock An isogeometric design-through-analysis methodology based on adaptive
  hierarchical refinement of nurbs, immersed boundary methods, and t-spline cad
  surfaces.
\newblock {\em Computer Methods in Applied Mechanics and Engineering},
  249--252:116--150, 2012.

\bibitem{Bornemann2013}
P.B. Bornemann and F.~Cirak.
\newblock A subdivision-based implementation of the hierarchical b-spline
  finite element method.
\newblock {\em Computer Methods in Applied Mechanics and Engineering},
  253:584--598, 2013.

\bibitem{Vuong2011}
A.-V. Vuong, Carlotta Giannelli, Bert J{\"u}ttler, and Bernd Simeon.
\newblock A hierarchical approach to adaptive local refinement in isogeometric
  analysis.
\newblock {\em Computer Methods in Applied Mechanics and Engineering},
  200(49):3554--3567, 2011.

\bibitem{MocciShearUnpublished}
Alice Mocci, Amir Abdollahi, and Irene Arias.
\newblock Quantification of shear flexoelectricity in ferroelectrics.
\newblock {\em In preparation}.

\bibitem{Baskaran2012}
Sivapalan Baskaran, Xiangtong He, Yu~Wang, and John~Y Fu.
\newblock Strain gradient induced electric polarization in $\alpha$-phase
  polyvinylidene fluoride films under bending conditions.
\newblock {\em Journal of Applied Physics}, 111(1):014109, 2012.

\bibitem{Chu2012}
Baojin Chu and D~R Salem.
\newblock Flexoelectricity in several thermoplastic and thermosetting polymers.
\newblock {\em Applied Physics Letters}, 101(10):103905, 2012.

\bibitem{Ma2002}
Wenhui Ma and L~Eric Cross.
\newblock Flexoelectric polarization of barium strontium titanate in the
  paraelectric state.
\newblock {\em Applied Physics Letters}, 81(18):3440--3442, 2002.

\bibitem{Ma2005}
Wenhui Ma and L~Eric Cross.
\newblock Flexoelectric effect in ceramic lead zirconate titanate.
\newblock {\em Applied Physics Letters}, 86(7):072905, 2005.

\bibitem{Ma2006}
Wenhui Ma and L~Eric Cross.
\newblock Flexoelectricity of barium titanate.
\newblock {\em Applied Physics Letters}, 88(23):232902, 2006.

\bibitem{Lu2016}
Jianfeng Lu, Jiangyan Lv, Xu~Liang, Minglong Xu, and Shengping Shen.
\newblock Improved approach to measure the direct flexoelectric coefficient of
  bulk polyvinylidene fluoride.
\newblock {\em Journal of Applied Physics}, 119(9):094104, 2016.

\bibitem{Baskaran2011a}
Sivapalan Baskaran, Narayanan Ramachandran, Xiangtong He, Sankar
  Thiruvannamalai, Ho~Joon Lee, Hyun Heo, Qin Chen, and John~Y Fu.
\newblock Giant flexoelectricity in polyvinylidene fluoride films.
\newblock {\em Physics Letters A}, 375(20):2082--2084, 2011.

\bibitem{Baskaran2011b}
Sivapalan Baskaran, Xiangtong He, Qin Chen, and John~Y Fu.
\newblock Experimental studies on the direct flexoelectric effect in
  $\alpha$-phase polyvinylidene fluoride films.
\newblock {\em Applied Physics Letters}, 98(24):242901, 2011.

\bibitem{Amanatidou2002}
E~Amanatidou and N~Aravas.
\newblock Mixed finite element formulations of strain-gradient elasticity
  problems.
\newblock {\em Computer Methods in Applied Mechanics and Engineering},
  191(15):1723--1751, 2002.

\bibitem{Altan1997}
B~S Altan and E~C Aifantis.
\newblock On some aspects in the special theory of gradient elasticity.
\newblock {\em Journal of the Mechanical Behavior of Materials}, 8(3):231--282,
  1997.

\bibitem{Embar2010}
Anand Embar, John Dolbow, and Isaac Harari.
\newblock Imposing dirichlet boundary conditions with nitsche's method and
  spline-based finite elements.
\newblock {\em International Journal for Numerical Methods in Engineering},
  83(7):877--898, 2010.

\bibitem{Griebel2003}
Michael Griebel and Marc~Alexander Schweitzer.
\newblock A particle-partition of unity method part v: boundary conditions.
\newblock In {\em Geometric analysis and nonlinear partial differential
  equations}, pages 519--542. Springer, 2003.

\end{thebibliography}
\end{document}